\documentclass[oneside]{amsart}
\usepackage{amsmath,amssymb,color,amsthm,graphicx,cite}
\usepackage{color,bm}
\graphicspath{{../../img_pdf/}{../img_pdf/}{./img_pdf/}}
\usepackage{enumerate} 
\usepackage[backref]{hyperref}
\newtheorem{theorem}{Theorem}[section]
\newtheorem{proposition}[theorem]{Proposition}
\newtheorem{lemma}[theorem]{Lemma}
\newtheorem{corollary}[theorem]{Corollary}
\newtheorem{claim}{Claim}

\newtheorem{question}{Question}

\theoremstyle{definition}
\newtheorem{definition}{Definition}
\newtheorem{main}{Theorem}

\def\Z{\mathbb{Z} }
\def\R{\mathbb{R} }

\def\nbd{neighborhood }
\def\nbds{neighborhoods }

\author{Tomoo Yokoyama}
\date{\today}
\address{Department of Mathematics, Faculty of Science, Saitama University, Shimo-Okubo 255, Sakura-ku, Saitama-shi, 338-8570 Japan\\}
\email{tyokoyama@rimath.saitama-u.ac.jp}
\thanks{The author was partially supported by JSPS Grant Number 24K06733}

\title[Structures of the space of gradient vector fields]{Combinatorial structures of the space of gradient vector fields on compact surfaces}
\makeatletter
\@namedef{subjclassname@2020}{%
  \textup{2020} Mathematics Subject Classification}
\makeatother
\subjclass[2020]{Primary 37G10; Secondary 37E35,57Q10,\\58B05,76A02}

\keywords{Homotopy type; beat point; gradient vector field; bifurcation}

\begin{document}
\maketitle

\begin{abstract}
Gradient vector fields are fundamental objects from both theoretical and practical perspectives, since various phenomena can be modeled within this framework. The ``moduli space'' of such vector fields provides the foundation for describing these phenomena. However, little is known about the topology of the space of gradient vector fields. For instance, it remains unknown whether a connected component of this space can fail to be simply connected. This paper aims to lay the foundation for describing the possible generic time evolution of gradient vector fields on surfaces, with or without constraints, under the assumption that no creation or annihilation of singular points occurs, by using combinatorics and simple homotopy theory. In fact, the space of gradient vector fields on a closed annulus contains a non-contractible connected component, which is weakly homotopy equivalent to a bouquet of two two-dimensional spheres.
\end{abstract}

\section{Introduction}\label{intro}

\begin{figure}[t]
\begin{center}
\includegraphics[scale=0.45]{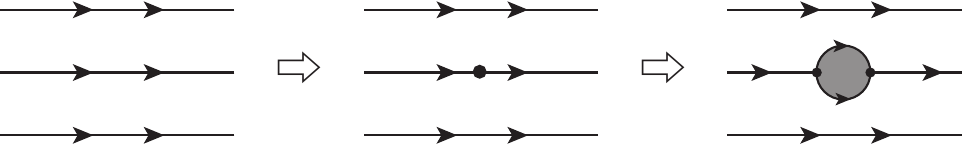}
\end{center}
\caption{Creation of a physical boundary.}
\label{fig:creations_bdry}
\end{figure}

A gradient vector field is one of the fundamental objects from a theoretical and practical point of view. 
In the time evolution of fluids on punctured spheres, some kinds of such fluids are modeled by gradient vector fields, and the creations and annihilations of singular points and physical boundaries can change the topologies of streamlines.
For instance, one can observe the creation of a physical boundary, which is a boundary of a stone on the surface of a river, when the river's water level goes down, as in Figure~\ref{fig:creations_bdry}. 
Notice that creations and annihilations of physical boundaries change the topologies of surfaces. 
On the other hand, the topologies of such fluids also can be changed by switching combinatorial structures of separatrices. 
Such combinatorial structures are studied from fluid mechanics \cite{aref1998stagnation,kidambi2000streamline,moffatt2001topology}, integrable systems \cite{bolsinov2004integrable}, and dynamical systems \cite{ma2005geometric,sakajo2014unique,sakajo2015transitions,sakajo2018tree,sakajo2020discrete,yokoyama2021complete,yokoyama2013word,yokoyama2021combinatorial,yokoyama2021cot}. 

From a dynamical system's point of view, Smale \cite{smale1961gradient} proved that any Morse vector field (i.e. Morse-Smale vector field without limit cycles) on a closed manifold is a gradient vector field without separatrices from a saddle to a saddle. 
By a work of Andronov-Pontryagin~\cite{andronov1937rough} and a work of Peixoto~\cite{peixoto1962structural}, it is known that the set of Morse-Smale $C^r$-vector fields ($r \geq 1$) on a closed orientable surface $S$ is open dense in the space of $C^r$-vector fields on $S$ and that Morse-Smale $C^r$-vector fields are structurally stable in the space of $C^r$-vector fields.
In particular, the set of Morse $C^r$-vector fields (i.e. Morse-Smale vector fields without limit cycles) on a closed orientable surface is open dense in the space of gradient vector fields. 
Similarly, in \cite{labarca1990stability}, Morse-Smale $C^r$-vector fields ($r \geq 1$) on a compact surface possibly with boundary is defined, and the structural stability is shown for Morse-Smale vector fields on orientable compact manifolds under the $C^2$ linearizable condition and the eigenvalue conditions. 
By these facts, one characterizes a ``generic'' non-Morse gradient vector field on an orientable compact surface to describe a generic time evaluation of gradient vector fields on compact surfaces (e.g. solutions of differential equations) which is an alternating sequence of Morse vector fields and instantaneous non-Morse gradient vector fields under no physical or symmetric restrictions \cite{kibkalo2022topological}.

On the other hand, ``non-generic'' intermediate vector fields (e.g. fluids with symmetric vortex pairs) naturally appear under physical or symmetric restrictions. 
For instance, degenerate multi-saddle connections are the reasons for ``non-genericity''. 
Though the hierarchical structure of the space of gradient vector fields is one of the foundations for describing generic time evaluations, only the low codimensional structures were studied. 
More globally, we ask the following question. 

\begin{question}\label{q:01}
Does the space of topological equivalence classes of gradient vector fields on a manifold have non-contractible connected components, under the non-existence of creations and annihilations of singular points? 
\end{question}

In this paper, we demonstrate that there is such a non-contractible connected component of the space of topological equivalence classes of gradient vector fields on a manifold as follows.

\begin{main}\label{re_prop:5.3}
For any $r \in \Z_{\geq 0} \sqcup \{ \infty \}$, the space of $C^r$ gradient vector fields with finitely many singular points without fake multi-saddles \rm{(i.e.} $\partial$-$0$-saddles and $0$-saddles\rm{)} or fake parabolic sectors (see Definition~\ref{def:fake} below) on a closed annulus whose sums of indices of sinks and $\partial$-sinks (resp. sources and $\partial$-sources) are $1$ (resp. $2$) has a connected component which is weakly homotopy equivalent to a bouquet of two two-dimensional spheres. 
\end{main}

This means that the structure of ``transitions of gradient vector fields on compact surfaces'' is not contractible. 
To calculate the non-contractible property, we study the hierarchical structure of the space of such vector fields under the non-existence of creations and annihilations of singular points. 
In fact, the space of topologically equivalence classes of such vector fields is a disjoint union of finite abstract cell complexes such that the codimension of a cell corresponds to the instability of a gradient vector field as follows. 

\begin{main}\label{main:02-}
The following statements hold for any $r \in \mathbb{Z}_{> 0} \sqcup \{ \infty \}$, any integers $k_-, k_+ \in \mathbb{Z}_{\geq 0}$ and any $q \in \Z_{\geq -1}$:
\\
{\rm(1)} The space $P$ of topologically equivalence classes of $C^r$ gradient vector fields with at most $l$ singular points but without fake multi-saddles or fake parabolic sectors on a compact surface is a finite abstract cell complex, a finite poset, and a finite $T_0$-space. 
\\
{\rm(2)} For any $q \in \Z_{\geq 0}$, the closure of the set of elements of height $q$ in the poset $P$ is the set  of elements of height at most $q$ in $P$ \rm{(i.e.} $\overline{P_{=q}} = \bigsqcup_{j=0}^q P_{=j}$, where $P_{=s}$ is the set of elements of height $s$ in the poset $X$ \rm{)}. 
\end{main}

Using the hierarchical structure, we can describe the possible generic time evolution of gradient vector fields on a compact surface with or without restriction conditions.

The present paper consists of seven sections and Appendix.
In the next section, as preliminaries, we introduce fundamental concepts.
In \S 3, the invariance of sums of indices of isolated singular points under small perturbations and the existence of ``good'' \nbds of the multi-saddle connection diagram are described (Proposition~\ref{lem:no_merge_msc}). 
In \S 4, we study the combinatorial structure of the space of gradient vector fields under the non-existence of creations and annihilations of singular points and boundary components (Theorem~\ref{prop:codimension}). 
In \S 5, the abstract cell complex structure and filtration of the space of gradient vector fields are described (Theorem~\ref{main:02} and Theorem~\ref{main:02-}). 
In \S 6, we demonstrate the non-contractibility of a connected component of the space of gradient vector fields (Theorem~\ref{lem:007} and Theorem~\ref{re_prop:5.3}). 
In \S 7, we describe the combinatorial structure of the space of Morse-Smale-like flows (Theorem~\ref{prop:codimension_ms} and Theorem~\ref{main:03}). 
The final section states a future work and an open question. 

\section{Preliminaries}

\subsection{Notion of topology}

A subset $A$ of a topological space is a {\bf collar} of a subset $B \subseteq A$ if there is a homeomorphism $h \colon [0,1] \times B \to A$ with $h(\{ 0\} \times B) = B$. 

A point $x$ of a topological space $X$ is $\bm{T_0}$ (or Kolmogorov) if for any point $y \neq x \in X$ there is an open subset $U$ of $X$ such that $\vert \{x, y \} \cap U \vert  =1$, where $\vert A \vert$ is the cardinality of a subset $A$.
A topological space is $T_0$ if each point is $T_0$.

The {\bf specialization preorder} $\leq_{\tau}$ on a topological space $(X, \tau)$ is defined as follows: $ x \leq_{\tau} y $ if  $ x \in \overline{\{ y \}}$, where $\overline{\{ y \}}$ is the closure of $\{ y \}$. 
A topology $\tau$ is $T_0$ if and only if $\leq_{\tau}$ is a partial order. 
The heights of a topological space, a point, and a subset are defined as the heights with respect to the specialization preorder. 
For any $k \in \Z_{\geq 0}$ and a topological space $X$, denote by $X_k$ the set of height $k$ points of $X$ and by $X_{\leq k}$ the set of points of $X$ whose height is less than or equal to $k$. 

\subsubsection{Surface graph}
By a {\bf surface}, we mean a two-dimensional paracompact manifold, that does not need to be orientable.
A {\bf graph} is a cell complex whose dimension is at most one and which is a geometric realization of an abstract multi-graph.
In other words, it can be drawn such that no edges cross each other.
Such a drawing is called a {\bf surface graph}.

\subsubsection{Double of a compact surface}

The {\bf double} $S_{\mathrm{dbl}}$ of a compact surface $S$ is the resulting closed surface by gluing two copies of $S$ as follows:  $S_{\mathrm{dbl}} := S \times \{+,-\} / \sim$ where $(x,+) \sim (x,-)$ if $x \in \partial S$.

\subsubsection{Curves and loops}
A {\bf curve} is a continuous mapping $C: I \to X$ where $I$ is a non-degenerate connected subset of a circle $\mathbb{S}^1$.
A curve is {\bf simple} if it is injective.
We also denote by $C$ the image of a curve $C$.
Denote by $\partial C := C(\partial I)$ the boundary of a curve $C$ if $C$ can be extended into a continuous map whose domain is $I \cup \partial I$, where $\partial I$ is the boundary of $I \subset \mathbb{S}^1$. 
Put $\mathop{\mathrm{int}} C := C \setminus \partial C$ if $\partial C$ is defined. 
A simple curve is a simple closed curve if its domain is $\mathbb{S}^1$ (i.e. $I = \mathbb{S}^1$).
A simple closed curve is also called a {\bf loop}. 
An {\bf arc} is a simple curve whose domain is an interval. 

\subsection{Notion of dynamical systems}
A {\bf flow} is a continuous $\R$-action on a topological space.
Let $v : \R \times X \to X$ be a flow on a topological space $X$.
For $t \in \R$, define $v^t : X \to X$ by $v^t := v(t, \cdot )$.
For a point $x$ of $X$, we denote by $O(x)$ the {\bf orbit} of $x$ (i.e. $O(x) := \{ v^t(x) \mid t \in \R \} = v(\R, x)$).
A positive (resp. negative) orbit of $x$ is $v(\R_{>0}, x)$ (resp. $v(\R_{<0}, x)$), denoted by $O^+(x)$ (resp. $O^-(x)$).
A point $x$ of $M$ is {\bf singular} if $x = v^t(x)$ for any $t \in \R$ and is {\bf periodic} if there is a positive number $T > 0$ such that $x = v^T(x)$ and $x \neq v^t(x)$ for any $t \in (0, T)$.
A point is {\bf closed} if it is singular or periodic. 
An orbit is singular (resp. periodic, closed) if it contains singular (resp. periodic, closed) points. 
Denote by $\bm{\mathop{\mathrm{Sing}}(v)}$ (resp. $\bm{\mathop{\mathrm{Per}}(v)}$, $\bm{\mathop{\mathrm{Cl}}(v)}$) the set of singular (resp. periodic, closed) points. 
Denote by $\overline{A}$ the closure of a subset $A$.

A subset of $X$ is said to be {\bf invariant} (or {\bf saturated}) if it is a union of orbits. 
The saturation $\mathrm{Sat}_v(A) = v(A)$ of a subset $A \subseteq X$ is the union of orbits intersecting $A$. 
A flow $v$ on a topological space $S$ is {\bf topologically equivalent} to a flow $w$ on a topological space $T$ if there is a homeomorphism $h \colon S \to T$ such that the image of any orbit of $v$ is an orbit of $w$ and that $h$ preserves the directions of orbits of $v$ and $w$.

\subsubsection{The $C^0$-topology on the set of flows on a topological space}

The $C^0$-topology (or compact-open topology) on the set of flows on a topological space $X$ is generated by a subbase $\{ \mathcal{U}(K,U) \mid K \subseteq \R \times X : \text{compact}, U \subseteq X : \text{open } \}$, where $\mathcal{U}(K, U) := \{ w : \text{flow on } X \mid w(K) \subset U \}$. 
By $C^r$-topology on the space of $C^r$-flows on a manifold, we mean the Whitney $C^r$-topology.

\subsubsection{Structural stability in a class for flows}

For a class $\mathcal{C}$ of $C^r$-flows on a manifold, a flow $v \in \mathcal{C}$ is {\bf $C^r$-structurally stable} in $\mathcal{C}$ if there is a \nbd $\mathcal{U}$ of $v$ in $\mathcal{C}$ with respect to the induced topology of the $C^r$-topology such that any flow in $\mathcal{U}$ is topologically equivalent to $v$. 
A flow is {\bf structurally stable} in $\mathcal{C}$ if it is $C^1$-structurally stable in $\mathcal{C}$. 

\subsubsection{Orbit arc and local topological equivalence}

A closed interval contained in an orbit is called an {\bf orbit arc}. 
The restriction of a flow $v$ to a subset $U$ of $S$ is {\bf locally topologically equivalent} to one of a flow $w$ to a subset $V$ on a topological space $T$ if there is a homeomorphism $h \colon U \to V$ such that the image of any orbit arc of $v$ in $U$ is an orbit arc of $w$, the inverse image of any orbit arc of $w$ in $V$ is an orbit arc of $v$, and $h$ and $h^{-1}$ preserve the directions of orbit arc of $v$ and $w$. 

\subsubsection{Double of a flow on a compact surface}

Since any flow on compact surfaces perverses the boundary components, a flow $v$ on a compact surface $S$ implies {\bf double} $v_{\mathrm{dbl}}$ on the double $S_{\mathrm{dbl}}$ of a compact surface $S$.

\subsubsection{Recurrence and relative concepts}
The {\bf $\bm{\omega}$-limit} (resp. {\bf $\bm{\alpha}$-limit}) {\bf set} of a point $x$ is $\omega(x) := \bigcap_{n\in \mathbb{R}}\overline{v(\R_{>n}, x)}$ (resp.  $\alpha(x) := \bigcap_{n\in \mathbb{R}}\overline{v(\R_{<n}, x)}$), where the closure of a subset $A$ is denoted by $\overline{A}$.
A point $x$ is {\bf recurrent} if $x \in \alpha(x) \cup \omega(x)$.
Denote by $\bm{\mathrm{P}(v)}$ (resp. $\bm{\mathrm{R}(v)}$) the set of (resp. non-closed recurrent) non-recurrent points.
Notice that $M = \mathop{\mathrm{Cl}}(v) \sqcup \mathrm{P}(v) \sqcup \mathrm{R}(v)$, where $\sqcup$ denotes a disjoint union.

\begin{definition}
An orbit is a {\bf separatrix} if it is a non-singular orbit from or to a singular point.
\end{definition}

\begin{definition}
A periodic orbit $O$ is a {\bf limit cycle} if there is a point $x \notin O$ with either $O = \omega(x)$ or $O = \alpha(x)$. 
\end{definition}

By a {\bf non-trivial circuit}, we mean either a periodic orbit or an image of an oriented circle by a continuous orientation-preserving mapping which is a directed graph but not a singleton and which is the union of separatrices and finitely many singular points. In other words, a non-trivial circuit which is not a periodic orbit is a directed path as a graph whose initial point is also terminal.

\subsubsection{Gradient vector fields}
For any $r \in \mathbb{Z}_{\geq 0} \sqcup \{ \infty \}$, a $C^r$-vector field $X$ on a Riemannian manifold $(M, g)$ is a {\bf gradient $\bm{C^r}$-vector field} if there is a function $f$ on $M$ with $g(X, \cdot) = d f$.
\begin{definition}
A flow on a manifold is {\bf gradient} if it is topologically equivalent to a flow generated by a $C^\infty$ gradient vector field.
\end{definition}

Notice that any gradient vector fields have no non-trivial circuits because of the existence of the height function. 
Moreover, we will show that any flow on a compact surface that is generated by a gradient $C^0$-vector field is topologically equivalent to a gradient flow (see Lemma~\ref{lem:grad_vf}). 
This means that the space of topologically equivalence classes of flows generated by gradient $C^s$-vector fields and the space of topologically equivalence classes of flows generated by gradient $C^r$-vector fields for any $s \leq r \in \mathbb{Z}_{\geq 0} \sqcup \{ \infty \}$ are canonically bijective.  
Notice that local Lipschitz continuity of vector fields guarantees the uniqueness and the existence of their flows, though the unique existence of the flows of continuous vector fields need not be guaranteed. 
%

\subsubsection{Flow box and transverse annulus}

The restriction of a flow to a disk is a {\bf trivial flow box} if it is locally topologically equivalent to the restriction of the flow generated by a vector field $(0,1) = \partial/\partial_y \in T (I \times J)$ to a square $I \times J$ where $I$ and $J$ are non-degenerate intervals. 
The disk is also called a {\bf trivial flow box}.
The restriction of a flow to an annulus is a {\bf transverse annulus} if it is locally topologically equivalent to the restriction of the flow generated by a vector field $(-x,-y)$ to an annulus $\{ (x,y) \in \R^2 \mid x^2 + y^2 \in I \}$, where $I \subset (0,1)$ is a non-degenerate interval. 
The annulus is also called a {\bf transverse annulus}.

\subsubsection{Types of singular points}
A singular point outside of the boundary is a {\bf sink} (resp. {\bf source}) if there is its open \nbd to which the restriction of the flow is locally topologically equivalent to a flow on some Euclidean space $\R^n$ generated by a vector field $X(x) = -x \in T \R^n$ (resp. $X(x) = x \in T \R^n$). 
A singular point on the boundary is a {\bf $\partial$-sink} (resp. {\bf $\partial$-source}) if there is its open \nbd to which the restriction of the flow is locally topologically equivalent to a flow on some $\R_{\geq 0} \times \R^{n-1}$ generated by a vector field $X(x) = -x$ (resp. $X(x) = x$). 
A singular point is {\bf attracting} if it is either a sink or a $\partial$-sink, and is {\bf repelling} if it is either a source and a $\partial$-source as in Figure~\ref{non-deg-sing}.

\subsubsection{Indices of isolated singular points of $C^1$-flows on manifolds}

We define the index of an isolated singular point for a $C^1$-flow as follows. 

\begin{definition}\label{index_C1_int}
The {\bf index} of an isolated singular point $x_0 \notin \partial M$ for a $C^1$-flow $w$ on a manifold $M$ is defined as the degree of the map $u: \partial U \to \mathbb{S}^{n-1}$ by 
\[
u(x) :=  
\dfrac{\partial w(0,x)}{\partial t} \cdot \left| \dfrac{\partial w(0,x)}{\partial t} \right|^{-1}
\]
for some isolated \nbd $U$ of $x_0$ which is a closed disk and whose boundary contains no singular points. 
\end{definition}

Notice that the index in the previous definition can be defined if 
\[
\dfrac{\partial w(0, \cdot)}{\partial t} \colon \partial U \to \R^n - \{0 \}
\]
is well-defined and continuous. 
Therefore, we can define the index of a point on the boundary $\partial M$ for a $C^1$-flow by taking the doble of a manifold $M$.  

\begin{definition}\label{index_C1_bd}
Let $w$ be a $C^1$-flow $w$ on a manifold $M$. 
The {\bf index} of an isolated singular point $x_0 \in \partial M$ for $w$ is half of the index of $[x_0] \in \widetilde{M}$ for the resulting flow $\widetilde{w}$ on $\widetilde{M}$, where $\widetilde{M}$ is the double of $M$ and $\widetilde{w}$ is the resulting flow on $\widetilde{M}$. 
\end{definition}

Since any connected sum of two spheres whose dimensions are the same is also a sphere, notice that, for a $C^1$-flow $w$ on a manifold and for any closed disk $U$ whose boundary contains no singular points and which contains at most finitely many singular points, the degree of the map $u: \partial U \to \mathbb{S}^{n-1}$ by 
\[
u(x) := \dfrac{\partial w(0,x)}{\partial t} \cdot \left| \dfrac{\partial w(0,x)}{\partial t} \right|^{-1}
\]
is the sum of indices contained in the singular points. 

\subsubsection{An continuous extension of $u(x)$}

Let $w$ be a $C^1$-flow on a manifold and $B$ a closed disk whose boundary contains no singular points and which contains at most finitely many singular points.
\begin{definition}\label{def:conti}
Define a partial map $V^w_B \colon \R \times \partial B \rightharpoonup \mathbb{S}^{n-1}$ as follows: 
\[
V^w_B(t,x) := 
\begin{cases}
     \dfrac{\partial w(0,x)}{\partial t}  \cdot \left| \dfrac{\partial w(0,x)}{\partial t} \right|^{-1} = u(x) & \text{if } t= 0 \\
    \dfrac{w(t,x) - x}{ \vert w(t,x) - x \vert } & \text{if } t > 0 \\
   - \dfrac{w(t,x) - x}{ \vert w(t,x) - x \vert } & \text{if } t < 0
\end{cases}
\]
\end{definition}

Notice that $V^w_B$ is continuous (see Corollary~\ref{cor:inv_index_sum} for details).

\subsection{Flows on surfaces}

From now on, we suppose that flows and vector fields are on surfaces unless otherwise stated.

\subsubsection{Types of singular points}

A singular point is a {\bf center} if there is its open \nbd to which the restriction of the flow is locally topologically equivalent to a flow on the open unit disk generated by a vector field $X(x_1,x_2) = (-x_2, x_1)$. 

\begin{figure}
\begin{center}
\includegraphics[scale=0.55]{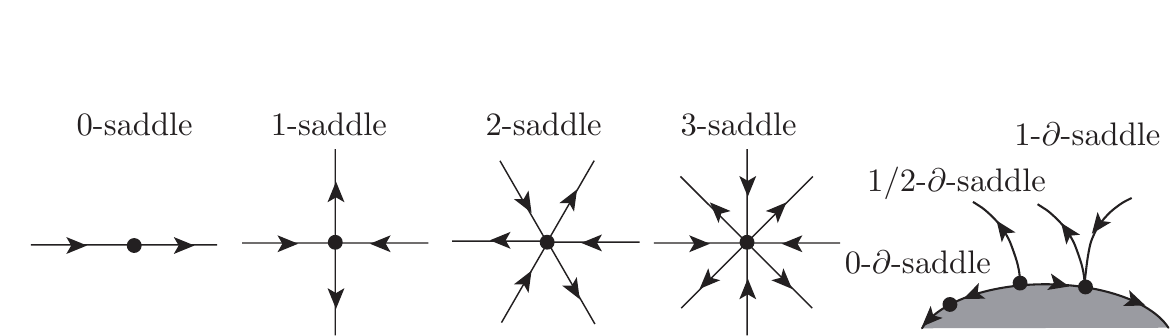}
\end{center}
\caption{Multi-saddles.}
\label{multi-saddles}
\end{figure}

A singular point is a {\bf $\bm{k}$-$\bm{\partial}$-saddle} (resp. {\bf $\bm{k}$-saddle}) if it is an isolated singular point on (resp. outside of) $\partial S$ with exactly $(2k + 2)$-separatrices, counted with multiplicity.
A singular point is a {\bf multi-saddle} if it is either $k$-saddle or a $k/2$-$\partial$-saddle for some $k \in \mathbb{Z}_{\geq 0}$. 
Note that each multi-saddle has at most finitely many separatrices, as in Figure~\ref{multi-saddles}.

\begin{definition}
A flow is {\bf quasi-regular} if it is topological equivalent to a flow such that each singular point is either a multi-saddle, a center, a sink, a source, a $\partial$-sink, or a $\partial$-source. 
\end{definition}

\begin{figure}
\begin{center}
\includegraphics[scale=0.35]{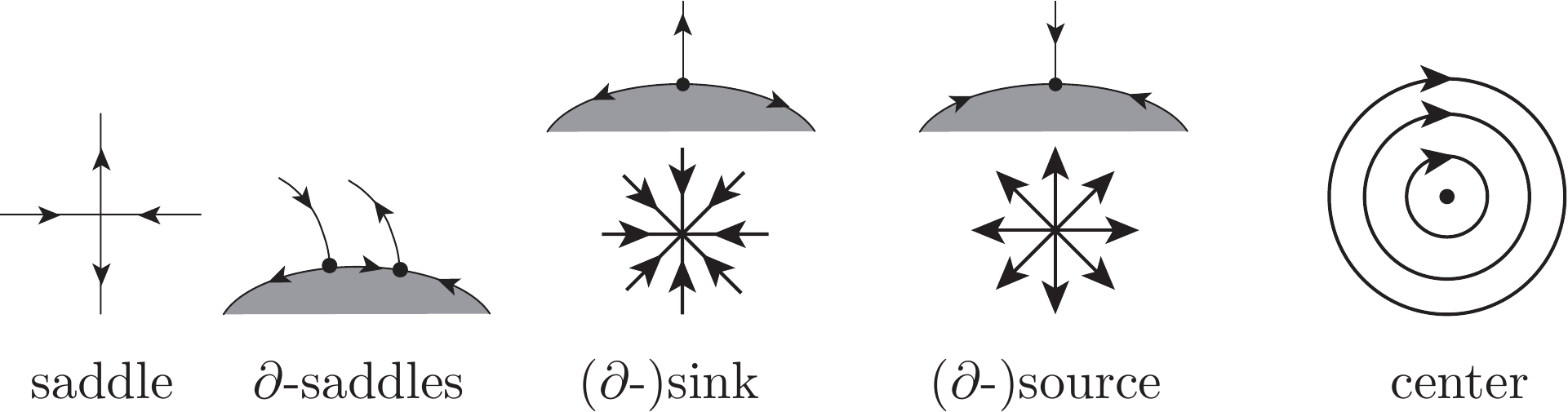}
\end{center}
\caption{A saddle, two $\partial$-saddles, a sink, a $\partial$-sink, a source, a $\partial$-source, and a center}
\label{non-deg-sing}
\end{figure}

\subsubsection{Multi-saddle separatrices}
A separatrix is a {\bf multi-saddle separatrix} if it is from and to multi-saddles.
A separatrix is a {\bf semi-multi-saddle separatrix} if it is from or to a multi-saddle.
The {\bf multi-saddle connection diagram} $D(v)$ is the union of multi-saddles and semi-multi-saddle separatrices. 
A {\bf multi-saddle connection} is a connected component of the multi-saddle connection diagram. 

\subsubsection{Sectors for flows on surfaces}

The restriction $v \vert _A$ for a subset $A$ is a {\bf sector} for a singular point $x$ if there are a non-degenerate interval $I \subseteq [0, 2 \pi )$ and a homeomorphism $h \colon \{ x\} \sqcup A \to \{ 0 \} \sqcup \{ ( r \cos \theta, r \sin \theta ) \in \R^2 \mid r \in (0,1), \theta \in I \}$ such that $h^{-1}(0) = x$.

\begin{figure}
\begin{center}
\includegraphics[scale=0.2]{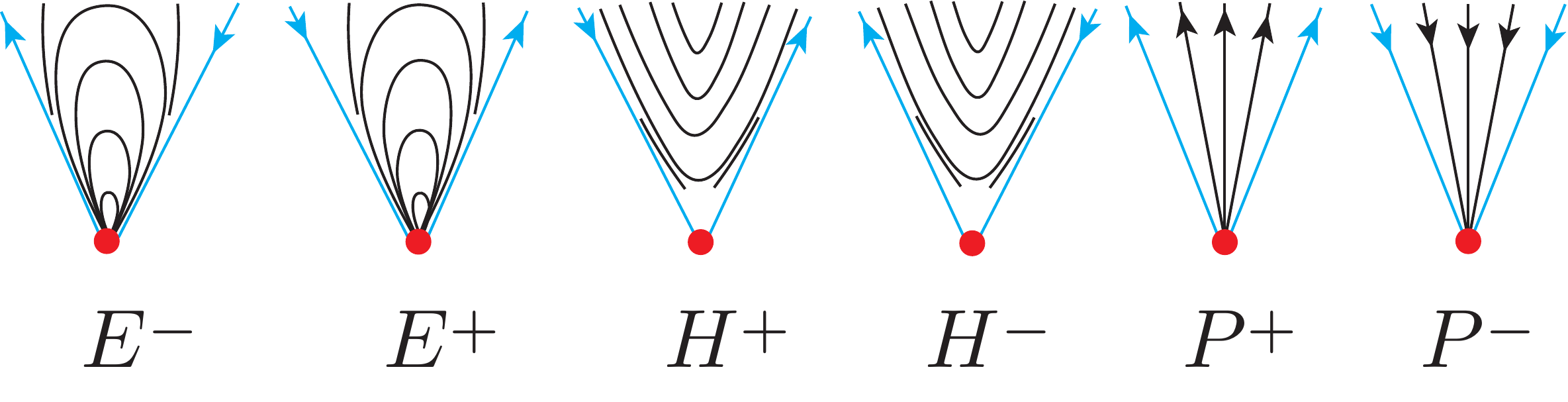}
\end{center}
\caption{Two parabolic sectors $P^-$ and $P^+$, two hyperbolic sectors $H^-$ and $H^+$ with clockwise and anti-clockwise orbit directions, and two elliptic sectors $E^+$ and $E^-$ with clockwise and anti-clockwise orbit directions, respectively.}
\label{sectors_all}
\end{figure}

\begin{definition}
A sector on an open subset $A$ for a singular point $x$ is a {\bf parabolic} if the restriction of $v$ to $\{ x\} \sqcup A$ is locally topologically equivalent to the restriction on $\{ 0 \} \sqcup \{ ( r \cos \theta, r \sin \theta ) \in \R^2 \mid r \in (0,1), \theta \in [0,\pi/2] \}$ of the flow generated by either $X(x,y) := (x,y)$ or $X(x,y) := -(x,y)$ as on the right of Figure~\ref{sectors_all}. 
\end{definition}

\begin{definition}
A sector on an open subset $A$ for a singular point $x$ is a {\bf hyperbolic} if the restriction of $v$ to $\{ x\} \sqcup A$ is locally topologically equivalent to the restriction on $\{ 0 \} \sqcup \{ ( r \cos \theta, r \sin \theta ) \in \R^2 \mid r \in (0,1), \theta \in [0,\pi/2] \}$ of the flow generated by  $X(x,y) := (-x, y)$ as on the middle of Figure~\ref{sectors_all}. 
\end{definition}

\begin{definition}
A sector on an open subset $A$ for a singular point $x$ is {\bf elliptic} if the restriction of $v$ to $\{ x\} \sqcup A$ is locally topologically equivalent to the restriction on $\{ (x,y) \in \R^2 \mid x \geq 0, x^2 + y^2 > 1 \} \sqcup \{ \infty \} \subset \R^2 \sqcup \{ \infty \}$ of the flow generated by a continuous vector field $X$ defined by $X = (0,1)$ on $\R^2$ and $X = 0$ at the point $\infty$ at infinity as on the left of Figure~\ref{sectors_all}. 
\end{definition}

Here the union $\R^2 \sqcup \{ \infty \}$ is the one-point compactification of $\R^2$, which is a sphere. 

A separatrix is a {\bf hyperbolic} (resp. {\bf elliptic}, {\bf parabolic}) {\bf border separatrix} if it is contained in the boundary of a hyperbolic sector (resp. a maximal open elliptic sector, a maximal open parabolic sector). 

\begin{definition}
A singular point in the interior of the surface is {\bf finitely sectored} if either it is a center or there is its open neighborhood which is an open disk and is a finite union of the point, parabolic sectors, hyperbolic sectors, and elliptic sectors such that a pair of distinct sectors intersects at most two orbit arcs. 
\end{definition}

A singular point on the boundary of the surface is {\bf finitely sectored} if it is finitely sectored for the resulting flow on the double of the compact surface. 
Notice that quasi-regularity implies that each singular point is a finitely sectored singular point without elliptic sectors. 
A singular point is a multi-saddle if and only if it is a finitely sectored singular point whose sectors are hyperbolic.
Similarly, a singular point is a sink, a $\partial$-sink, a source, or a $\partial$-source if and only if it is a finitely sectored singular point whose sectors consist of exactly one parabolic sector as in Figure~\ref{non-deg-sing}.

In \cite[Theorem~A]{kibkalo2022topological}, any isolated singular point of a gradient vector field on a surface is characterized as a non-trivial finitely sectored singular point without elliptic sectors.

\subsubsection{Differentiability of gradient flows on surfaces}

We have the following observation. 

\begin{lemma}\label{lem:grad_vf}
Every flow on a compact surface which is generated by a gradient $C^0$-vector field with finitely many singular points is a gradient flow. 
\end{lemma}

\begin{proof}
 Let $v$ be a flow with finitely many singular points on a compact surface $S$ which is generated by a gradient $C^0$-vector field. 
By \cite[Proposition~3.3]{yokoyama2023dependency} and the dual statement, each of $\omega$-limit and $\alpha$-limit sets of points of $v$ is a singular point.
Therefore, any recurrent points are singular points and there are no non-trivial circuits. 
The finiteness of singular points implies that the set of non-recurrent points is open dense in $S$.  
By \cite[Lemma~3.2]{kibkalo2022topological}, the existence of the height function implies every singular point of the flow $v$ is a finitely sectored singular point without elliptic sectors.
By \cite[Theorem~B]{kibkalo2022topological}, the flow $v$ is a gradient flow. 
\end{proof}

The previous lemma implies the following equivalence. 

\begin{corollary}\label{cor:grad_vf}
The following statements are equivalent for any flow with finitely many singular points on a compact surface: 
\\
{\rm(1)} The flow is gradient. 
\\
{\rm(2)} The flow is topologically equivalent to a flow generated by a gradient $C^0$-vector field.  
\end{corollary}

\subsubsection{Transversality}
A differentiable subset $A$ on $S$ is transverse to a $C^1$-flow $v$ if $T_{x} S = T_{x} A + T_{x} O(x)$ for any $x \in A$ (i.e. the tangent space of any point in $A$ and the tangent space of its orbit span the tangent spaces for $S$). 

We define transversality for flows on surfaces using tangential spaces of surfaces because each flow on a compact surface is topologically equivalent to a $C^1$-flow by Gutierrez's smoothing theorem~\cite{gutierrez1978structural}.

A subset $A$ on $S$ is {\bf transverse} to a flow $v$ if there is a $C^1$-flow $w$ on $S$ which is topologically equivalent to $v$ via a homeomorphism $h \colon S \to S$ such that the image $h(A)$ is transverse to $w$. 

\subsubsection{Indices of isolated singular points of flows on surfaces}

\begin{figure}
\begin{center}
\includegraphics[scale=0.15]{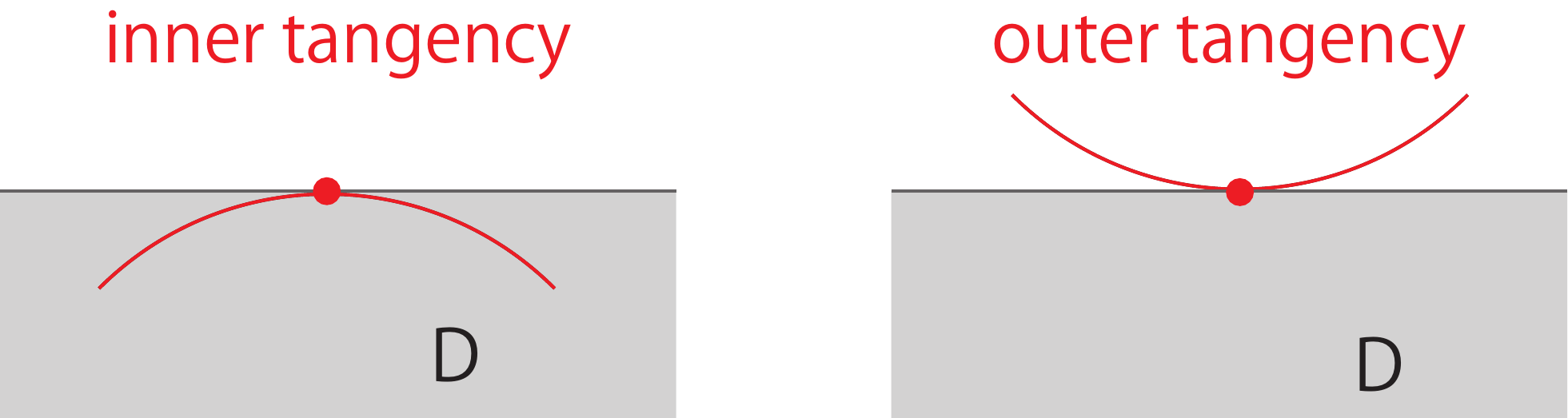}
\end{center}
\caption{Inner and outer tangencies.}
\label{fig:tangencies}
\end{figure}

A tangency $x$ of a curve $C$ to a loop bounding a closed disk $D$ is {\bf inner} (resp. {\bf outer}) if there is a small arc $I$ in $C$ containing $x$ such that the difference $I - \{ x \}$ is contained in the interior $\mathop{\mathrm{int}} D$ (resp. $(I - \{ x \}) \cap D = \emptyset$) as in left the (resp. right) on Figure~\ref{fig:tangencies}.
 
Recall the indices of isolated singular points of flows on surfaces (cf. \cite[\S~2.1.2]{kibkalo2022topological}). 

\begin{definition}\label{def:index_int}
The {\bf index} of an isolated singular point outside of the boundary $\partial S$ is the difference $(2 + n_{i} - n_{o})/2$, where $n_i$ is the number of inner tangencies and $n_o$ is the number of outer tangencies of a loop which is transverse at all but finitely many points and bounds an open disk containing the singular point. 
\end{definition}

Notice that the index of an isolated singular point outside of the boundary of a surface is independent of the choice of such a loop in the previous definition. 

\begin{definition}\label{def:index_bd}
For an isolated singular point $x \in \partial S$, the {\bf index} $\bm{\mathrm{ind}_v(x)}$ of $v$ at $x$ is defined by $\mathrm{ind}_v(x) := \mathrm{ind}_{v_{\mathrm{dbl}}}(x)/2$, where $v_{\mathrm{dbl}}$ is the induced flow on the double $S_{\mathrm{dbl}}$. 
\end{definition}

The following equivalence holds. 

\begin{lemma}[cf. Theorem 2.1\cite{izydorek1996note}]\label{lem:equiv_index}
For any isolated singular point $x \in S$ for a $C^1$-flow on a surface $S$, the index of $x$ in Definition~\ref{def:index_int} {\rm(resp.} Definition~\ref{def:index_bd}{\rm)} corresponds to one in Definition~\ref{index_C1_int} and {\rm(resp.} Definition~\ref{index_C1_bd}{\rm)}. 
\end{lemma}

The previous lemma also means that the index of a $C^1$-flow at an isolated singular point $x$ is independent of the choice of such an open neighborhood $U$. 

\subsubsection{Semi-attracting cycles and semi-repelling cycles}

A limit cycle $O$ is {\bf semi-attracting} on a side if there is a collar $U$ of $O$ such that $O = \omega(x)$ for any $x \in U$. 
A limit cycle $O$ is {\bf semi-repelling} on a side if there is a collar $U$ of $O$ such that $O = \alpha(x)$ for any $x \in U$. 
Recall that a limit cycle is {\bf topologically hyperbolic} if it is topologically attracting or repelling (i.e. it is either semi-attracting on each side or semi-repelling on each side). 

\subsubsection{Fake limit cycle, a fake multi-saddle, and a fake parabolic sector}

We define a fake limit cycle, a fake multi-saddle, and a fake parabolic sector as follows. 

\begin{figure}
\begin{center}
\includegraphics[scale=0.35]{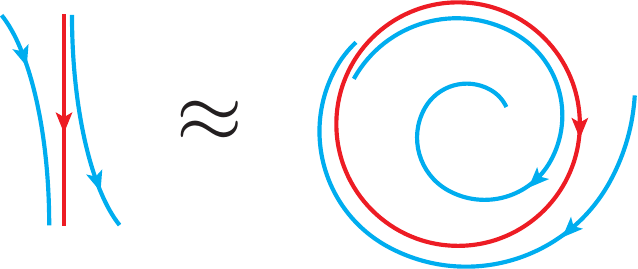}
\end{center}
\caption{Fake limit cycle}
\label{fake_limit_cycle}
\end{figure}

\begin{definition}
A limit cycle is a {\bf fake limit cycle} \cite{kibkalo2022topological} if it is semi-attracting on one side and semi-repelling on another side as in Figure~\ref{fake_limit_cycle}.
\end{definition}

\begin{definition}
A singular point is a {\bf fake multi-saddle} if it is either a $0$-saddle or a $0$-$\partial$-saddle as in Figure~\ref{non-deg-sing}. 
\end{definition}

\begin{figure}
\begin{center}
\includegraphics[scale=0.15]{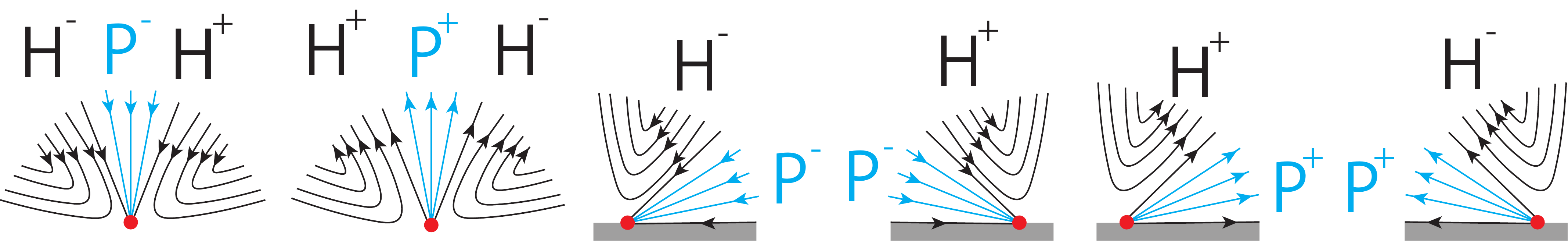}
\end{center}
\caption{Fake parabolic separatrices}
\label{fake_parabolic_sector}
\end{figure}

\begin{definition}\label{def:fake}
A parabolic sector $A$ for a singular point $x$ is {\bf fake} \cite{kibkalo2022topological} if there are separatrices $\gamma, \mu$ from/to $x$ satisfying the following conditions (as in Figure~\ref{fake_parabolic_sector}): 
\\
{\rm(1)} The separatrix $\gamma$ is a hyperbolic border separatrix. 
\\
{\rm(2)} The separatrix $\mu$ is either a hyperbolic border separatrix or is contained in the boundary $\partial S$. 
%
%
\end{definition}
By combination, any fake parabolic sector is locally topological equivalent to one of the sectors in Figure~\ref{fake_parabolic_sector}. 
Roughly speaking, a fake parabolic sector is a parabolic sector either between two hyperbolic sectors or between one hyperbolic sector and a boundary component. 
\cite[Theorem~A]{kibkalo2022topological} implies that a gradient flow with isolated singular points on a surface is quasi-regular if and only if there are no fake parabolic sectors. 

\subsubsection{Generic non-gradient on surfaces}

Notice that the set of gradient $C^1$-flows is not open in the set of $C^1$-flows because singular points need not be non-degenerate.
On the other hand, forbidding the existence of fake multi-saddles, fixing the sum of indices of repelling singular points (i.e. sources and $\partial$-sources) and the sum of indices of attracting singular points (i.e. sinks and $\partial$-sinks), we topologically characterize codimension $k$ flows of gradient flows in this paper.

\section{Invariance of sums of indices of isolated singular points under small perturbations}


For a subset $A$ of a metric space $(M,d)$, put 
\[
B_r (A) := \{x \in M \mid \inf_{a \in A } d (x, a) < r \}
\]
for any $A \subseteq M$ and $r\geq 0$. 
Recall that the subset 
\[
\mathcal{U}(K, U) := \{ w : \text{flow on } M \mid w(K) \subset U \}
\]
for any compact subset $K \subseteq \R \times M$ and any open subset $U \subseteq M$ is a $C^0$-open subset of the space of flows on $M$. 
We have the following statement. 

\begin{lemma}\label{lem:inv_index02}
Let $v$ be a flow on a compact metric space $(M,d)$. 
For any number $\delta > 0$, there is a $C^0$-\nbd $\mathcal{U}$ of $v$ such that $\bigcup_{w \in \mathcal{U}} \mathop{\mathrm{Sing}}(w) \subset B_{\delta} (\mathop{\mathrm{Sing}}(v))$. 
\end{lemma}

\begin{proof}
Fix a positive number $\delta > 0$. 
The complement $M_0 := M - B_{\delta} (\mathop{\mathrm{Sing}}(v))$ is a compact metric space. 
For any point $x \in M_0$, there are positive numbers $\varepsilon_0(x), \varepsilon_1(x), \varepsilon_2(x) >0$ such that $v(K_x) \subset U_x$ and $\overline{B_{\varepsilon_0(x)}(x)} \cap U_x = \emptyset$, where $K_x := \{ \varepsilon_1(x)\} \times \overline{B_{\varepsilon_0(x)}(x)}$ and $U_x := B_{\varepsilon_2(x)} (v^{\varepsilon_1(x)}(x))$.
Then 
\[
\bigcup_{w \in \mathcal{U}(K_x, U_x)} \mathop{\mathrm{Sing}}(w) \cap \overline{B_{\varepsilon_0(x)}(x)} = \emptyset
\]
for any $x \in M_0$.
By construction, we have $v \in \mathcal{U}(K_x, U_x)$ for any $x \in M_0$. 
Since $\{ B_{\varepsilon_0(x)}(x) \mid x \in M_0 \}$ is an open cover of the compact subset $M_0$, there are finitely many points $x_1, \ldots , x_k \in M_0$ such that $M_0 \subset \bigcup_{i=1}^k B_{\varepsilon_0(x_i)}(x_i)$. 
By $v \in \bigcap_{i=1}^k \mathcal{U}(K_{x_i}, U_{x_i})$, the intersection $\mathcal{U} := \bigcap_{i=1}^k \mathcal{U}(K_{x_i}, U_{x_i})$ is an open \nbd of $v$ such that $\emptyset = \bigcup_{w \in \mathcal{U}} \mathop{\mathrm{Sing}}(w) \cap \bigcup_{i=1}^k \overline{B_{\varepsilon_0(x_i)}(x_i)} \supseteq \bigcup_{w \in \mathcal{U}} \mathop{\mathrm{Sing}}(w) \cap M_0$ and so that $\bigcup_{w \in \mathcal{U}} \mathop{\mathrm{Sing}}(w) \subseteq B_{\delta} (\mathop{\mathrm{Sing}}(v))$. 
\end{proof}

For a subset $C$ of a metric space $(M,d)$, write the $\varepsilon$-\nbd $B_{\varepsilon}(C)$ of $C$ (i.e. $B_{\varepsilon}(C) := \bigcup_{c \in C} \{ y \in M \mid d(y, c) < \varepsilon \}$). 
We have the following invariance of indices under small perturbations, which partially generalizes \cite[Lemma~2.1]{yokoyama2021combinatorial}. 

\begin{lemma}\label{lem:inv_index}
Let $v$ be a flow on a closed surface $S$ and $x \in S$ an isolated singular point. 
For any loop $C$ which is transverse at all but finitely many points and bounds an open disk $D$ whose closure is an isolated \nbd of $x$, there is a $C^0$-\nbd $\mathcal{U}$ of $v$ satisfying the following statements for any flow $w \in \mathcal{U}$: 
\\
{\rm(1)} $\mathrm{ind}_v(x) = \sum_{y \in D \cap \mathop{\mathrm{Sing}}(w)} \mathrm{ind}_w(y)$ if $w$ has only isolated singular points on $\overline{D}$. 
\\
{\rm(2)} For any small positive number $\varepsilon >0$, there is a loop $C' \subset B_{\varepsilon}(C)$ which is transverse at all but finitely many points such that the number of the inner {\rm(resp.} outer{\rm)} tangencies on $C$ of $v$ equals one of the inner {\rm(resp.} outer{\rm)} tangencies on $C'$ of $w$, where $B_{\varepsilon}(C)$ is the $\varepsilon$-\nbd of $C$. 
\end{lemma}

Moreover, we have the folllowing invariance. 

\begin{corollary}\label{cor:inv_index_02}
Let $v$ be a gradient flow on a compact surface $S$, $x \in S$ an attracting {\rm(resp.} repelling{\rm)} singular point, and $D$ an open disk whose closure is an isolated \nbd of $x$ such that $\partial D \cap \partial S = \emptyset$ unless $x \in \partial S$, and that $\partial D \setminus \partial S$ is transverse. 
Then there is a $C^0$-\nbd $\mathcal{U}$ of $v$ satisfying the following statements for any flow $w \in \mathcal{U}$: 
\\
{\rm(1)} If $D \cap \mathop{\mathrm{Sing}}(w)$ contains exactly one singular point $x_w$, then the singular point $x_w$ is attracting {\rm(resp.} repelling{\rm)}. 
\end{corollary}

Since the proofs of the previous lemma and corollary are straightforward but technical, the proofs are stated in Appendix~\ref{prf:lem_inv_ind}.

\subsection{Subspaces}

For any $r \in \mathbb{Z}_{\geq 0} \sqcup \{ \infty \}$, let $\bm{\mathcal{G}^r_{> -1}(S)}$ be the set of gradient $C^r$-flows with finitely many singular points without fake multi-saddles or fake parabolic sectors on a compact surface $S$.
Equip $\mathcal{G}^r_{> -1}(S)$ with the $C^r$ topology for any $r \in \mathbb{Z}_{\geq 0} \sqcup \{ \infty \}$. 
By \cite[Theorem~A]{kibkalo2022topological}, each singular point of any gradient flows in $\mathcal{G}^r_{> -1}(S)$ is either attracting singular point (i.e. either a $\partial$-source, a source), repelling singular point (i.e. a $\partial$-)sink, a sink), or a multi-saddle. 
In particular, any gradient flows in $\mathcal{G}^r_{> -1}(S)$ is quasi-regular.  
%

\subsubsection{Complete invariance of the union of the multi-saddle connection diagram and positively indexed singular points}

An isolated singular point of a gradient flow is {\bf positively indexed} if its index is positive.
We observe the complete invariance of the union of the multi-saddle connection diagram and attracting/repelling singular points as follows by a similar argument in the proof of \cite[Lemma~3.1]{yokoyama2021combinatorial}.  

\begin{lemma}\label{lem:fin_comb}
Let $S$ be a compact surface. 
The union of the multi-saddle connection diagram and positively indexed singular points is a complete invariant of $\mathcal{G}^r_{> -1}(S)$ as a surface graph for any $r \in \Z_{\geq 0} \sqcup \{ \infty \}$. 
Moreover, the complement of such a union is a finite disjoint union of invariant open trivial flow boxes and invariant open transverse annuli. 
\end{lemma}

\begin{proof}
Let $v$ be a gradient $C^r$-flows with finitely many singular points without fake multi-saddles or fake parabolic sectors on the compact surface $S$ (i.e. $v \in \mathcal{G}^r_{> -1}(S)$), and denote by $D(v)$ the multi-saddle connection diagram. 
We may assume that $S$ is connected. 
From $S = \mathop{\mathrm{Sing}}(v) \sqcup \mathop{\mathrm{Per}}(v) \sqcup \mathrm{P}(v) \sqcup \mathrm{R}(v)$, by the existence of the height function of a gradient vector field, there are no non-singular recurrent orbits and so $S =\mathop{\mathrm{Sing}}(v)  \sqcup \mathrm{P}(v)$. 
Because the double of any invariant trivial flow box is either a trivial flow box or a transverse annulus, taking the double of the surface $S$ if necessary, we may assume that $S$ is closed. 

\begin{claim}\label{claim:01}
The multi-saddle connection diagram $D(v)$ consists of finitely many orbits such that the set difference $\mathop{\mathrm{Sing}}(v) \setminus D(v)$ is the set of attracting or repelling singular points. 
\end{claim}
\begin{proof}[Proof of Claim~\ref{claim:01}]
By \cite[Theorem~A]{kibkalo2022topological}, each singular point of $v$ is either an attracting or repelling singular point or a multi-saddle. 
In particular, the flow $v$ is quasi-regular. 
Note that any attracting or repelling singular points are positively indexed. 
By \cite[Proposition 3.3]{yokoyama2023dependency} and the dual statement, each of the $\omega$-limit set and the $\alpha$-limit set of any point is a singular point. 
Since $v$ is a gradient flow on a compact surface, there are attracting singular points and repelling singular points. 
By the finite existence of singular points, there are at most finitely many semi-multi-saddle separatrices. 
\end{proof}

Then the complement $C := S - (\mathop{\mathrm{Sing}}(v) \cup D(v)) \subseteq \mathrm{P}(v)$ is a finite disjoint union of invariant open connected subsets such that the orbit of any point in the complement $C$ connects a repelling point and an attracting point. 
Therefore, it suffices to show that every connected component of the complement of such a union is either an invariant open transverse annulus or an invariant open trivial flow box. 
(Notice that the assertion can be reduced from \cite[Theorem~A]{yokoyama2017decompositions}. However, we state the proof because the proof is simple and straightforward in our case.)

\begin{claim}\label{claim:02}
We may assume that any \nbd of an attracting or repelling singular point intersects $D(v)$. 
\end{claim}
\begin{proof}[Proof of Claim~\ref{claim:02}]
Assume that there are \nbds of attracting or repelling singular points that do not intersect $D(v)$. 
Then there is a closed transversal $\gamma \subset S - (\mathop{\mathrm{Sing}}(v) \cup D(v))$ near a repelling singular point. 
The saturation $v(\gamma)$ is an invariant transverse annulus such that the $\omega$-limit (resp. the $\alpha$-limit) set of any point of $v(\gamma)$ is a sink (resp. source). 
This means that the surface $S$ is a sphere and so that $C$ is an invariant open transverse annulus.
\end{proof}

\begin{claim}\label{claim:03}
Any connected component of $C$ is a trivial flow box. 
\end{claim}
\begin{proof}[Proof of Claim~\ref{claim:03}]
For any point $x \in C$, there is a small closed disk $B$ which is an isolated \nbd of an attracting point, and there is a connected component $I$ of the set difference $\partial B \setminus \cup D(v)$ which is an open transverse arc such that $I$ intersects $O(x)$. 
Therefore, for any connected component $U$ of $C$, there is a small closed disk $B$ which is an isolated \nbd of an attracting point and there is a connected component $I$ of the set difference $\partial B \setminus \cup D(v)$ which is an open transverse arc such that $v(I) = U$. 
This means that any connected component of $C$ is a trivial flow box. 
\end{proof}

Therefore, the surface graph $\mathop{\mathrm{Sing}}(v) \cup D(v)$ is a complete invariant as a surface graph. 
\end{proof}

\subsubsection{Invariance of positively indexed singular points under small perturbations}

For gradient flows $v, w$ with finitely many singular points on a compact surface $S$ and for any disjoint open subsets $U_1, \ldots , U_l$ of $S$, the gradient flow $w$ {\bf has same type of positively indexed singular points} of $v$ with respect to $U_1, \ldots , U_l$ if following conditions hold: 
\\
{\rm(1)} The disjoint union $\bigsqcup_{i=1}^l U_i$ is a \nbd of $\mathop{\mathrm{Sing}}(v) \cup \mathop{\mathrm{Sing}}(w)$. 
\\
{\rm(2)} $ \vert U_i \cap \mathop{\mathrm{Sing}}(v) \vert  = 1 \leq  \vert U_i \cap \mathop{\mathrm{Sing}}(w) \vert $ for any $i \in \{1, \ldots , l \}$. 
\\
{\rm(3)} For any $i \in \{1, \ldots , l \}$, the singular point $x_i$ in $U_i \cap \mathop{\mathrm{Sing}}(v)$ is a positively indexed singular point if and only if $ \vert U_i \cap \mathop{\mathrm{Sing}}(w) \vert  = 1$ and the singular point $x'_i$ in $U_i \cap \mathop{\mathrm{Sing}}(w)$ is a positively indexed singular point. 
Then we call $x'_i$ the {\bf continuation} of $x_i$ with respect to $w$. 

\subsubsection{Codimension $k$ subspaces}

For any integer $k \in \mathbb{Z}_{\geq 0}$, denote by $\bm{\mathcal{G}^r_{k, > -1}(S)} \subseteq \mathcal{G}^r_{> -1}(S)$ the subset of (quasi-regular) gradient flows whose sums of indices of positively indexed singular points 
are $k$. 
%
%
We observe that any small perturbation can not merge a pair of distinct singular points for a flow in $\mathcal{G}^r_{> -1}(S)$. 
More precisely, we have the following statement by a similar argument in the proof of \cite[Lemma~3.2]{yokoyama2021combinatorial}. 

\begin{lemma}\label{lem:perturbation_closed}
Suppose that $S$ has no boundary. 
For any $v \in \mathcal{G}^r_{k, > -1}(S)$ and any disjoint open disks $D_i$ which are isolated \nbds of $\mathop{\mathrm{Sing}}(v)$, there is its $C^0$-\nbd in $\mathcal{G}^r_{k, > -1}(S)$ each element of which has same type of positively indexed singular points of $v$ with respect to $D_i$. 
\end{lemma}

\begin{proof}
Let $v$ be a gradient flow with finitely many singular points on a closed surface $S$ and $H_v$ the height function of $v$. 
By \cite[Theorem~A]{kibkalo2022topological}, each singular point of $v$ is either an attracting or repelling singular point or a multi-saddle. 
Then the set $\mathop{\mathrm{Sing}}_{\mathrm{c}}(v)$ of positively indexed singular points is the complement in $\mathop{\mathrm{Sing}}(v)$ of multi-saddles. 
Let $x_1, x_2, \ldots , x_{k'}$ be the positively indexed singular points of the flow $v$ and $x_{k'+1}, x_{k'+2}, \ldots , x_l$ the multi-saddles of $v$.  
Lemma~\ref{lem:inv_index02} implies that, for any number $\varepsilon > 0$, there is a $C^0$-\nbd $\mathcal{U} \subseteq \mathcal{G}^r_{k, > -1}(S)$ of $v$ such that the open \nbds $D_i := B_{\varepsilon}(x_i)$ for any $i \in \{ 1, \ldots , l \}$ are open disks whose boundaries are loops, and that $\bigcup_{w \in \mathcal{U}} \mathop{\mathrm{Sing}}(w) \subset B_{\varepsilon} (\mathop{\mathrm{Sing}}(v)) = \bigsqcup_{i=1}^l D_i$. 
Since $\mathcal{U} \subseteq \mathcal{G}^r_{k, > -1}(S)$, any flow $w \in \mathcal{U}$ has at most finitely many singular points and such singular points are isolated. 
Applying Lemma~\ref{lem:inv_index} to loops $ \partial B_{\varepsilon/2}(x_i)$, there is a \nbd $\mathcal{V}$ of $v$ in $\mathcal{G}^r_{k, > -1}(S)$ such that $\mathrm{ind}_v(x_i) = \sum_{y \in D_i \cap \mathop{\mathrm{Sing}}(w)} \mathrm{ind}_w(y)$ for any $w \in \mathcal{V}$. 
Fix any $w \in \mathcal{V}$. 
Since $x_1, x_2, \ldots , x_{k'}$ are the positively indexed singular points of the flow $v$, we have that $1 = \sum_{y \in D_i \cap \mathop{\mathrm{Sing}}(w)} \mathrm{ind}_w(y)$ for any $i \in \{ 1,2, \ldots , k'\}$. 
This means that any $D_1, D_2, \ldots , D_{k'}$ contain at least one positively indexed singular points of $w$. 
Since $w$ contains exactly $k'$ positively indexed singular points, the open disks $D_{1}, D_{2}, \ldots , D_{k'}$ contain exactly one positively indexed singular point of $w$ and the open disks $D_{k'+1}, D_{k'+2}, \ldots , D_l$ contain no positively indexed singular points of $w$. 
This implies the assertion. 
\end{proof}

\begin{lemma}\label{lem:perturbation}
For any $v \in \mathcal{G}^r_{k, > -1}(S)$ and any disjoint open disks $D_i$ which are isolated \nbds of $\mathop{\mathrm{Sing}}(v)$, there is its $C^0$-\nbd in $\mathcal{G}^r_{k, > -1}(S)$ each element of which has same type of positively indexed singular points of $v$ with respect to $D_i$. 
\end{lemma}

\begin{proof}
Notice that the $C^0$-topology of $\mathcal{G}^r_{k, > -1}(S)$ corresponds to the induced topology of the $C^0$-topology of $\mathcal{G}^0_{k, > -1}(S)$ to the subspace $\mathcal{G}^r_{k, > -1}(S)$. 
By Lemma~\ref{lem:perturbation_closed}, we may assume that $S$ has a nonempty boundary. 
For any $w \in \mathcal{G}^r_{k, > -1}(S)$, write the double $w_{\mathrm{dbl}}$ on the double $S_{\mathrm{dbl}}$ of $S$. 
Fix any $v \in \mathcal{G}^r_{k, > -1}(S)$. 

\begin{claim}\label{claim:04}
For any $C^0$-open \nbd $\mathcal{V}$ of $v_{\mathrm{dbl}}$ in $\mathcal{G}^0_{k, > -1}(S_{\mathrm{dbl}})$, the subset $\{ w \in \mathcal{G}^r_{k, > -1}(S) \mid w_{\mathrm{dbl}} \in \mathcal{V} \}$ is a $C^0$-open \nbd of $v$ in $\mathcal{G}^r_{k, > -1}(S)$. 
\end{claim}

\begin{proof}[Proof of Claim~\ref{claim:04}]
Fix any $C^0$-open subset $\mathcal{V}$ in $\mathcal{G}^0_{k, > -1}(S_{\mathrm{dbl}})$. 
By construction, the subset $\mathcal{V}' := \{ w \in \mathcal{G}^0_{k, > -1}(S) \mid w_{\mathrm{dbl}} \in \mathcal{V} \}$ is a $C^0$-open subset $\mathcal{G}^0_{k, > -1}(S)$. 
Since the $C^0$-topology of $\mathcal{G}^r_{k, > -1}(S)$ corresponds to the induced topology of the $C^0$-topology of $\mathcal{G}^0_{k, > -1}(S)$ restricted to the subspace $\mathcal{G}^r_{k, > -1}(S)$, the restriction $\mathcal{V}' \cap \mathcal{G}^r_{k, > -1}(S) = \{ w \in \mathcal{G}^r_{k, > -1}(S) \mid w_{\mathrm{dbl}} \in \mathcal{V} \}$ is a $C^0$-open subset in $\mathcal{G}^r_{k, > -1}(S)$. 
\end{proof}

By Lemma~\ref{lem:perturbation_closed}, there is a $C^0$-\nbd $\mathcal{V}$ in $\mathcal{G}^0_{k, > -1}(S)$ of $v_{\mathrm{dbl}}$ each element of which has same type of positively indexed singular points of $v_{\mathrm{dbl}}$. 
From the previous claim, the subset $\{ w \in \mathcal{G}^r_{k, > -1}(S) \mid w_{\mathrm{dbl}} \in \mathcal{V} \}$ is a $C^0$-open \nbd of $v$ in $\mathcal{G}^r_{k, > -1}(S)$. 
\end{proof}

The previous lemma means that the subspaces $\mathcal{G}^r_{k, > -1}(S)$ correspond to the spaces without creations and annihilations of singular points because of Poincar{\'e}-Hopf theorem. 

\subsection{Existence of associated neighborhoods of the multi-saddle connections}

For any vector field $v \in \mathcal{G}^r_{> -1}(S)$, write $\bm{N(v)}$ the number of multi-saddle connections of $v$, and denote by $\bm{\mathop{\mathrm{Sing}}_+(v)}$ the set of positively indexed singular points of $v$ and put $\mathop{\mathrm{Sing}}_-(v) := \mathop{\mathrm{Sing}}(v) - \mathop{\mathrm{Sing}}_+(v)$. 

\subsubsection{Codimensions of multi-saddles}
Let $w$ be a gradient flow on a surface $T$ with finitely many singular points.  
The {\bf codimension} of a $k$-saddle $y$ ($k \geq 1$) is $\bm{\mathrm{codim}_{\mathrm{m},w}(y)} := 2(k-1) = -2(1+ \mathrm{ind}(y))$. 
The {\bf codimension} of a $\partial$-$(l/2)$-saddle $y$ ($l \geq 1$) is $\bm{\mathrm{codim}_{\mathrm{m},w}(y)} := l-1 = -2(1/2 + \mathrm{ind}(y))$. 
Equivalently, the codimension of a $\partial$-$k$-saddle $y$ ($k \geq 1/2$) is $2k-1 = -2(1/2 + \mathrm{ind}(y))$.

\subsubsection{Codimensions of gradient flows}
Let $w$ be a gradient flow on a surface $T$ with finitely many singular points.  
We define codimensions as follows. 

\begin{definition}
The codimension $\bm{\mathrm{codim}_{\mathrm{m}}(w)} \in \Z_{\geq 0}$ of $w$ with respect to multiplicity is defined as the sum of the codimensions of multi-saddles of $w$. 
\end{definition}

\begin{definition}
The codimension $\bm{\mathrm{codim}_{\mathrm{h}}(w)} \in \Z_{\geq 0}$ with respect to heteroclinicity is defined as the number of multi-saddle separatrices outside of the boundary $\partial T$. 
\end{definition}

\begin{definition}
The codimension $\bm{\mathrm{codim}_{\mathrm{f}}(w)} \in \Z_{\geq 0}$ with respect to fakeness is defined as the sum of the number of fake multi-saddles and fake parabolic sectors. 
\end{definition}



\begin{definition}
The codimension $\bm{\mathrm{codim}(w)}$ of the gradient flow $w$ is defined as follows: 
\[
\mathrm{codim}(w) := \mathrm{codim}_{\mathrm{m}}(w) + \mathrm{codim}_{\mathrm{h}}(w) + \mathrm{codim}_{\mathrm{f}}(w) 
\]
\end{definition}

Since $\mathcal{G}^r_{> -1}(S)$ contains no fake multi-saddles, we have the following observation. 

\begin{lemma}
For any flow $w \in \mathcal{G}^r_{> -1}(S)$, we have the following equality: 
\[
\mathrm{codim}(w) = \mathrm{codim}_{\mathrm{m}}(w) + \mathrm{codim}_{\mathrm{h}}(w)
\]
\end{lemma}

\subsubsection{Annihilations of non-self-connected separatrices}

For a non-self-connected separatrix $\gamma$ from a multi-saddle $\alpha$ to a multi-saddle $\omega$ of a flow $v$ and a \nbd $U$ of the closure $\overline{\gamma} = \gamma \sqcup \{ \alpha, \omega \}$, by replacing a trivial flow box in $U$ intersecting only the separatrix $\gamma$ and no other separatrices into a trivial flow box $B_w$, we call that the resulting flow $w$ is the flow obtained by an {\bf annihilation operation of a non-self-connected separatrix} if the separatrices of $\alpha$ intersecting the trivial flow box $B_w$ do not intersect any separatrix of $\omega$. 

\subsubsection{Whitehead moves}

\begin{figure}
\begin{center}
\includegraphics[scale=0.6]{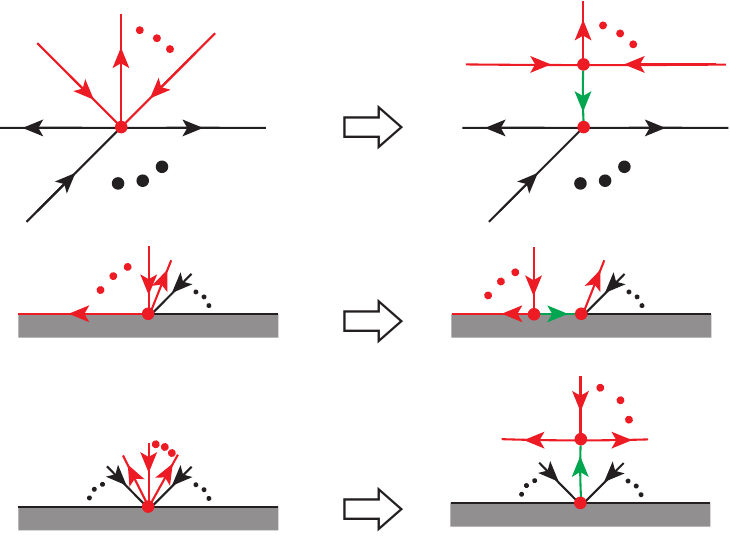}
\end{center}
\caption{Whitehead moves: splittings of multi-saddles}
\label{splitting}
\end{figure}

Every collapsing of a heteroclinic separatrix and the inverse operation as in Figure~\ref{splitting} is called the {\bf Whitehead move}. 
We observe the following isotopy realizing Whitehead moves. 

\begin{lemma}\label{lem:example01}
For any $q,s \in \Z_{>0}$, a smooth $1$-parameter family $\{X_{t} \mid t \in [0,1] \}$ on the complex plane $\mathbb{C}$ defined by 
\[
X_t := (\mathrm{Re}(\bar{z}^{q-s}(\bar{z}-t)^{s}), \mathrm{Im}(\bar{z}^{q-s}(\bar{z}-t)^{s}))
\]
generates a $1$-parameter family $v_{X_t}$ of smooth gradient flows satisfying the following statements for any $t \in (0,1]$:
\\
{\rm(1)} $\mathrm{codim}(v_{X_0}) - \mathrm{codim}(v_{X_t}) = 1$. 
\\
{\rm(2)} $\mathrm{codim}_{\mathrm{m}}(v_{X_0}) - \mathrm{codim}_{\mathrm{m}}(v_{X_t}) = 2$. 
\\
{\rm(3)} $\mathrm{codim}_{\mathrm{h}}(v_{X_0}) - \mathrm{codim}_{\mathrm{h}}(v_{X_t}) = -1$. 
\end{lemma}

The previous lemma implies the following ``$C^r$-density'' up to topological equivalent. 

\begin{lemma}\label{lem:WH}
For any flow $v \in \mathcal{G}^r_{k, > -1}(S)$ with a multi-saddle $x$ whose codimension is nonzero, there is a flow $w \in \mathcal{G}^r_{\kappa, > -1}(S)$ which is topological equivalent to $v$ via a topological conjugacy $h \colon S \to S$ such that any \nbd $\mathcal{U} \subset \mathcal{G}^r_{k, > -1}(S)$ of $w$ contains a flow $w' \in \mathcal{G}^r_{k, > -1}(S)$ with $\mathrm{codim}(w) - \mathrm{codim}(w') = 1 = \mathrm{codim}_{\mathrm{m}}(w') - \mathrm{codim}_{\mathrm{m}}(w)$. 
\end{lemma}

\begin{proof}
Let $v \in \mathcal{G}^r_{\kappa, > -1}(S)$ be a flow with a multi-saddle $x$ whose codimension is nonzero. 

Suppose that $x \notin \partial S$. 
By Gutierrez's smoothing theorem~\cite{gutierrez1978structural}, take a flow $w \in \mathcal{G}^r_{k, > -1}(S)$ which is generated by a smooth vector field $X_w$ and is topological equivalent to $v$ via a topological conjugacy $h \colon S \to S$ such that a small \nbd $V$ of $h(x)$ can be identified with a small \nbd of the origin in $\mathbb{C}$ where the orbits are generated by the vector field $X_0 = (\mathrm{Re}(\bar{z}^{q}), \mathrm{Im}(\bar{z}^{q}))$ as in Lemma~\ref{lem:example01}. 
Choose an open subset $U \Subset V$ whose boundary consists of finitely many outer tangencies and finitely many open transverse arcs, where $U \Subset V$ means $\overline{U} \subset \mathop{\mathrm{int}} V$. 
Denote by $F_U$ (resp. $F_V$) the intersection of $\partial U$ (resp. $\partial V$) and separatrices of $h(x)$. 
Fix any \nbd $\mathcal{U} \subset \mathcal{G}^r_{k, > -1}(S)$ of $w$. 
Take a bump function $\varphi \colon S \to [0,1]$ whose support is contained in $U$ and which is rotationally symmetric and weakly monotonically decreasing with respect to the radius of polar coordinates in $\mathbb{C}$ such that $\varphi^{-1}(1)$ is a \nbd of $h(x)$. 
Then $\varphi X_0 = \varphi X_w$. 
Moreover, a vector field $X_{w,t} := (1 - \varphi) X_w + \varphi X_t$ generates a smooth $1$-parameter $\{ X_{w,t} \mid t \in [0,1]  \}$ such that the $1$-parameter $\{ w_t \mid t \in [0,1]  \}$ of flows generated by $\{ X_{w,t} \mid t \in [0,1]  \}$ is contained in $\mathcal{G}^r_{k, > -1}(S)$, and that $X_{w,0} = (1 - \varphi) X_w + \varphi X_0 = (1 - \varphi) X_w + \varphi X_w = X_w$. 
Then $V$ contains exactly two singular points $x_{t+}$ and $x_{t-}$ of $w_t$ for any $t>0$. 
For any point $y \in F_U$, consider a trivial flow box $B_{y} \Subset V$ which is a \nbd of $y$. 
Let $H_{w,t}$ be the height function whose gradient vector field is $X_{w,t}$. 
For any small $t >0$, perturbing $H_{w,t}$ in $\bigcup_{y \in F_U} B_{y}$, we can obtain the resulting height function $H_{w',t}$ whose gradient vector field $X_{w',t}$ generates a flow $w'_t$ such that two singular points $x_{t+}$ and $x_{t-}$ are exactly singular points of $w'_t$ in $V$ and are contained in the same multi-saddle connection of $w'_t$. 
Moreover, we may assume that $\lim_{t \to 0} \|X_{w,t} - X_{w',t} \|_{C^r} = 0$. 
Since $\{ X_{w,t} \mid t \in [0,1]  \}$ is smooth, there is a small number $\varepsilon \in (0,1]$ such that $w'_{\varepsilon} \in \mathcal{U}$, $\mathrm{codim}_{\mathrm{m}}(w) - \mathrm{codim}_{\mathrm{m}}(w'_\varepsilon) = 2$, and $\mathrm{codim}_{\mathrm{h}}(w) - \mathrm{codim}_{\mathrm{h}}(w'_\varepsilon) = -1$. 
Therefore, $\mathrm{codim}(w) - \mathrm{codim}(w'_\varepsilon) = 1$. 
This means that $w' := w'_{\varepsilon}$ is desired. 

Suppose that $x \in \partial S$. 
Replacing $\mathbb{C}$ with $\{z \in C \mid \mathrm{Im}(z) \geq 0 \}$ in the previous argument, we can obtain a flow $w'_{\varepsilon} \in \mathcal{U}$ such that $\mathrm{codim}_{\mathrm{m}}(w) - \mathrm{codim}_{\mathrm{m}}(w'_\varepsilon) = 1$ and $\mathrm{codim}_{\mathrm{h}}(w) - \mathrm{codim}_{\mathrm{h}}(w'_\varepsilon) = 0$ if $x \in \partial S$. 
Therefore, $\mathrm{codim}(w) - \mathrm{codim}(w'_\varepsilon) = 1$. 
This means that $w' := w'_{\varepsilon}$ is desired. 
\end{proof}

\subsubsection{Associated neighborhoods of the multi-saddle connections}

We have the local stability of basins of repellers. 

\begin{lemma}\label{lem:inv_sink/source}
Let $v$ be a flow on a compact surface $S$, $x \in S$ a source {\rm(resp.} sink{\rm)}, and $B$ and $B'$ closed disks with $x \in \mathop{\mathrm{int}} B'$ and $B' \Subset B$. 
Suppose that there is a positive  {\rm(resp.} negative{\rm)} number $T \neq 0$ such that $v(T,B) \cup v([0,T],B') \Subset B$, and $v(T,B') \Subset B'$.
There is a $C^0$-\nbd $\mathcal{U}$ of $v$ satisfying the following statements:  
\\
{\rm(1)} $1 = \mathrm{ind}_v(x) = \sum_{y \in B \cap \mathop{\mathrm{Sing}}(w)} \mathrm{ind}_w(y)$ for any flow $w \in \mathcal{U}$ with $ \vert B \cap \mathop{\mathrm{Sing}}(w) \vert  < \infty$.
\\
{\rm(2)} $w(T,B') \Subset w(T,B) \Subset B'$ and $w(\R_{\geq 0},B') \Subset B$ {\rm(resp.} $w(\R_{\leq 0},B') \Subset B$ {\rm)} for any flow $w \in \mathcal{U}$. 
\\
{\rm(3)} For any flow $w \in \mathcal{U}$, there is a closed transversal $\mu$ of $w$ in $\mathop{\mathrm{int}} B$  such that $\mathop{\mathrm{Sing}}(w) \cap B = \mathop{\mathrm{Sing}}(w) \cap D_\mu$, where $D_\mu \subset \mathop{\mathrm{int}}B$ is the closed disk whose boundary is $\mu$. 
\end{lemma}

Since the previous lemma naturally holds in the case of smooth vector fields and the proof of the previous lemma is straightforward but technical, the proof is stated in Appendix~\ref{sec:loc_stb}.
We have the following key description of the associated neighborhood of the union of the multi-saddle connection diagram and positively indexed singular points.

\begin{figure}
\begin{center}
\includegraphics[scale=0.8]{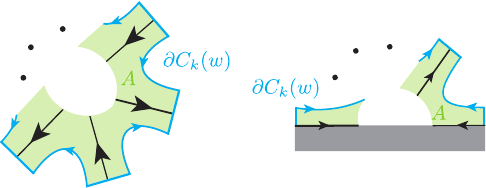}
\end{center}
\caption{A small collar of the boundary of a \nbd of multi-saddles outside of the boundary and one on the boundary.}
\label{multi-saddle_nbd02}
\end{figure}

\begin{figure}
\begin{center}
\includegraphics[scale=0.9]{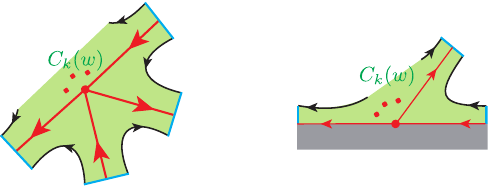}
\end{center}
\caption{A connected component $C_k(w)$ of the set difference $W(w) \setminus \bigsqcup_{i} \gamma_{i,w}$ containing a singular point of $v$ outside of the boundary and one on the boundary.}
\label{multi-saddle_nbd00}
\end{figure}

\begin{proposition}\label{lem:no_merge_msc}
Let $S$ be a compact connected surface. 
For any $v \in \mathcal{G}^r_{\kappa, > -1}(S)$ with $D(v) \neq \emptyset$, there are disjoint open disks $B_1, B_2, \ldots , B_l$ whose union is a \nbd of the union of positively indexed singular points and there is a closed \nbd $W_v$ of the multi-saddle connection diagram $D(v)$,  
 and there is a \nbd $\mathcal{U}$ in $\mathcal{G}^r_{\kappa, > -1}(S)$ of $v$ satisfying the following properties for any element $w$ of $\mathcal{U}$: 
\\
{\rm (1)} The flow $w$ has same type of positively indexed singular points of $v$ with respect to $B_1, B_2, \ldots , B_l$. 
\\
{\rm(2)} The \nbd $W_v$ is a deformation retract of $D(v)$ and  contains an open invariant \nbd $W(w)$ of the multi-saddle connection diagram $D(w)$ of $w$ which is a deformation retract of $D(v)$ such that the set difference $W(w) - D(w)$ consists of finitely many invariant open trivial flow boxes of $w$. 
\\
{\rm(3)} Any connected components of the complement $S - (W(w) \sqcup \mathop{\mathrm{Sing}}_+(w))$ are invariant open trivial flow boxes of $w$. 
\\
{\rm(4)} $N(w) \geq N(v)$. 
\\
{\rm(5)} There are loops $\gamma_{i,w} \subset W_v$ which are transverse to $W(w)$ and intersect separatrics of $w$ such that any connected component $C_k(w)$ of the set difference $W(w) \setminus \bigsqcup_{i} \gamma_{i,w}$ is either a trivial flow box or a disk in which the restriction $w' \vert _{A}$ of $w'$ to some collar $A$ of $\partial C_k(w)$ in $C_k(w)$ is locally topologically equivalent to the restriction to a small collar of the boundary of an isolated \nbd of either a $n$-saddle or $\partial$-$n/2$-saddle of the vector field $Z = (Z_x, Z_y) := ( r \cos (n \theta), r \cos (n \theta))$ as in Figure~\ref{multi-saddle_nbd02}, where $(r,\theta)$ is the polar coordinate system. 
\\
{\rm(6)} For any $C_k(w)$ containing singular points, the connected component $C_k(w)$ contains exactly one singular point $x_{v,k}$ of $v$ as in Figure~\ref{multi-saddle_nbd00}, which is a multi-saddle, and finitely many singular points $x_{w,k,1}, \ldots , x_{w,k,m_k}$ of $w$, which are multi-saddles, such that if $x_{v,k} \notin \partial S$ then 
\[
\mathrm{codim}_{\mathrm{m},v}(x_{v,k}) + 2(1- m_k) = \sum_{s= 1}^{m_k} \mathrm{codim}_{\mathrm{m},w}(x_{w,k,s}), 
\]
and if $x_{v,k} \in \partial S$ then 
\[
\mathrm{codim}_{\mathrm{m},v}(x_{v,k})  + 1+ m'_k - 2m_k = \sum_{s= 1}^{m_k} \mathrm{codim}_{\mathrm{m},w}(x_{w,k,s}), 
\]
where $m'_k$ is the number of multi-saddles of $w$ on the boundary $\partial S$ in the multi-saddles $x_{w,k,1}, \ldots , x_{w,k,m_k}$ of $w$ in $C_k(w)$. 
\\
{\rm(7)} Any restriction $w \vert _{C_k(w)}$ on $C_k(w)$ of $w$ can be obtain from the restriction $v \vert _{C_k(w)}$ by Whitehead moves and annihilation operations of non-self-connected separatrices finitely many times. 
\end{proposition}

\begin{proof}
Let $x_1, x_2, \ldots , x_{k'}$ be the positively indexed repelling singular points of the flow $v$, $x_{k'+1}, x_{k'+2}, \ldots , x_{k''}$  the positively indexed attracting singular points, and $x_{k''+1}, x_{k''+2}, \ldots , x_l$ the multi-saddles of $v$.  
Fix any small positive number $\delta > 0$. 
Since $\delta > 0$ is small, the closures of the open disks $B_i := B_\delta(x_i)$ are disjoint. 
Take any small disjoint open disks $U_1, U_2, \ldots , U_{l}$ with $B_{\delta/4}(x_i) \subset U_i \setminus B_{\delta/4}(\partial U_i)$ and $B_{\delta/4}(U_i) \subset B_\delta(x_i) = B_i$ such that $\partial U_1, \partial U_2, \ldots , \partial U_{k''}$ are closed transversals, and that $\partial U_{k''+1}$, $\partial U_{k''+2}, \ldots , \partial U_{l}$ are loops which are transverse except finitely many outer tangencies.  
By Lemma~\ref{lem:fin_comb}, there is an open invariant \nbd $W_v$ of $D(v)$ which is a deformation retract of $D(v)$ such that the set difference $W_v - D(v)$ consists of finitely many invariant open trivial flow boxes of $v$. 
Then the union $V_v := W_v \sqcup \bigsqcup_{i=1}^{k''} B_i$ is a closed \nbd of the union $\mathop{\mathrm{Sing}}(v) \cup D(v)$ which is a deformation retract of $\mathop{\mathrm{Sing}}(v) \cup D(v)$. 
For any connected $C_j$ component of $W_v - D(v)$, choose a point $c_j \in C_j \cap \bigsqcup_{i=1}^{k' }U_i$. 
Then there are positive numbers $T'_j$ with $v(T'_j, c_j) \in \bigsqcup_{i=k'+1}^{k''}U_i$. 
For any $c_j$, there is exactly one index $i'(c_j)$ with $c_j \in U_{i'(c_j)}$ and there is exactly one index $i(c_j)$ with $v(T'_j, c_j) \in U_{i(c_j)}$. 
Fix any small positive number $\varepsilon > 0$ such that $B_{\varepsilon}(O_v(c_j)) \cap \mathop{\mathrm{Sing}}_-(v) = \emptyset$ for any $c_j$.

\begin{claim}\label{claim:001}
There is a neighborhood $\mathcal{U}_4 \subset \mathcal{G}^r_{> -1}(S)$ of $v$ such that the following statements hold for any element $w$ of $\mathcal{U}_4$: 
\\
{\rm (1)} The flow $w$ has same type of positively indexed singular points of $v$ with respect to $B_1, B_2, \ldots , B_l$. 
In particular, the flow $w$ contains no positively indexed singular points of $w$ on $S - \bigsqcup_{i=1}^k B_i$. 
\\
{\rm (2)} There are positive number $T_1, \ldots ,T_{k''}$ such that $w(\R_{\leq 0}, U_i) \Subset B_i$ and $w(-T_i, B_i) \Subset B_i$ {\rm(resp.} $w(\R_{\geq 0}, U_i) \Subset B_i$ and $w(T_i, B_i) \Subset B_i${\rm)} for any $i \in \{ 1, 2, \ldots , k' \}$ {\rm(resp.} $i \in \{ k'+1, k'+2, \ldots , k'' \}${\rm)}.
\\
{\rm (3)} There are loops $\gamma_{i,w} \subset B_i$ each of which is transverse at all but finitely many outer tangencies such that the number of the outer tangencies on $\partial U_i$ of $v$ equals one of the outer tangencies on $\gamma_{i,w}$ of $w$. 
\\
{\rm (4)} For any $j$, we obtain $O_w(c_j) \subset (W_v \cap B_{{\varepsilon}}(O_v(c_j))) \cup \bigsqcup_{i=1}^k B_i$, $\alpha_w(c_j) = x_{w,i'(c_j)}$, and $\omega_w(c_j) = x_{w,i(c_j)}$, where $x_{w,\lambda}$ is the continuation of the positively indexed singular point $x_\lambda$.
\end{claim}
\begin{proof}[Proof of Claim~\ref{claim:001}]
Lemma~\ref{lem:inv_index02} implies that there is a $C^0$-\nbd $\mathcal{U}_0 \subset \mathcal{G}^r_{> -1}(S)$ of $v$ such that $\bigcup_{w \in \mathcal{U}_0} \mathop{\mathrm{Sing}}(w) \subset B_{\delta/4} (\mathop{\mathrm{Sing}}(v)) \subset \bigsqcup_{i=1}^l U_i \setminus B_{\delta/4}(\partial U_i)$. 
From Lemma~\ref{lem:perturbation}, there is a $C^0$-\nbd $\mathcal{U}_1 \subseteq \mathcal{U}_0$ of $v$ each element of which has same type of positively indexed singular points of $v$  with respect to $B_1, B_2, \ldots , B_l$. 
Therefore, assertion {\rm (1)} holds. 
%
By Lemma~\ref{lem:inv_sink/source}, there is a \nbd $\mathcal{U}_2 \subset \mathcal{U}_1$ of $v$ such that there are nonzero number $T_1, \ldots ,T_k$ such that $w(\R_{\leq 0}, U_i) \Subset B_i$ and $w(T_i, B_i) \Subset B_i$ for any flow $w \in \mathcal{U}_2$. 
 This means that assertion {\rm (2)} holds. 
Applying Lemma~\ref{lem:inv_index} to $\partial U_i$, there is a \nbd $\mathcal{U}_3 \subset \mathcal{U}_2$ of $v$, for any $w \in \mathcal{U}_3$, there are loops $\gamma_{i,w} \subset B_{\delta/4}(\partial U_i) \subset B_i$ each of which is transverse at all but finitely many outer tangencies such that the number of the outer tangencies on $\partial U_i$ of $v$ equals one of the outer tangencies on $\gamma_{i,w}$ of $w$. 
This means that assertion {\rm (3)} holds. 
By definition, the intersection $\bigcap_j \mathcal{U}(\{T'_j\} \times \{c_j\}, B_{i(c_j)}) = \bigcap_j \{ w' \in \mathcal{G}^r_{> -1}(S) \mid w'(\{T'_j\} \times \{c_j\}) \subset B_{i(c_j)} \}$ of open subsets is an open \nbd of $v$. 
For any orbit arc $v([0,T'_j], c_j)$, fix the connected \nbd $W(c_j) \subset (W_v - D(v)) \cap B_{{\varepsilon}}(O_v(c_j))$ containing $v([0,T'_j], c_j)$. 
Then the intersection 
\[
\bigcap_j \mathcal{U}([0,T'_j] \times \{c_j\}, W(c_j)) = \bigcap_j \{ w' \in \mathcal{G}^r_{> -1}(S) \mid w'([0,T'_j] \times \{c_j\}) \subset W(c_j) \}
\]
of open subsets is an open \nbd of $v$. 
Therefore, the intersection 
\[
\mathcal{U}_4 := \mathcal{U}_3 \cap \bigcap_j \mathcal{U}([0,T'_j] \times \{c_j\}, W(c_j)) \cap \bigcap_j \mathcal{U}(\{T'_j\} \times \{c_j\}, B_{i(c_j)})
\]
is an open \nbd of $v$ such that $w(\{T'_j\} \times \{c_j\}) \subset B_{i(c_j)}$ and $w([0,T'_j] \times \{c_j\}) \subset W(c_j) \subset W_v \subset V_v$ for any $w \in \mathcal{U}_4$. 
By assertions {\rm (1)} and {\rm (2)}, we have that $w(\R_{\leq 0}, c_j) \subset w(\R_{\leq 0}, U_{i(c_j)}) \Subset B_{i(c_j)} \subset V_v$ and $w(\R_{\geq T'_j}, c_j) = w(\R_{\geq 0}, w(T'_j,c_j)) \subset w(\R_{\geq 0}, U_{i'(c_j)}) \Subset B_{i'(c_j)} \subset V_v$. 
Therefore, we have that $\alpha_w(c_j) = x_{w,i'(c_j)}$, $\omega_w(c_j) = x_{w,i(c_j)}$, and 
\[
\begin{split}
O_w(c_j) &= w(\R_{\leq 0}, c_j) \cup v([0,T'_j], c_j) \cup w(\R_{\geq T'_j}, c_j)
\\
& \subset B_{i(c_j)} \cup (W_v \cap B_{{\varepsilon}}(O_v(c_j)) \cup B_{i'(c_j)} \subseteq (W_v \cap B_{{\varepsilon}}(O_v(c_j)) \cup \bigsqcup_{i=1}^k B_i
\end{split}
\]
 for any $j$. 
\end{proof}

\begin{claim}\label{claim:001++}
For any element $w \in \mathcal{U}_4$, there is an open $w$-invariant \nbd $W'(w) \subseteq V_v$ of the multi-saddle connection diagram $D(w)$ of $w$ satisfying the following statements: 
\\
{\rm (1)} The \nbd $W'(w)$ is a deformation retract of $D(v)$ each of whose connected components contains exactly one multi-saddle connection of $v$. 
\\
{\rm (2)} The set difference $W'(w) \setminus D(w)$ consists of finitely many invariant open trivial flow boxes of $w$.  
\end{claim}
\begin{proof}[Proof of Claim~\ref{claim:001++}]
Fix any $w \in \mathcal{U}_4$. 
By Lemma~\ref{lem:fin_comb} and Claim~\ref{claim:001}, the finite union of the connected component $W'(w)$ of $W_v - \bigsqcup_{j} O_w(c_j)$ intersecting $D(v)$ is desired. 
\end{proof}

\begin{claim}\label{claim:001+}
For any element $w \in \mathcal{U}_4$, we obtain $N(w) \geq N(v)$. 
\end{claim}
\begin{proof}[Proof of Claim~\ref{claim:001+}]
By Claim~\ref{claim:001++}, any connected component of $W'(v)$
 contains exactly one multi-saddle connection of $v$. 
On the other hand, any connected component of $W'(w)$ 
 contains at least one multi-saddle connections of $w$. 
\end{proof}

\begin{claim}\label{claim:002}
There is an open invariant \nbd $W(w) \subset W'(w)$ of the multi-saddle connection diagram $D(w)$ of $w$: 
\\
{\rm(1)} The \nbd $W(w)$ is a deformation retract of $D(v)$ such that $W(w) - D(w)$ consists of finitely many invariant open trivial flow boxes of $w$. 
\\
{\rm(2)} Any connected component of the complement $S - (W(w) \sqcup \mathop{\mathrm{Sing}}_+(w))$ is an invariant open trivial flow box of $w$. 
\\
{\rm(3)} The loops $\gamma_{i,w} \subset B_i$ as in Claim~\ref{claim:001} are transverse to $W(w)$ and intersect separatrics of $w$ such that any connected component $C_k(w)$ containing multi-saddles of the set difference $W(w) \setminus \bigsqcup_{i} \gamma_{i,w}$ is a disk in which the restriction $w' \vert _{A}$ of $w'$ to some collar $A$ of $\partial C_k(w)$ in $C_k(w)$ is locally topologically equivalent to the restriction to a small collar of the boundary of an isolated \nbd of either a $n$-saddle or $\partial$-$n/2$-saddle as in Figure~\ref{multi-saddle_nbd02}. 
\end{claim}
\begin{proof}[Proof of Claim~\ref{claim:002}]
Fix any $w \in \mathcal{U}_4$ and loops $\gamma_{1,w}, \ldots , \gamma_{l,w}$  as in Claim~\ref{claim:001}. 
By perturbing the loops $\gamma_{1,w}, \ldots , \gamma_{l,w}$, we may assume that any tangencies intersect no semi-multi-saddle separatrices. 
Choose small closed arcs $I_1, \ldots , I_m \subset \bigsqcup_{i=1}^l \gamma_{i,w}$ each of whose interiors contains exactly one tangency and intersects no semi-multi-saddle separatrices. 
By construction, the union $W(w)$ of the connected components of $W'(w) - \bigcup_{j=1}^m w(I_j)$ intersecting $D(w)$ is desired. 
\end{proof}

\begin{claim}\label{claim:004}
Assertion {\rm (6)} holds. 
\end{claim}
\begin{proof}[Proof of Claim~\ref{claim:004}]

By Claim~\ref{claim:001++}, any connected component $C_k(w)$ contains exactly one singular point $x_{v,k}$ of $v$, which is a multi-saddle of $v$. 
Claim~\ref{claim:001}(1) implies that any connected component $C_k(w)$ contains only multi-saddles $x_{w,k,1}, \ldots , x_{w,k,m_k}$. 
The Poincar{\'e}-Hopf theorem implies that the sum of indices of singular points of $w$ in $C_k(w)$ corresponds to one of $v$ in $C_k(w)$.  
Since the index of $k$-saddle (resp. $\partial$-$k/2$-saddle) is $-k$ (resp. $-k/2$), we obtain the following equality: 
\[
\mathop{\mathrm{ind}_v}(x_{v,k}) = \sum_{s= 1}^{m_k} \mathop{\mathrm{ind}_w}(x_{w,k,s})
\] 
Because the codimension of a $k$-saddle $y$ of $w$ ($k \geq 1$) is $\mathrm{codim}_{\mathrm{m},w}(y) = 2(k-1) = -2(1+ \mathrm{ind}(y))$, if $x_{v,k} \notin \partial S$, then any multi-saddles $x_{w,k,1}, \ldots , x_{w,k,m_k}$ are outside of the boundary $\partial S$ and so the following equality holds: 
\[
\begin{split}
\mathrm{codim}_{\mathrm{m},v}(x_{v,k}) &= -2(1+ \mathrm{ind}(x_{v,k})) 
\\
&= -2 \left(1+ \sum_{s= 1}^{m_k} \mathop{\mathrm{ind}_w}(x_{w,k,s}) \right)
\\
&= -2(1- m_k) -2 m_k - \sum_{s= 1}^{m_k} \mathop{\mathrm{ind}_w}(x_{w,k,s})
\\
&= -2(1- m_k) + \sum_{s= 1}^{m_k} -2(1+  \mathop{\mathrm{ind}_w}(x_{w,k,s}))
\\
&= -2(1- m_k) + \sum_{s= 1}^{m_k} \mathrm{codim}_{\mathrm{m},w}(x_{w,k,s})
\end{split}
\]
This means that assertion (6) holds if $x_{v,k} \notin \partial S$. 

Thus, we may assume that $x_{v,k} \in \partial S$. 
Then at least one of multi-saddles $x_{w,k,s}$ is on the boundary $\partial S$. 
By renumbering if necessary, we may assume that the multi-saddles $x_{w,k,1}, \ldots , x_{w,k,m'_k}$ are  on the boundary $\partial S$ and the multi-saddles $x_{w,k,m_{k'+1}}, \ldots , x_{w,k,m_k}$ are outside of the boundary $\partial S$. 
Since the codimension of a $\partial$-$k$-saddle $y$ of $w$ ($k \geq 1/2$) is $\mathrm{codim}_{\mathrm{m},w}(y) = 2k-1 = -2(1/2 + \mathrm{ind}(y))$, the following equality holds: 
\[
\begin{split}
&\mathrm{codim}_{\mathrm{m},v}(x_{v,k}) + 1+ m'_k - 2m_k
\\
= &\mathrm{codim}_{\mathrm{m},v}(x_{v,k}) + (1-m'_k) - 2(m_k - m'_k)
\\
= &-2(1/2+ \mathrm{ind}(x_{v,k})) + (1-m'_k) - 2(m_k - m'_k)
\\
= &-2\left(1/2+ \sum_{s= 1}^{m_k} \mathop{\mathrm{ind}_w}(x_{w,k,s})\right) + (1-m'_k) - 2(m_k - m'_k)
\\
= &-2(m'_k/2 + (m_k - m'_k)) -2 \sum_{s= 1}^{m_k} \mathop{\mathrm{ind}_w}(x_{w,k,s})
\\
= &\sum_{s= 1}^{m'_k} -2(1/2+  \mathop{\mathrm{ind}_w}(x_{w,k,s})) + \sum_{s= m'_k + 1}^{m_k} -2(1+  \mathop{\mathrm{ind}_w}(x_{w,k,s}))
\\
= &\sum_{s= 1}^{m_k} \mathrm{codim}_{\mathrm{m},w}(x_{w,k,s})
\end{split}
\]
This means that assertion (6) holds if $x_{v,k} \in \partial S$. 
\end{proof}

By the inverse operations of annihilation operations of non-self-connected separatrices and Whitehead moves to the restriction $w \vert _{C_k(w)}$ on $C_k(w)$ of $w$ which preserves near $\partial C_k(w)$ finitely many times, we can obtain a multi-saddle on $C_k(w)$ such that the restriction near the multi-saddle to the resulting flow is locally topological equivalent to the  restriction $v \vert _{C_k(w)}$. 
The inverse operation implies assertion {\rm(7)}. 
%
%
%
%
%
%
This completes the proof. 
\end{proof}

\section{Combinatorial structure of the space of gradient flows}


For any $q \in \mathbb{Z}_{\geq 0}$, denote by $\bm{\mathcal{G}^r_{k_-/2, k_+/2, q}(S)} \subseteq \mathcal{G}^r_{k_-/2, k_+/2, > -1}(S)$ the subspace of the gradient $C^r$-flows on a compact surface $S$ whose codimensions are $q$. 
Put $\bm{\mathcal{G}^r_{k_-/2, k_+/2, > q}(S)} := \bigcup_{q' > q}\mathcal{G}^r_{k_-/2, k_+/2, q'}(S)$. 
Notice that the subset $\mathcal{G}^r_{k_-/2, k_+/2, > q}(S)$ is the subspace of the gradient flows whose codimensions are more than $q$. 


For any $k_-, k_+ \in \mathbb{Z}_{\geq 0}$ with $k := k_-/2 + k_+/2$, denote by 
\[
\bm{\mathcal{G}^r_{k_-/2, k_+/2, > -1}(S)} \subseteq \mathcal{G}^r_{k, > -1}(S)
\]
the subset of quasi-regular gradient flows whose sums of indices of sinks and $\partial$-sinks (resp. sources and $\partial$-sources) are $k_-/2$ (resp. $k_+/2$). 
Notice that the subspace $\mathcal{G}^r_{k_-/2, k_+/2, > -1}(S)$ is a connected component of $\mathcal{G}^r_{k, > -1}(S)$.


\subsection{Hierarchical structure of gradient flows}

In this section, we have the following hierarchical structure. 


\begin{theorem}\label{prop:codimension}
The following statements hold for any $r \in \mathbb{Z}_{\geq 0} \sqcup \{ \infty \}$, any integers $k_-, k_+ \in \mathbb{Z}_{\geq 0}$ and any $q \in \Z_{\geq -1}$: 
\\
{\rm(1)} The subspace $\mathcal{G}^r_{k_-/2, k_+/2, q+1}(S)$ is open in the space $\mathcal{G}^r_{k_-/2, k_+/2, > q}(S)$
\\
{\rm(2)} The subspace $\mathcal{G}^r_{k_-/2, k_+/2, q+1}(S)$ is $C^0$-dense in the space $\mathcal{G}^r_{k_-/2, k_+/2, > q}(S)$. 
\\
{\rm(3)} The subspace $\mathcal{G}^r_{k_-/2, k_+/2, q+1}(S)$ consists of $C^r$-structurally stable flows in the space $\mathcal{G}^r_{k_-/2, k_+/2, > q}(S)$. 
\end{theorem}

To show the previous theorem, we state the following two lemmas. 
First, we have the following observation. 

\begin{lemma}\label{lem:codimension_open}
For any $r \in \mathbb{Z}_{> 0} \sqcup \{ \infty \}$, any integer $k \in \Z_{\geq 0}$ and any $q \in \Z_{\geq -1}$, the subspace $\mathcal{G}^r_{k, q+1}(S)$ is open in the space $\mathcal{G}^r_{k, > q}(S)$ and consists of $C^r$-structurally stable flows. 
Moreover, for any $k_-, k_+ \in \mathbb{Z}_{\geq 0}$ and any $s \in \mathbb{Z}_{> 0} \sqcup \{ \infty \}$, the subspace $\mathcal{G}^r_{k_-, k_+, q+1}(S)$ is open in the space $\mathcal{G}^r_{k_-, k_+, > q}(S)$ and consists of $C^s$-structurally stable flows in the space $\mathcal{G}^r_{k_-, k_+, > q}(S)$. 
\end{lemma}

\begin{proof}
Fix any flow $v \in \mathcal{G}^r_{k, > q}(S)$. 
Denote by $k_-, k_+ \in \mathbb{Z}_{\geq 0}$ integers with $v \in \mathcal{G}^r_{k_-/2, k_+/2, -1}(S)$. 
Choose any small \nbd $\mathcal{U}$ of $v$ as in Proposition~\ref{lem:no_merge_msc}. 
By Proposition~\ref{lem:no_merge_msc}, the \nbd $\mathcal{U}$ is contained in $\mathcal{G}^r_{k_-/2, k_+/2, -1}(S)$. 
 
\begin{claim}\label{claim:05}
$\mathrm{codim}(v) - \mathrm{codim}(w) \geq 0$ for any flow $w \in \mathcal{U}$, and that the equality holds only that $w$ is topologically equivalent to $v$. 
\end{claim}

\begin{proof}[Proof of Claim~\ref{claim:05}]
 Fix any flow $w \in \mathcal{U}$. 
Since the closure of any multi-saddle connection is a tree (i.e. the closure of the multi-saddle connection diagram is a forest), the difference $\mathrm{codim}_{\mathrm{h}}(w) - \mathrm{codim}_{\mathrm{h}}(v)$ is less than or equals the difference of the numbers $w$ and $v$ of multi-saddles outside of $\partial S$ and so the following inequality holds: 
\[
\mathrm{codim}_{\mathrm{h}}(w) - \mathrm{codim}_{\mathrm{h}}(v) \leq \sum_{x_{v,k} \notin \partial S} (m_k -1) + \sum_{x_{v,k} \in \partial S} (m_k - m'_k)
\]
Here, $x_{v,k}$ is the singular point of $v$ as in Lemma~\ref{lem:WH} (6), $m_k$ is the number of multi-saddles in the \nbd $C_k(w)$ of $x_{v,k}$, and $m'_k$ is the number of multi-saddles of $w$ on the boundary $\partial S$ in the multi-saddles $x_{w,k,1}, \ldots , x_{w,k,m_k}$ of $w$ in $C_k(w)$. 
On the other hand, the following equality holds: 
\[
\mathrm{codim}_{\mathrm{m}}(w) - \mathrm{codim}_{\mathrm{m}}(v) = \sum_{x_{v,k} \notin \partial S} 2(1- m_k) +  \sum_{x_{v,k} \in \partial S} 1+ m'_k - 2m_k
\]
Therefore, we have the following positivity: 
\[
\begin{split}
&\mathrm{codim}(v) - \mathrm{codim}(w) 
\\
= &(\mathrm{codim}_{\mathrm{m}}(v) + \mathrm{codim}_{\mathrm{h}}(v)) - (\mathrm{codim}_{\mathrm{m}}(w) + \mathrm{codim}_{\mathrm{h}}(w))
\\
= &(\mathrm{codim}_{\mathrm{m}}(v) - \mathrm{codim}_{\mathrm{m}}(w)) - (\mathrm{codim}_{\mathrm{h}}(v) - \mathrm{codim}_{\mathrm{h}}(w))
\\
\geq &\sum_{x_{v,k} \notin \partial S} (- 2(1- m_k)-(m_k -1)) +  \sum_{x_{v,k} \in \partial S} (- (1+ m'_k - 2m_k) -(m_k - m'_k))
\\ 
= &\sum_{x_{v,k} \notin \partial S} (- 1 + m_k) +  \sum_{x_{v,k} \in \partial S} (- 1 + m_k)
\\ 
= &\sum_{x_{v,k}} (- 1 + m_k) \geq 0
\end{split}
\]
Moreover, the equality $\mathrm{codim}(v) = \mathrm{codim}(w)$ holds if and only if any multi-saddles and multi-saddle separatrices are preserved. 
By Proposition~\ref{lem:no_merge_msc}, this means that $\mathrm{codim}(v) = \mathrm{codim}(w)$ holds if and only if $v$ and $w$ are topologically equivalent. 
\end{proof}

The previous claim implies that, for any $q \in \Z_{\geq -1}$,  the subspace $\mathcal{G}^r_{k, q+1}(S)$ is open in the space $\mathcal{G}^r_{k, > q}(S)$ and the subspace $\mathcal{G}^r_{k_-, k_+, q+1}(S)$ is open and $C^s$-structurally stable in the space $\mathcal{G}^r_{k_-, k_+, > q}(S)$ for any $s \in \mathbb{Z}_{> 0} \sqcup \{ \infty \}$. 
\end{proof}

\begin{figure}
\begin{center}
\includegraphics[scale=0.625]{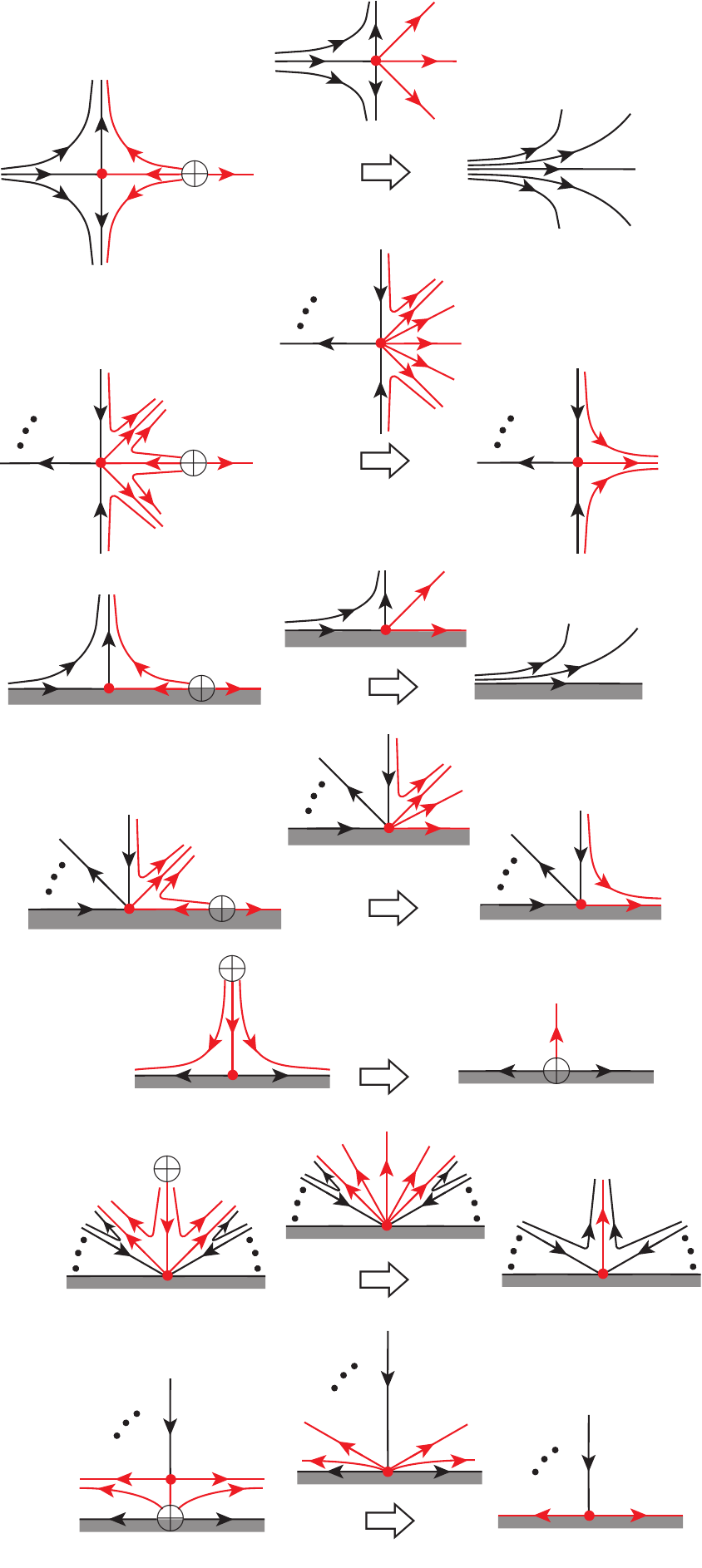}
\end{center}
\caption{Forbidden merges of multi-saddles and attracting singular points in $\mathcal{G}^r_{k_-/2, k_+/2, >-1}(S)$}
\label{fig:merge_attracting}
\end{figure}

The previous lemma implies that any merges as in Figure~\ref{fig:merge_attracting} are forbidden.  
For a class $\mathcal{C} \subseteq \mathcal{G}^r_{k, -1}(S)$, denote by $[\mathcal{C}]$ the equivalent class of $\mathcal{C}$ by the topological equivalence. 
We have the following denseness. 

\begin{lemma}\label{lem:codimension_density}
The following statements hold for any $r \in \mathbb{Z}_{> 0} \sqcup \{ \infty \}$, any integer $k ,k_-, k_+ \in \mathbb{Z}_{\geq 0}$ and any $q \in \Z_{\geq -1}$: 
\\
{\rm(1)} The subspace $\mathcal{G}^r_{k, q+1}(S)$ is $C^0$-dense in the space $\mathcal{G}^r_{k, > q}(S)$. 
\\
{\rm(2)} The subspace $\mathcal{G}^r_{k_-, k_+, q+1}(S)$ is $C^0$-dense in the space $\mathcal{G}^r_{k_-, k_+, > q}(S)$. 
\\
{\rm(3)} The subspace $[\mathcal{G}^r_{k_-, k_+, q+1}(S)]$ is dense in the quotient space $[\mathcal{G}^r_{k_-, k_+, > q}(S)]$ with respect to the quotient topology of the $C^r$-topology. 
\end{lemma}

\begin{proof}
Fix any flow $v \in \mathcal{G}^r_{k_-, k_+, > q}(S)$. 
It suffices to show that there is the resulting flow from $v$ by a small perturbation belonging to $\mathcal{G}^r_{k_-, k_+, q+1}(S)$ for the density of $\mathcal{G}^r_{k_-, k_+, q+1}(S)$. 
We may assume that $\mathrm{codim} (v) > q+1$. 
Write $q' := \mathrm{codim} (v) - (q+1) > 0$.
If there are $q'$ multi-saddle separatrices of $v$ outside of $\partial S$, then the resulting flow from $v$ by a small $C^r$-perturbation that cuts such multi-saddle separatrices belongs to $\mathcal{G}^r_{k_-, k_+, q+1}(S)$. 
Thus, we may assume that $v$ contains at most $q'-1$ multi-saddle separatrices of $v$ outside of $\partial S$. 
From Lemma~\ref{lem:WH}, assertion {\rm(3)} holds. 

\begin{claim}\label{claim:06}
We may assume that $v$ contains no multi-saddle separatrices of $v$ outside of $\partial S$. 
\end{claim}

\begin{proof}[Proof of Claim~\ref{claim:06}]
By a $C^r$-small perturbation that cuts multi-saddle separatrices of $v$ outside of $\partial S$, we have that the resulting flow is arbitrarily near $v$ and has no such multi-saddle separatrices. 
Replacing $v$ with the resulting flow, the claim holds. 
\end{proof}

Then $q' = \mathrm{codim} (v) - (q+1) = \mathrm{codim}_{\mathrm{m}}(v) - (q+1) > 0$. 
Applying the Whitehead move $q'$ times, we obtain that the resulting flow is $C^0$-near $v$ and belongs to $\mathcal{G}^r_{k_-, k_+, q+1}(S)$. 
This implies the $C^0$-density of $\mathcal{G}^r_{k, q+1}(S)$ in $\mathcal{G}^r_{k, > q}(S)$, which implies assertion (1). 
Since $\mathcal{G}^r_{k, q+1}(S) = \bigsqcup \{ \mathcal{G}^r_{k_-, k_+, q+1}(S) \mid k_- +  k_+ = k \}$, the subspace $\mathcal{G}^r_{k_-, k_+, q+1}(S)$ is $C^0$-dense in the space $\mathcal{G}^r_{k_-, k_+, > q}(S)$, which implies assertion (2).  
\end{proof}

Though the previous lemma shows only $C^0$-density, Lemma~\ref{lem:WH} suggests that it works well in the proof for the $C^r$ case.
Thus, the author would like to know whether the subspace $\mathcal{G}^r_{k_-, k_+, q+1}(S)$ for any $k_-, k_+ \in \mathbb{Z}_{\geq 0}$ is $C^r$-dense in the space $\mathcal{G}^r_{k_-, k_+, > q}(S)$. 

Lemma~\ref{lem:codimension_open} and Lemma~\ref{lem:codimension_density} imply Theorem~\ref{prop:codimension}. 
We will show a similar statement for ``gradient flow with limit cycles'' in Section~\ref{sec:07}.

\section{Abstract cell complex structure and filtration of $\mathcal{G}^r_{k_-/2, k_+/2, >-1}(S)$}

In this section, we construct the abstract cell complex structure of the space of gradient flows. 
The compactness of the surface implies the following finite property. 

\begin{lemma}\label{lem:top_eq_class}
For any $r \in \mathbb{Z}_{\geq 0} \sqcup \{ \infty \}$, any integers $k_-, k_+ \in \mathbb{Z}_{\geq 0}$, and any $q \in \Z_{\geq -1}$, the subspace $\mathcal{G}^r_{k_-/2, k_+/2, q}(S)$  contains at most finitely many topological equivalence classes.  
\end{lemma}

\begin{proof}
Since the number of multi-saddles is bounded, there are at most finitely many multi-saddle connections that appear in $\mathcal{G}^r_{k_-/2, k_+/2, q}(S)$. 
This implies the finite possible combination obtained by the inverse operations of annihilation operations of non-self-connected separatrices and the Whitehead move finitely many times. 
\end{proof}

%

\subsection{Abstract cell complex structure}

\subsubsection{Abstract cell complex}
For a set $S$ with a transitive relation $\prec$ and a function $\dim \colon S \to \Z_{\geq 0}$, the triple $(S, \prec, \dim )$ is an {\bf abstract cell complex} if $x \prec y$ implies $\dim x < \dim y$. 
Then $\dim x$ is called the {\bf dimension} of $x$, and $x$ is called a {\bf cell}. 
A $k$-cell is a cell whose dimension is $k$. 
The {\bf codimension} of $x$ is $\sup_{y \in S} \dim y - \dim x$. 
Note that the triple $(S, \prec, \mathrm{codim} )$ can determine the dimension and so the abstract cell complex structure. 
For a finite preordered set $(X, \leq)$, the triple $(X, <, \mathop{\mathrm{ht}})$ is an abstract cell complex. 

\subsubsection{Abstract cell complex structure of the space of gradient flows}

For any $r \in \mathbb{Z}_{\geq 0} \sqcup \{ \infty \}$, any integers $k_-, k_+ \in \mathbb{Z}_{\geq 0}$ and any integer $q \in \Z_{\geq -1}$, denote by $[\mathcal{G}^r_{k_-/2, k_+/2, >q}(S)]$ the quotient space of $\mathcal{G}^r_{k_-/2, k_+/2, >q}(S)$ by the topologically equivalence classes, and by $[\mathcal{G}^r_{k_-/2, k_+/2, q}(S)]$ the quotient space of $\mathcal{G}^r_{k_-/2, k_+/2, q}(S)$ by the topological equivalence classes equipped with the quotient topology of $C^r$ topology. 
We have the following statement.

\begin{theorem}\label{main:02}
The following statements hold for any $r \in \mathbb{Z}_{> 0} \sqcup \{ \infty \}$, any integers $k_-, k_+ \in \mathbb{Z}_{\geq 0}$ and any $q \in \Z_{\geq -1}$:
\\
{\rm(1)} The subset $[\mathcal{G}^r_{k_-/2, k_+/2, >-1}(S)]$ is a finite $T_0$-space. 
\\
{\rm(2)} The subset $[\mathcal{G}^r_{k_-/2, k_+/2, >-1}(S)]$ is an abstract cell complex with respect to the opposite order of the specialization preorder and the codimension. 
\\
{\rm(3)} $\overline{[\mathcal{G}^r_{k_-/2, k_+/2, >-1}(S)]_{q+1}} = [\mathcal{G}^r_{k_-/2, k_+/2, >-1}(S)]_{\geq q+1} = [\mathcal{G}^r_{k_-/2, k_+/2, > q}(S)]$. 
\end{theorem}

\begin{proof}
Lemma~\ref{lem:top_eq_class} implies that $[\mathcal{G}^r_{k_-/2, k_+/2, >q}(S)]$ is a finite $T_0$-space, and that the coheight of $[\mathcal{G}^r_{k_-/2, k_+/2, q}(S)]$ corresponds to the codimension. 
By Lemma~\ref{lem:codimension_density} (3), we have $\overline{[\mathcal{G}^r_{k_-/2, k_+/2, >-1}(S)]_{q+1}} = [\mathcal{G}^r_{k_-/2, k_+/2, >-1}(S)]_{\geq q+1} = [\mathcal{G}^r_{k_-/2, k_+/2, > q}(S)]$. 
\end{proof}

The previous result implies Theorem~\ref{main:02-}. 


\section{Non-contractibility of the space of gradient flows}

In this section, we demonstrate the non-contractibility of a connected component of the space of gradient flows. 
In other words, we will show that the quotient space $[\mathcal{G}^r_{1, 2, >-1}(\Sigma_{0,2})]$ has a connected component which is the weak homotopy type of a bouquet $\mathbb{S}^2 \vee \mathbb{S}^2$ of two two-dimensional spheres.
Here $\Sigma_{g,p}$ is a compact surface whose genus is $g$ and which has $p$ boundary components.
On the other hand, we have the following observation.

\begin{proposition}\label{prop:cont}
For any $r \in \mathbb{Z}_{> 0} \sqcup \{ \infty \}$ and any integers $k_-, k_+ \in \mathbb{Z}_{\geq 0}$, the subspace $[\mathcal{G}^r_{k_-/2, k_+/2, >-1}(\Sigma_{g,0})]$ is contractible or empty. 
\end{proposition}

%

To show the statements, we recall the theory of homotopy types of finite $T_0$-spaces as follows.

\subsection{Homotopy types of finite $T_0$-spaces}

From now on, we equip a finite $T_0$-space $(X, \tau)$ with the specialization preorder $\leq_{\tau}$. 
Notice that $ x \leq_{\tau} y $ if and only if  $ x \in \overline{\{ y \}}$. 
Write the upset $\mathop{\uparrow}_{\tau} x  := \{ y \in X \mid x \leq_{\tau} y \}$ and the downset $\mathop{\downarrow}_{\tau} x  := \{ y \in X \mid y \leq_{\tau} x \}$. 
Notice that $\mathop{\downarrow}_{\tau} x - \{ x \} = \{ y \in X \mid y <_{\tau} x \}$ and $\mathop{\uparrow}_{\tau} x - \{ x \} = \{ y \in X \mid x <_{\tau} y \}$. 
To state the characterization of contractibility, we recall beat points as follows \cite{may2003finite,stong1966finite}. 

\begin{definition}
A point $x \in X$ is a {\bf down beat point} {\rm(}or a {\bf colinear} point {\rm)} if there is a point $z \in X$ with  $\mathop{\downarrow}_{\tau} x - \{ x \} = \mathop{\downarrow}_{\tau} z$. 
\end{definition}

\begin{definition}
A point $x \in X$ is a {\bf up beat point} {\rm(}or a {\bf linear} point {\rm)} if there is a point $z \in X$ with  $\mathop{\uparrow}_{\tau} x - \{ x \} = \mathop{\uparrow}_{\tau} z$. 
\end{definition}

\begin{definition}
A point is a {\bf beat point} if it is a down or up beat point. 
\end{definition}

Notice that a point $x \in X$ is a beat point if and only if there is a point $z \in X$ such that either $\mathop{\uparrow}_{\tau} x - \{ x \} = \mathop{\uparrow}_{\tau} z$ or $\mathop{\downarrow}_{\tau} x - \{ x \} = \mathop{\downarrow}_{\tau} z$, and that the inclusion $X - \{x\} \to X$ for any beat point $x \in X$ is a strong deformation retract. 

A finite $T_0$-space without beat points is called a {\bf minimal finite space}.
A minimal finite space is a {\bf core} of a finite $T_0$-space $X$ if it is a strong deformation retract of $X$. 
In \cite[Theorem~4]{stong1966finite}, Stong proved that any finite $T_0$-space has the unique core up to homeomorphism and that two finite $T_0$-spaces are homotopy equivalent to each other if and only if their cores are homeomorphic to each other. 

\subsection{Contractibility of $[\mathcal{G}^r_{k_-/2, k_+/2, >-1}(\Sigma_{g,0})]$}

We demonstrate the contractibility of some connected components as follows. 

\begin{proof}[Proof of Proposition~\ref{prop:cont}]
Fix any $r \in \mathbb{Z}_{> 0} \sqcup \{ \infty \}$ and any integers $k_-, k_+ \in \mathbb{Z}_{\geq 0}$ such that $[\mathcal{G}^r_{k_-/2, k_+/2, >-1}(\Sigma_{g,0})]$ is not empty. 
Since $\Sigma_{g,0}$ is a closed surface, every singular point is either a sink, a source, or an $k$-saddle for some $k \in \Z_{>0}$. 
Put $k := \chi(\Sigma_{g,0}) - (k_+ + k_-)/2$. 
If $k = 0$, then there are no multi-saddles and so the component $[\mathcal{G}^r_{k_-/2, k_+/2, >-1}(\Sigma_{g,0})]$ is a singleton. 
Thus, we may assume that $k>0$. 
Then the topological equivalence class of a flow with a $k$-saddle is the unique maximal element of $[\mathcal{G}^r_{k_-/2, k_+/2, >-1}(\Sigma_{g,0})]$. 
Since the second maximal elements are up beat points, by induction, the space $[\mathcal{G}^r_{k_-/2, k_+/2, >-1}(\Sigma_{g,0})]$ is contractible. 
\end{proof}

\subsection{Non-contractibility of $[\mathcal{G}^r_{1, 2, >-1}(\Sigma_{0,2})]$}

Theorem~\ref{re_prop:5.3} is rewritten as follows.  

\begin{proposition}\label{prop:5.3}
The quotient space $[\mathcal{G}^r_{1, 2, >-1}(\Sigma_{0,2})]$ has a connected component that is weakly homotopy equivalent  to a bouquet $\mathbb{S}^2 \vee \mathbb{S}^2$ of two two-dimensional spheres. 
\end{proposition}

To demonstrate the previous proposition, we introduce a notation and show some technical statements as follows. 
%

\subsubsection{Weak {\rm(}beat{\rm)} points}

We also recall weak beat points \cite{barmak2008simple,barmak2008one} as follows. 

\begin{definition}
A point $x \in X$ is a {\bf down weak point} if $\mathop{\downarrow}_{\tau} x - \{ x \}$ is contractible. 
\end{definition}

\begin{definition}
A point $x \in X$ is a {\bf up beat point} if $\mathop{\uparrow}_{\tau} x - \{ x \}$ is contractible. 
\end{definition}

\begin{definition}
A point is a {\bf weak beat point} {\rm(}or a {\bf weak} point {\rm)} if it is a down or up weak point. 
\end{definition}

Notice that a point $x \in X$ is a weak beat point if and only if either $\mathop{\uparrow}_{\tau} x - \{ x \}$ or $\mathop{\downarrow}_{\tau} x - \{ x \}$ is contractible, and that the inclusion $X - \{x\} \to X$ for any weak beat point $x \in X$ is a weak homotopy equivalence \cite[Proposition~3.3]{barmak2008simple}. 

\subsubsection{Connected components}

For any non-negative integer $k_{-,1}$, $k_{-,2}$, $k_{+,1}$, $k_{+,2} \in \mathbb{Z}_{\geq 0}$, denote by 
\[
\bm{\mathcal{G}^r_{k_{-,1}, k_{-,2}/2, k_{+,1}, k_{+,2}/2, > -1}(S)} \subset \mathcal{G}^r_{k_-/2, k_+/2, > -1}(S)
\] 
the subset of flows in $\mathcal{G}^r_{k_-/2, k_+/2, > -1}(S)$ whose sums of indices of sinks (resp. $\partial$-sinks, sources, $\partial$-sources) are $k_{-,1}$ (resp. $k_{-,2}/2, k_{+,1}, k_{+,2}/2$). 
%
By Corollary~\ref{cor:inv_index_02} and Proposition~\ref{lem:no_merge_msc}, the subspace $\mathcal{G}^r_{k_{-,1},k_{-,2}/2, k_{+,1}, k_{+,2}/2, > -1}(S)$ is a connected component of $\mathcal{G}^r_{k_-/2, k_+/2, > -1}(S)$. 

\subsubsection{Non-contractible connected component}

Because any flow in $\mathcal{G}^r_{1,0,2,0, 0}(\Sigma_{0,2})$ has a sink, we can consider that the sink is the point at infinity and so that such flows can be identified with flows on the plane $\R^2$. 
To demonstrate the previous proposition, we have the following statements.  

\begin{figure}
\begin{center}
\includegraphics[scale=0.3]{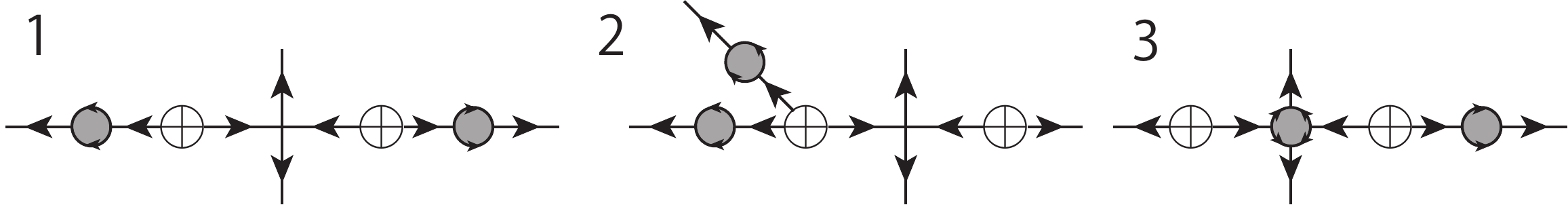}
\end{center}
\caption{Codimension zero elements in $[\mathcal{G}^r_{1,0,2,0, >-1}(\Sigma_{0,2})]$.}
\label{fig:codimension_zero}
\end{figure}

\begin{lemma}\label{lem:000}
The following statements hold for any flow in $\mathcal{G}^r_{1,0,2,0, >-1}(\Sigma_{0,2})$: 
\\
{\rm(1)} Two separatrices to a saddle connect between two different sources. 
\\
{\rm(2)} Two separatrices to a same boundary component connect between two different sources. 
\\
{\rm(3)} The resulting surface of $\Sigma_{0,2}$ removing a sink can be considered as a two-punctured plane. 
\\
{\rm(4)} The sum of indices of multi-saddles is $-3$. 
\\
{\rm(5)} The sum of indices of multi-saddles on a boundary component is either $-1$ or $-2$. 
\end{lemma}

\begin{proof}
Assume that two separatrices to a saddle (resp. a boundary component connected from a source) as on the left (resp. right) of Figure~\ref{fig:codimension_zero+} connect a source with itself. 
\begin{figure}
\begin{center}
\includegraphics[scale=0.3]{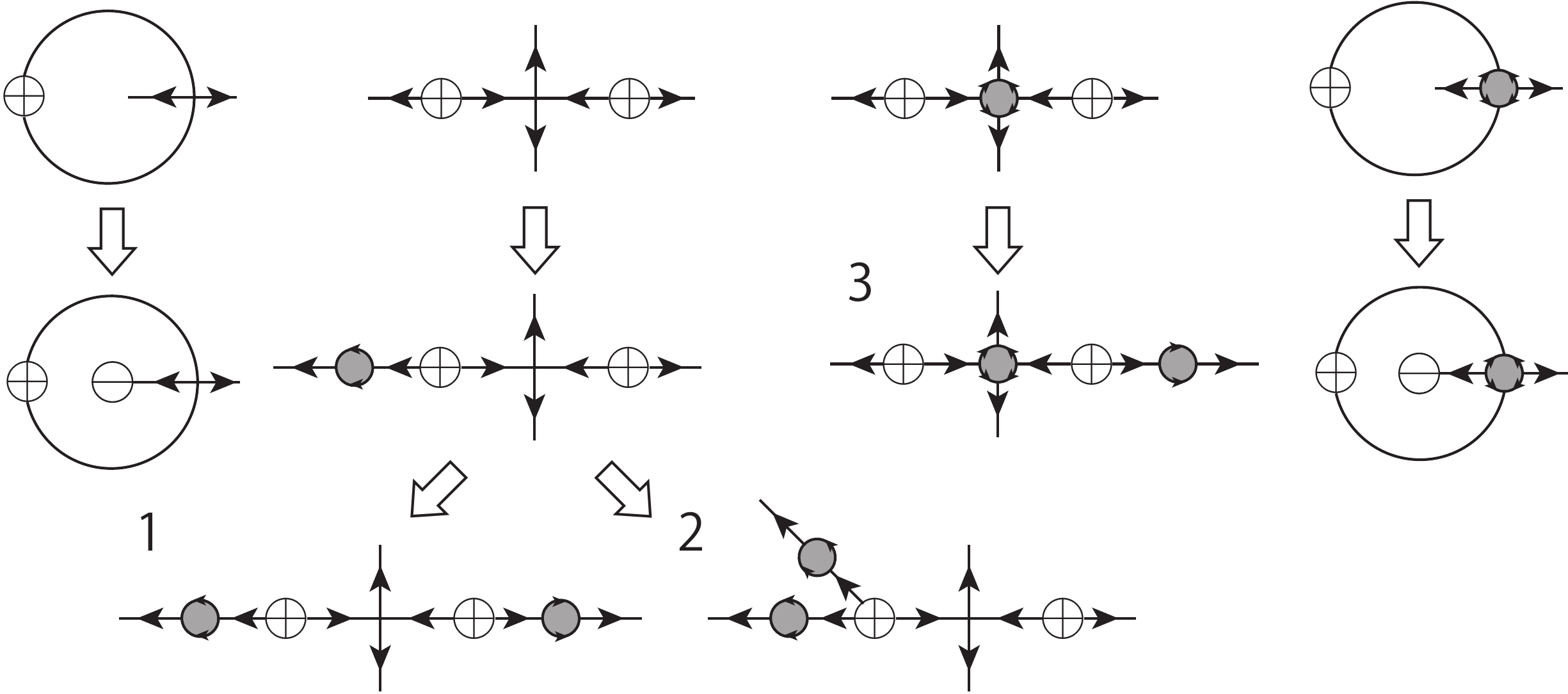}
\end{center}
\caption{The constructions of combinatorial structures of saddle connections of codimension zero elements.}
\label{fig:codimension_zero+}
\end{figure} 
Then a sink exists in the plane, which contradicts that a sink exists only at the point at infinity. 
Therefore, the assertion {\rm(1)} (resp. {\rm(2)}) holds. 

The surface $\Sigma_{0,2}$ is a closed annulus whose Euler characteristic is zero. 
The resulting surface of $\Sigma_{0,2}$ removing a sink can be considered as a two-punctured plane. 
This means that the assertion {\rm(3)} holds. 

The singular point set of any flow in $\mathcal{G}^r_{1,2, 0}(\Sigma_{0,2})$ consists of sinks, $\partial$-sinks, sources, $\partial$-sources, and multi-saddles. 
Since the sum of attracting or repelling singular points is $3$ and the Euler characteristic of any closed annulus is zero, by Poincar{\'e}-Hopf theorem, the sum of indices of multi-saddles is $-3$.
Therefore, the assertion {\rm(4)} holds. 

Because the sum of indices of any boundary component without $\partial$-sinks or $\partial$-sources of a flow in $\mathcal{G}^r_{1,0,2,0, >-1}(\Sigma_{0,2})$ is at least one and is an integer, since the closed annulus $\Sigma_{0,2}$ has exactly two boundary components, assertion (4) implies that the sum of indices of multi-saddles on a boundary component without $\partial$-sinks or $\partial$-sources is either $-1$ or $-2$. 
Therefore, the assertion {\rm(5)} holds. 
\begin{figure}
\begin{center}
\includegraphics[scale=0.6]{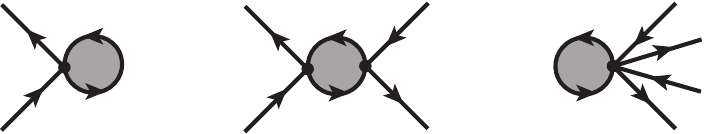}
\end{center}
\caption{Forbidden boundary components with pinchings in gradient flows.}
\label{fig:codimension_one_proof_0a}
\end{figure} 
\end{proof}

\begin{lemma}\label{lem:000+}
Any boundary component for a flow in the space $\mathcal{G}^r_{1,0,2,0, >-1}(\Sigma_{0,2})$ is one of the following structures as in Figure~\ref{fig:codimension_one_proof_0+}: 
\\
{\rm(1)} The singular points on it are exactly two $\partial$-saddles. 
\\
{\rm(2)} The singular points on it are exactly four $\partial$-saddles. 
\\
{\rm(3)} The singular points on it are exactly two $\partial$-saddles and one $1$-$\partial$-saddle. 
\\
{\rm(4)} The singular points on it are exactly one $\partial$-saddle and one $3/2$-$\partial$-saddle. 
\end{lemma}

\begin{proof}
Fix a flow in $[\mathcal{G}^r_{1,0,2,0, >-1}(\Sigma_{0,2})]$. 
Since the surface $\Sigma_{0,2}$ has exactly two boundary components, by Lemma~\ref{lem:000}, the sum of indices of multi-saddles on any boundary component is $-1$ or $-2$. 
Because a boundary component with one singular point and one non-recurrent orbit as on the left of Figure~\ref{fig:codimension_one_proof_0a} does not appear in any gradient flows, if the sum is $-1$ then the singular points on it are two $\partial$-saddles. 
Thus, we may assume that the sum is $-2$. 
Since a boundary component with one or two singular points and one or two non-recurrent orbits as on the middle or right of Figure~\ref{fig:codimension_one_proof_0a} does not appear in any gradient flows, the singular points on the boundary are one of the forms in Figure~\ref{fig:codimension_one_proof_0+}. 
\begin{figure}
\begin{center}
\includegraphics[scale=0.6]{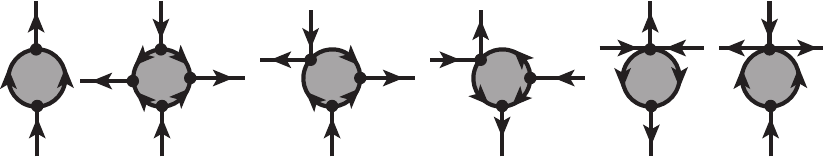}
\end{center}
\caption{Possible local structures on boundary components in gradient flows in $[\mathcal{G}^r_{1,0,2,0, >-1}(\Sigma_{0,2})]$.}
\label{fig:codimension_one_proof_0+}
\end{figure} 
\end{proof}

\begin{lemma}\label{lem:001}
The subspace $[\mathcal{G}^r_{1,0,2,0, 0}(\Sigma_{0,2})]$ consists of three topological equivalence classes, as in Figure~\ref{fig:codimension_zero+}. 
\end{lemma}

\begin{proof}
Lemma~\ref{lem:000}(4) implies that the sum of indices of multi-saddle is $-3$. 
%
From Lemma~\ref{lem:000+}, the set of multi-saddles consists either of six $\partial$-saddles or of one saddle and four $\partial$-saddles. 
By Lemma~\ref{lem:000}(1)--(2), two separatrices from a saddle or a same boundary component connect from two different sources. 
By definition of codimension, there are no multi-saddle separatrices outside of the boundary $\partial \Sigma_{0,2}$. 
Then, we have exactly three possibilities of multi-saddle connections as in the middle of Figure~\ref{fig:codimension_zero+}.
Therefore, three structurally stable gradient flows represent the codimension zero topological equivalence classes of flows as in Figure~\ref{fig:codimension_zero+}.  
\end{proof}

\begin{figure}
\begin{center}
\includegraphics[scale=0.275]{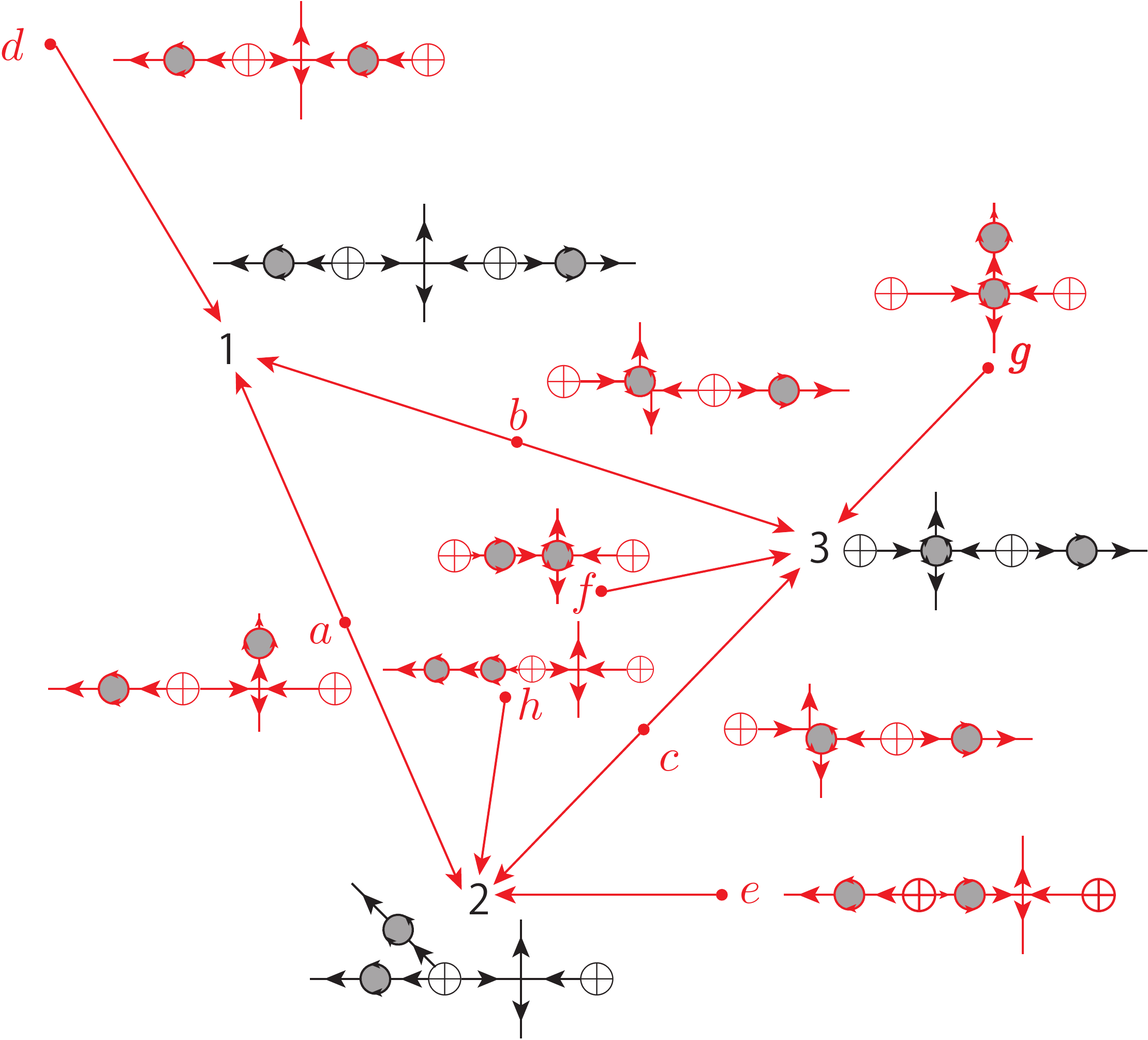}
\end{center}
\caption{Codimension one elements in $[\mathcal{G}^r_{1,0,2,0 >-1}(\Sigma_{0,2})]$.}
\label{fig:codimension_one}
\end{figure}

\begin{lemma}\label{lem:001+}
The subspace $[\mathcal{G}^r_{1,0,2,0, 1}(\Sigma_{0,2})]$ consists of eight topological equivalence classes as in Figure~\ref{fig:codimension_one}. 
\end{lemma}

\begin{proof}
From Lemma~\ref{lem:000}(5), the sum of indices of multi-saddles on a boundary component is at most $-2$. 
By definition of codimension, any codimension one flow has either one pinching or one multi-saddle separatrix outside of the boundary $\partial \Sigma_{0,2}$. 
From Lemma~\ref{lem:000}(4), the sum of indices of multi-saddles is $-3$. 

Suppose that there is exactly one multi-saddle separatrix outside of the boundary. 
From Lemma~\ref{lem:000}(5), the existence of two boundary components implies that there are at least four $\partial$-saddles. 
Therefore, the set of multi-saddles consists either of six $\partial$-saddles or of one saddle and four $\partial$-saddles. 
Then there are exactly six codimension one elements with one multi-saddle separatrix outside of the boundary in $[\mathcal{G}^r_{1,0,2,0 >-1}(\Sigma_{0,2})]$, which are listed in Figure~\ref{fig:codimension_one_proof}. 
\begin{figure}
\begin{center}
\includegraphics[scale=0.24]{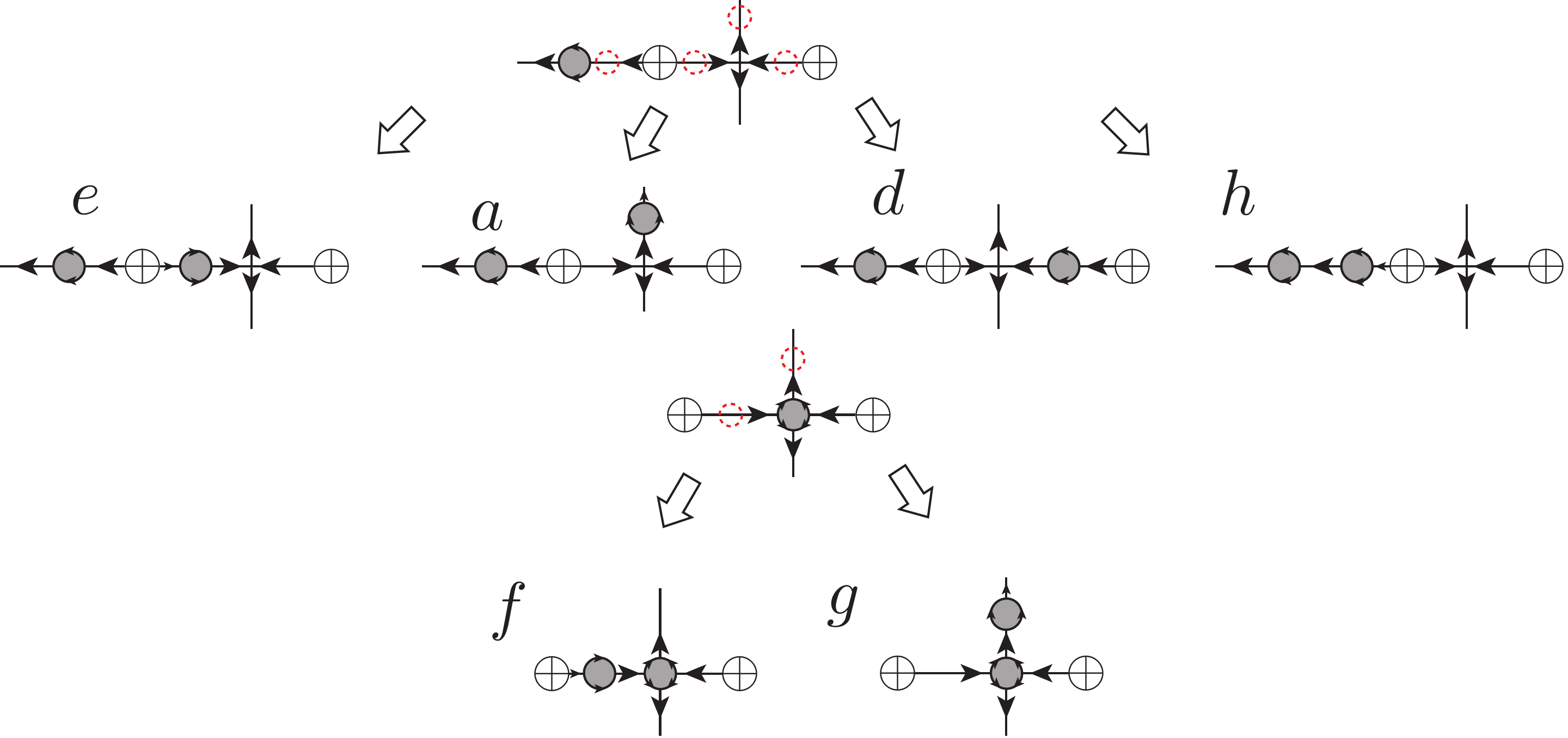}
\end{center}
\caption{The constructions of codimension one elements with one multi-saddle separatrix outside of the boundary in $[\mathcal{G}^r_{1,0,2,0 >-1}(\Sigma_{0,2})]$.}
\label{fig:codimension_one_proof}
\end{figure} 

Suppose that there is exactly one pinching. 
Then there are four $\partial$-saddles and one $1$-$\partial$-saddle. 
Therefore, there are two codimension one elements with one pinching in $[\mathcal{G}^r_{1,0,2,0 >-1}(\Sigma_{0,2})]$ as in Figure~\ref{fig:codimension_one_proof_02} by symmetry. 
\begin{figure}
\begin{center}
\includegraphics[scale=0.25]{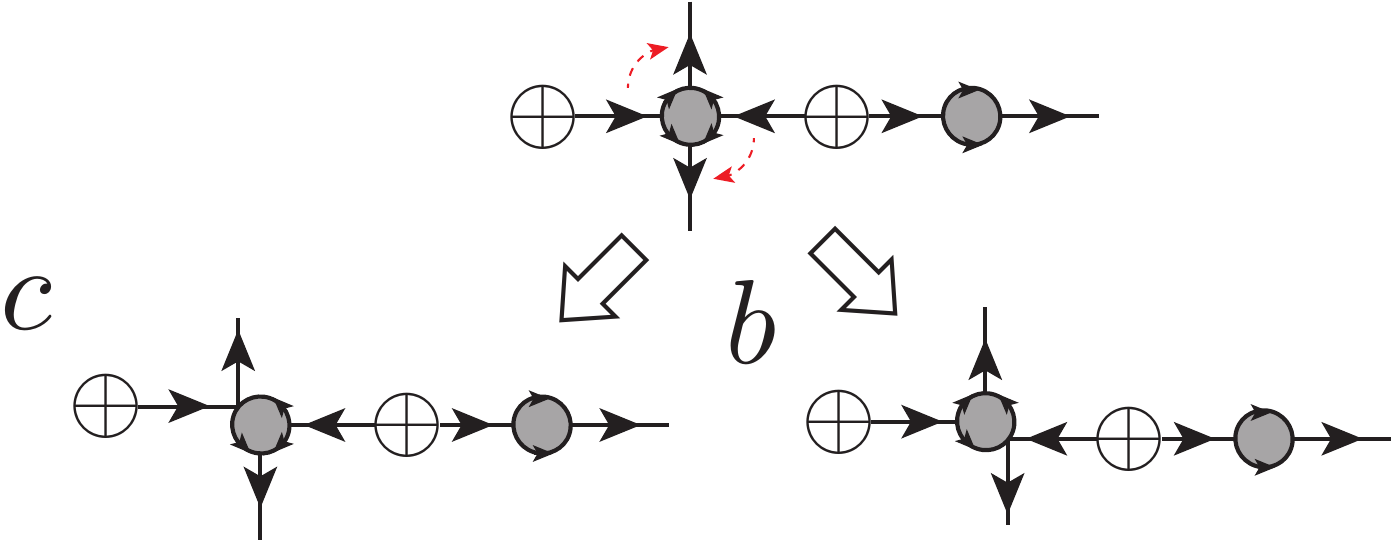}
\end{center}
\caption{The constructions of codimension one elements with one pinching in $[\mathcal{G}^r_{1,0,2,0 >-1}(\Sigma_{0,2})]$.}
\label{fig:codimension_one_proof_02}
\end{figure} 

Local perturbations for pinchings and multi-saddle separatices from or to $\partial$-saddles in the subspace $[\mathcal{G}^r_{1,0,2,0 >-1}(\Sigma_{0,2})]$ as on the upper and middle of Figure~\ref{fig:local_perturb} imply codimension zero elements. 
\begin{figure}
\begin{center}
\includegraphics[scale=0.4]{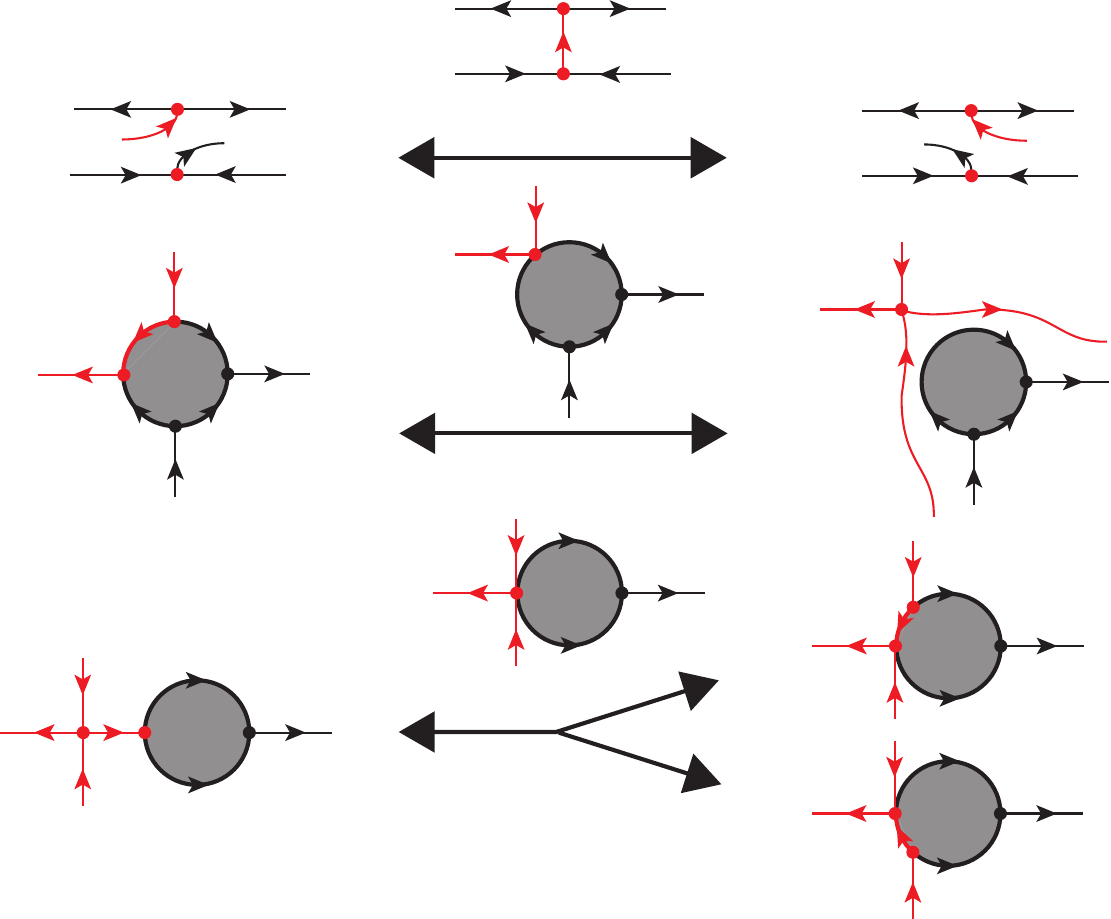}
\end{center}
\caption{Local perturbations for multi-saddle separatices from or to $\partial$-saddles and local perturbations for pinchings in the subspace $[\mathcal{G}^r_{1,0,2,0 >-1}(\Sigma_{0,2})]$.}
\label{fig:local_perturb}
\end{figure} 
Then we have exactly eight possibilities of multi-saddle connections as in Figure~\ref{fig:codimension_one}.
\end{proof}

We list all codimension two elements in $[\mathcal{G}^r_{1,0,2,0 >-1}(\Sigma_{0,2})]$.
%
\begin{figure}
\begin{center}
\includegraphics[scale=0.225]{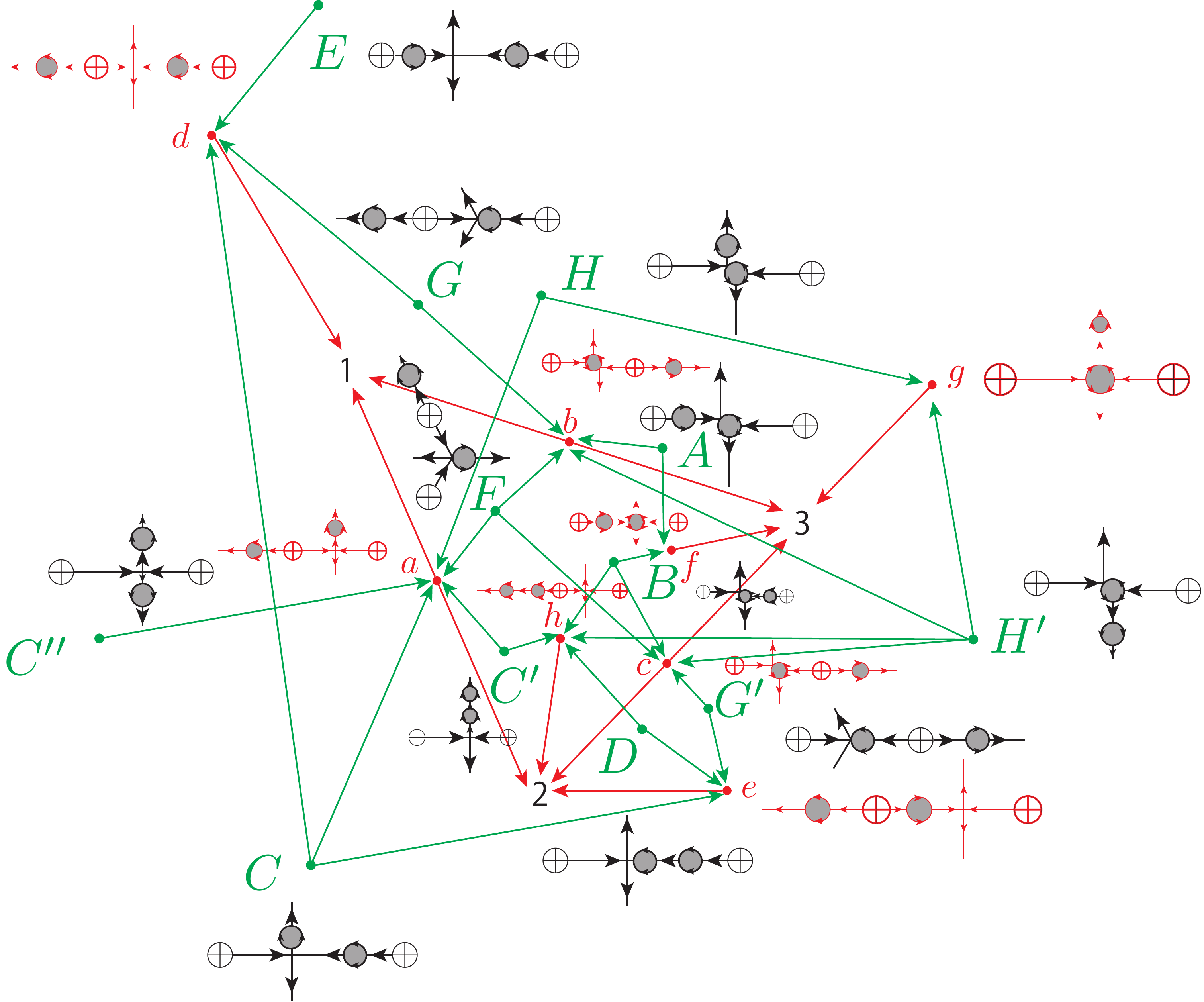}
\end{center}
\caption{Codimension two elements in $[\mathcal{G}^r_{1,0,2,0 >-1}(\Sigma_{0,2})]$.}
\label{fig:codimension_two}
\end{figure} 

\begin{lemma}\label{lem:002}
The subspace $[\mathcal{G}^r_{1,0,2,0, 2}(\Sigma_{0,2})]$ consists of twelve topological equivalence classes as in Figure~\ref{fig:codimension_two}. 
\end{lemma}

\begin{figure}
\begin{center}
\includegraphics[scale=0.275]{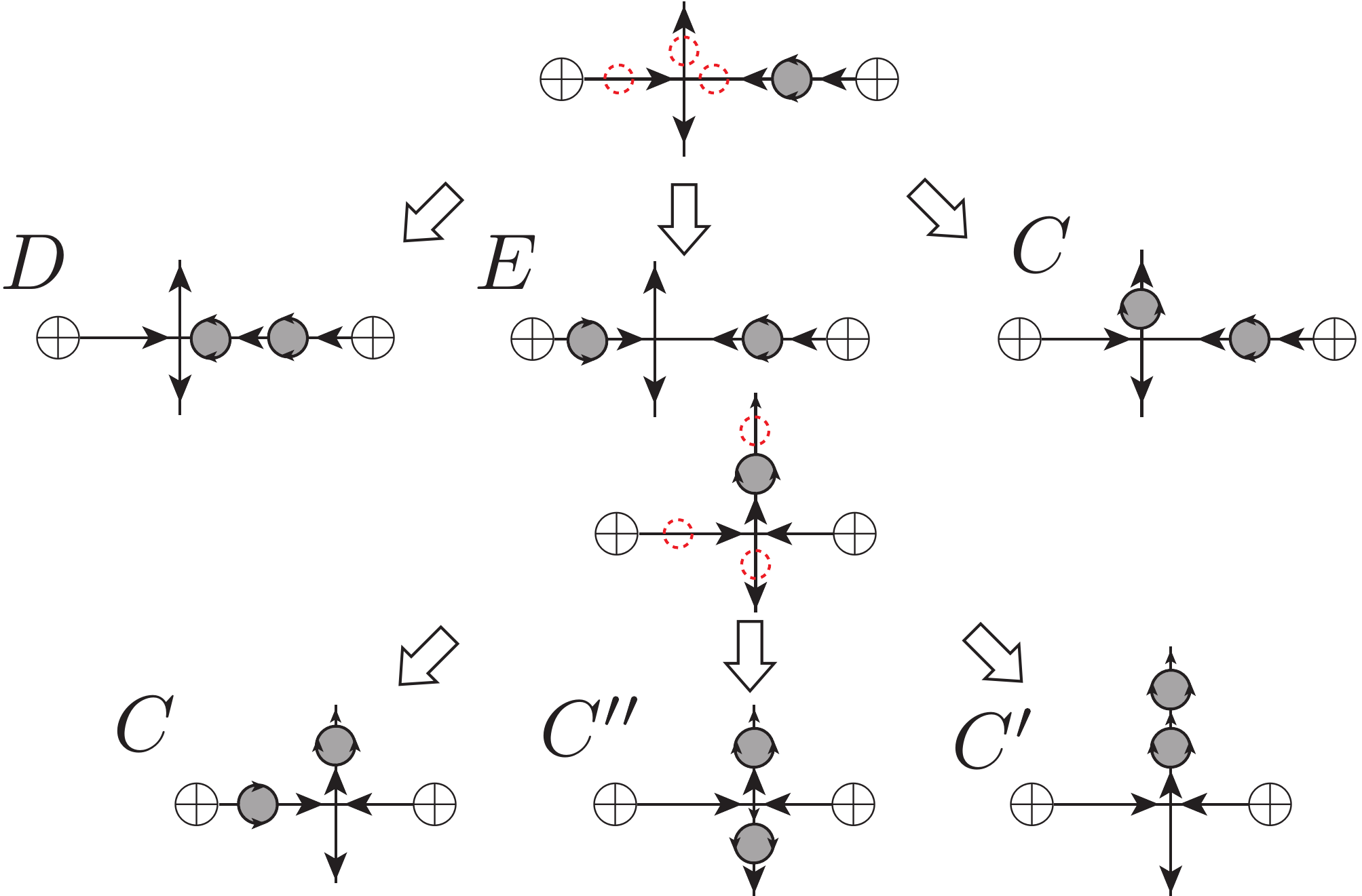}
\end{center}
\caption{The constructions of codimension two elements with two multi-saddle separatrices outside of the boundary in the subspace $[\mathcal{G}^r_{1,0,2,0 >-1}(\Sigma_{0,2})]$.}
\label{fig:codimension_two_proof_02}
\end{figure} 

\begin{figure}
\begin{center}
\includegraphics[scale=0.275]{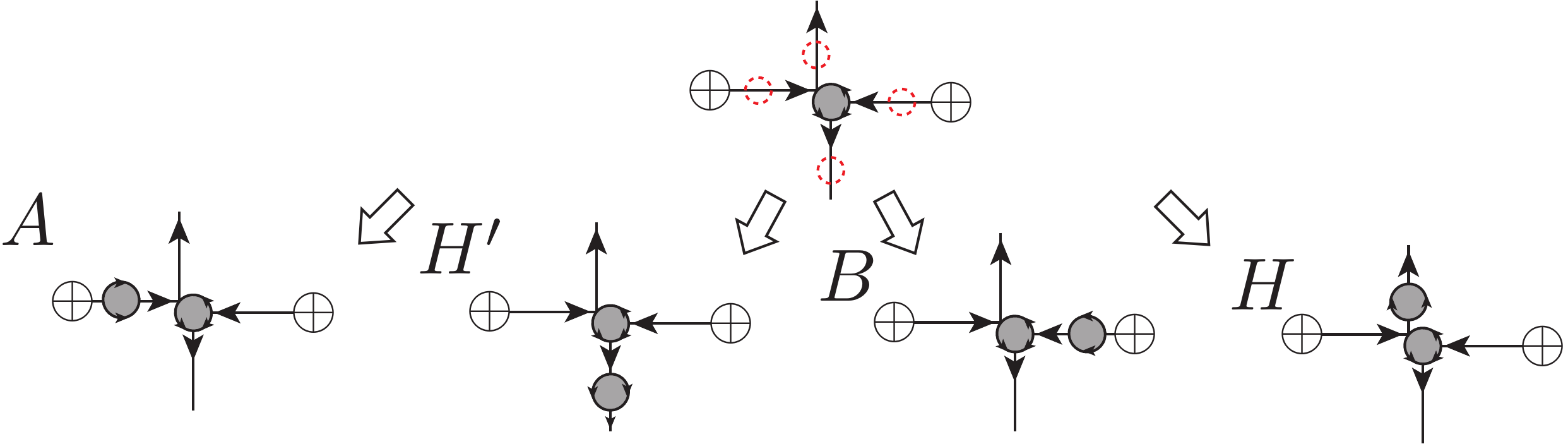}
\end{center}
\caption{The constructions of codimension two elements with one pinching and one multi-saddle separatrix outside of the boundary in $[\mathcal{G}^r_{1,0,2,0 >-1}(\Sigma_{0,2})]$.}
\label{fig:codimension_two_proof_01}
\end{figure} 

\begin{figure}
\begin{center}
\includegraphics[scale=0.25]{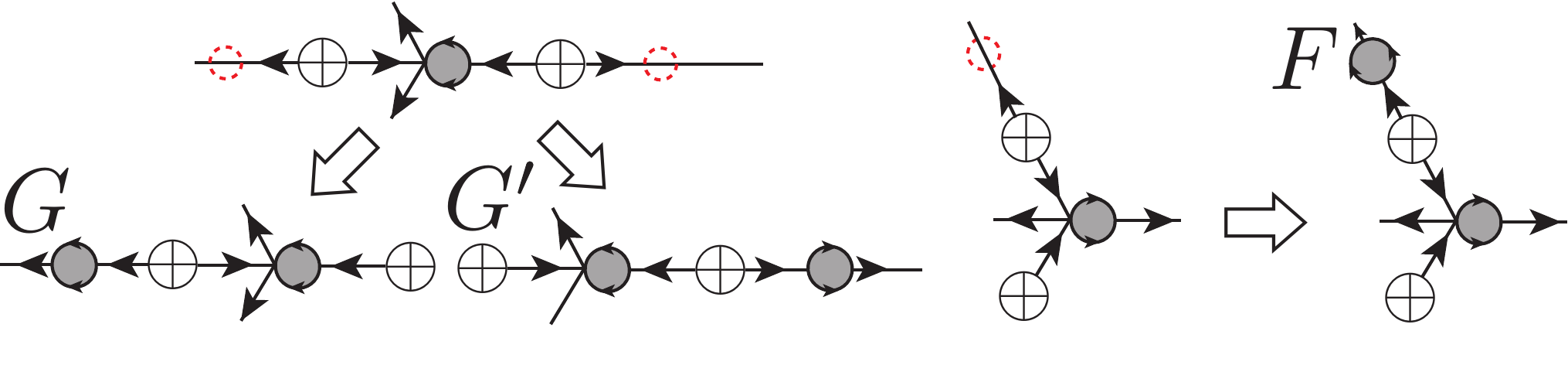}
\end{center}
\caption{The constructions of codimension two elements with one $3/2$-$\partial$-saddle in $[\mathcal{G}^r_{1,0,2,0 >-1}(\Sigma_{0,2})]$.}
\label{fig:codimension_two_proof_00}
\end{figure} 

\begin{proof}
%
By definition of codimension, from Lemma~\ref{lem:000+}, any codimension two flow has either exactly one $3/2$-$\partial$-saddle, a pair of one pinching and one multi-saddle separatrix outside of the boundary $\partial \Sigma_{0,2}$, or two multi-saddle separatrices outside of the boundary $\partial \Sigma_{0,2}$. 
From Lemma~\ref{lem:000}(4), the sum of indices of multi-saddles is $-3$. 

Suppose that there are exactly two multi-saddle separatrices outside of the boundary. 
According to Lemma~\ref{lem:000+}, the existence of two boundary components implies that there are at least four $\partial$-saddles. 
Then there are exactly five codimension two elements with two multi-saddle separatrices outside of the boundary in the subspace $[\mathcal{G}^r_{1,0,2,0 >-1}(\Sigma_{0,2})]$ which are listed in Figure~\ref{fig:codimension_two_proof_02}. 

Suppose that there is exactly one pinching and one multi-saddle separatrix outside of the boundary. 
By Lemma~\ref{lem:000+}, the existence of two boundary components implies that there are four $\partial$-saddles and one $1$-$\partial$-saddle. 
Then there are exactly four codimension two elements with one pinching and one multi-saddle separatrix outside of the boundary in $[\mathcal{G}^r_{1,0,2,0 >-1}(\Sigma_{0,2})]$ which are listed in Figure~\ref{fig:codimension_two_proof_01}. 

Suppose that there is exactly one $3/2$-$\partial$-saddle. 
Then there are exactly three codimension two elements with one $3/2$-$\partial$-saddle in $[\mathcal{G}^r_{1,0,2,0 >-1}(\Sigma_{0,2})]$ which are listed in Figure~\ref{fig:codimension_two_proof_00}. 

Local perturbations for pinchings and multi-saddle separatices from or to $\partial$-saddles in the subspace $[\mathcal{G}^r_{1,0,2,0 >-1}(\Sigma_{0,2})]$ as in Figure~\ref{fig:local_perturb} imply codimension two elements. 
Then we have exactly twelve possibilities of multi-saddle connections as in Figure~\ref{fig:codimension_two}.
\end{proof}

We list all codimension three elements in $[\mathcal{G}^r_{1,0,2,0 >-1}(\Sigma_{0,2})]$.

\begin{figure}
\begin{center}
\includegraphics[scale=0.31]{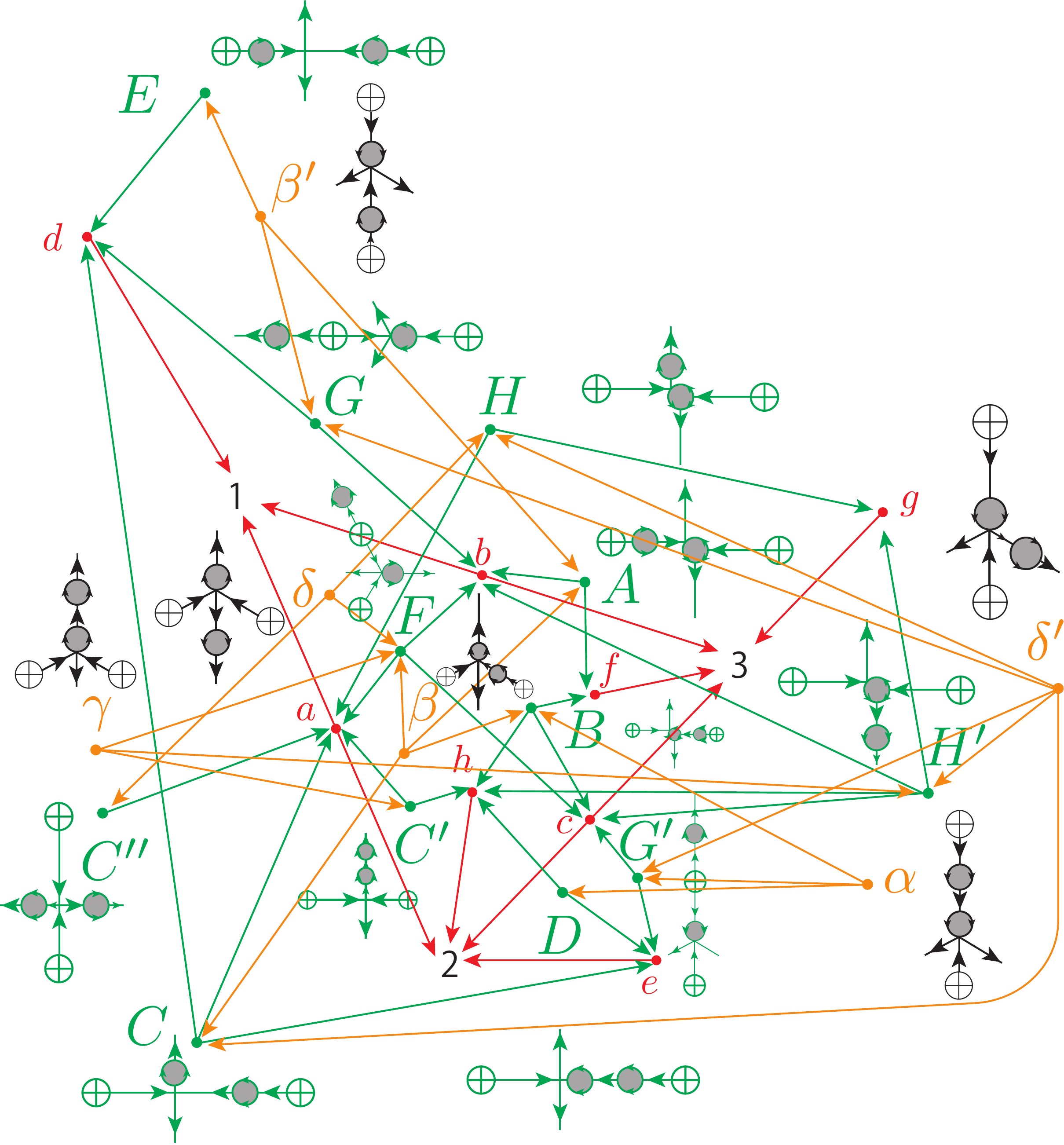}
\end{center}
\caption{Codimension three elements in $[\mathcal{G}^r_{1,0,2,0 >-1}(\Sigma_{0,2})]$.}
\label{fig:codimension_three}
\end{figure}

\begin{lemma}\label{lem:003}
The subspace $[\mathcal{G}^r_{1,0,2,0, 3}(\Sigma_{0,2})]$ consists of six topological equivalence classes as in Figure~\ref{fig:codimension_three}, each of which is represented by a flow with one $3/2$-$\partial$-saddle and one pinching and one multi-saddle separatrix outside of the boundary. 
\end{lemma}

\begin{proof}
By definition of codimension, from Lemma~\ref{lem:000+}, any gradient flows whose topological equivalence classes are codimension three has exactly one $3/2$-$\partial$-saddle and one pinching and one multi-saddle separatrix outside of the boundary as in Figure~\ref{fig:codimension_three_proof}. 
Local perturbations for pinchings and multi-saddle separatices from or to $\partial$-saddles in the subspace $[\mathcal{G}^r_{1,0,2,0 >-1}(\Sigma_{0,2})]$ as in Figure~\ref{fig:local_perturb} imply codimension two elements. 
By listing all the possible combinations, we obtain six multi-saddle connection diagrams as in Figure~\ref{fig:codimension_three}. 
\begin{figure}
\begin{center}
\includegraphics[scale=0.325]{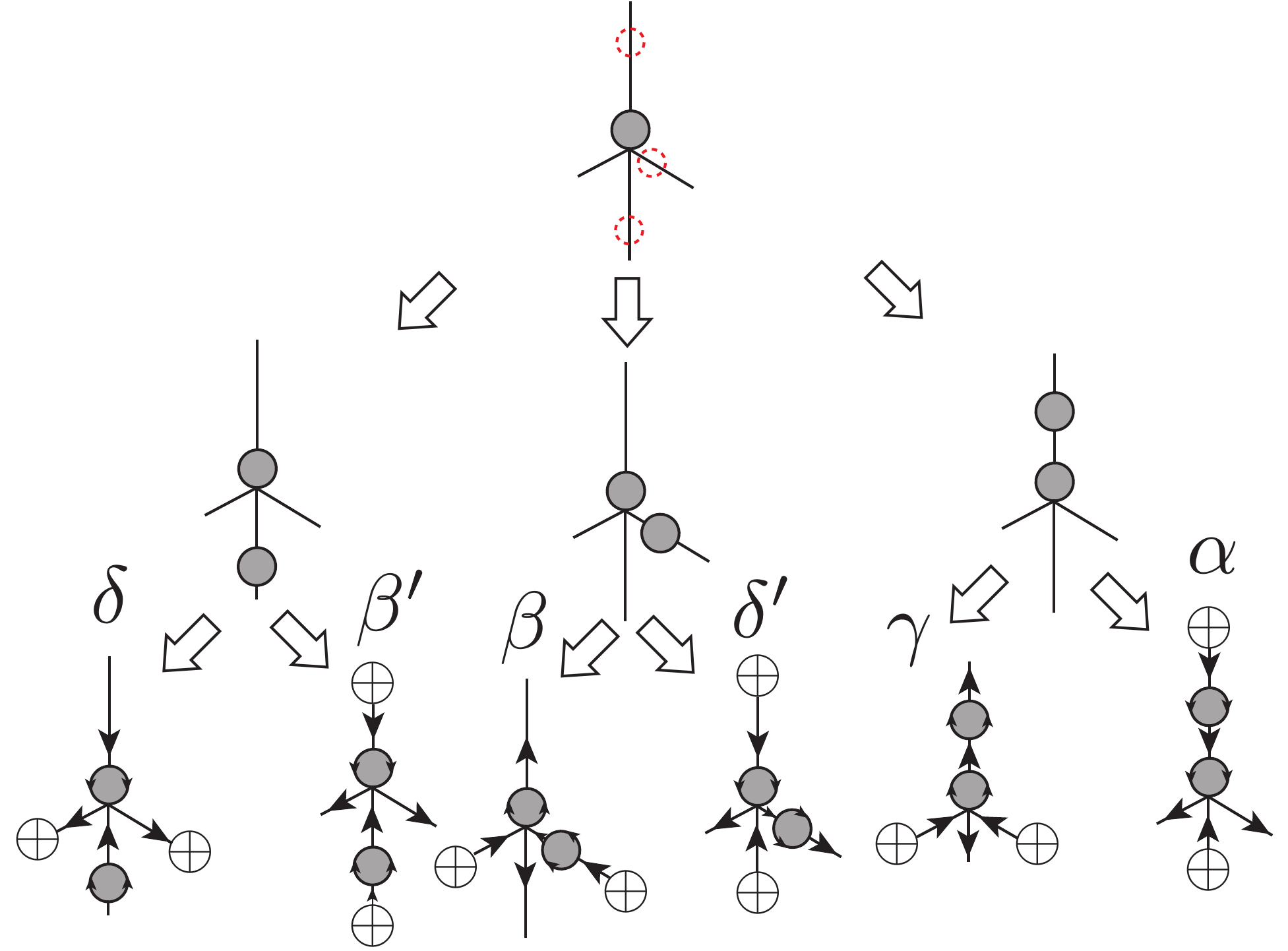}
\end{center}
\caption{The constructions of codimension three elements in $[\mathcal{G}^r_{1,0,2,0 >-1}(\Sigma_{0,2})]$.}
\label{fig:codimension_three_proof}
\end{figure} 
\end{proof}

By the previous three lemmas, we have the following statement. 

\begin{figure}
\begin{center}
\includegraphics[scale=0.275]{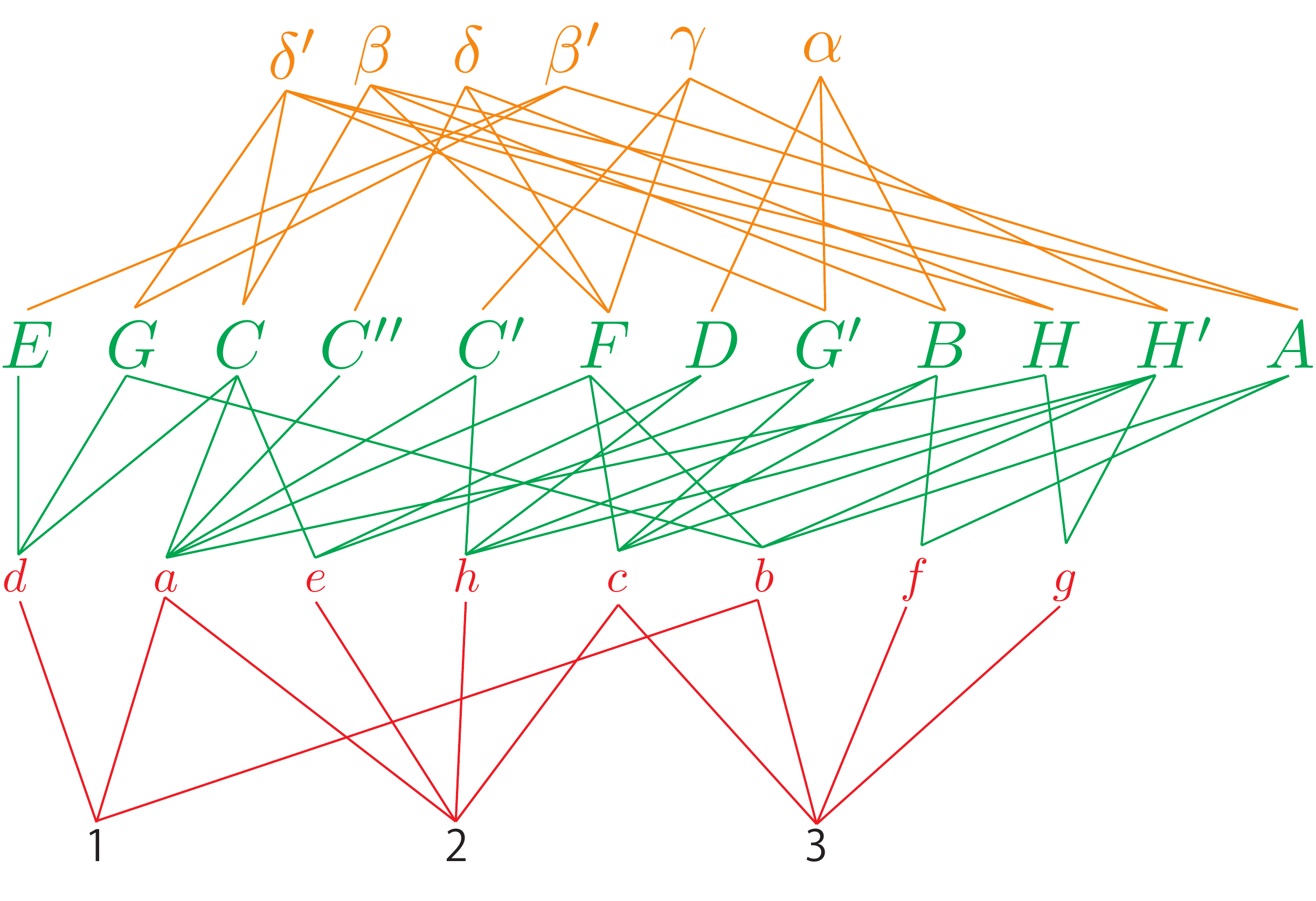}
\end{center}
\caption{The Hesse diagram $X_0$ of the opposite order of the specialization preorder of the finite $T_0$-space $[\mathcal{G}^r_{1,0,2,0 >-1}(\Sigma_{0,2})]$.}
\label{fig:Hessian}
\end{figure}

\begin{lemma}\label{lem:004}
The Hesse diagram of the opposite order of the specialization preorder of the finite $T_0$-space $[\mathcal{G}^r_{1,0,2,0 >-1}(\Sigma_{0,2})]$ is shown in Figure~\ref{fig:Hessian}. 
\end{lemma}

\begin{figure}
\begin{center}
\includegraphics[scale=0.15]{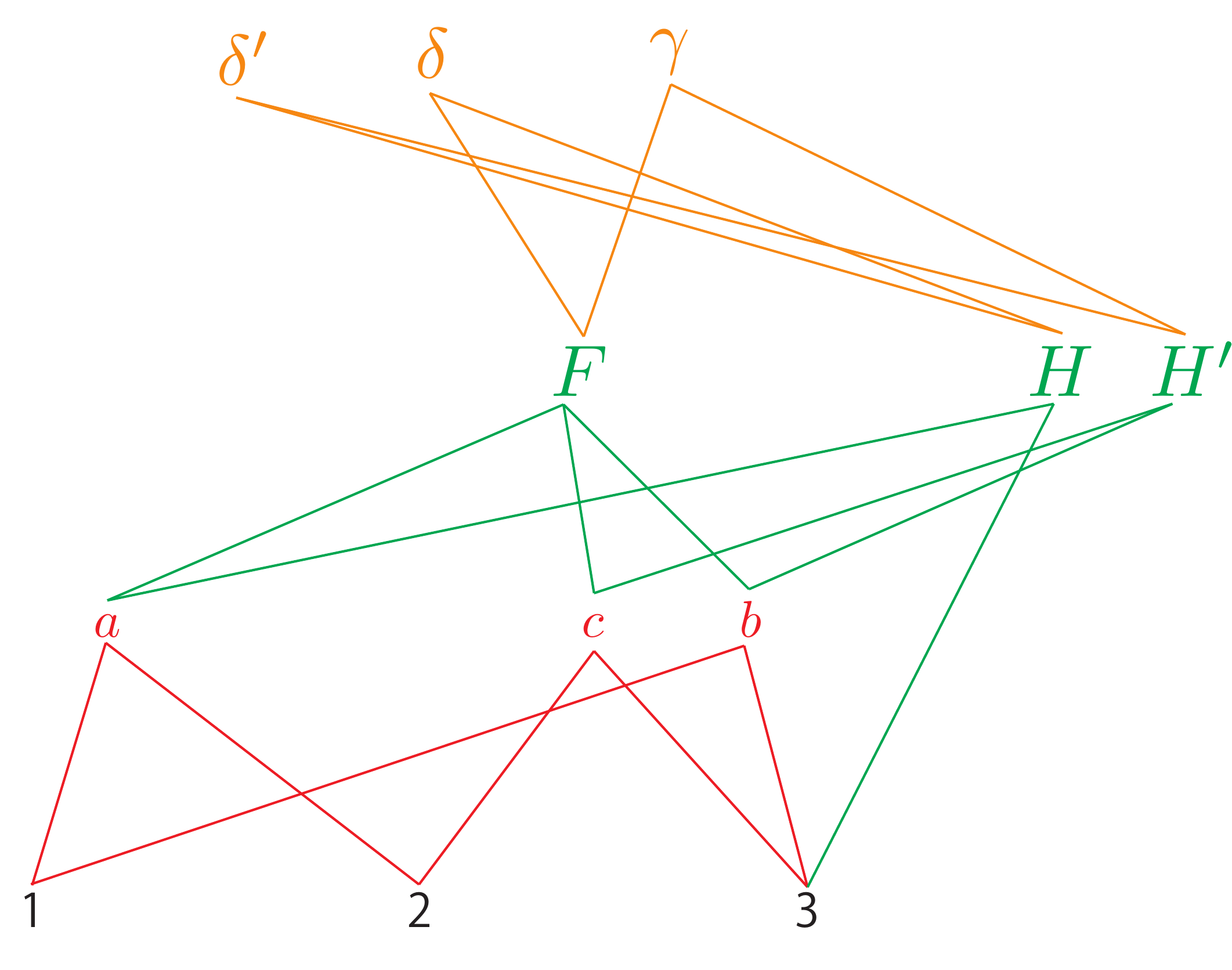}
\includegraphics[scale=0.15]{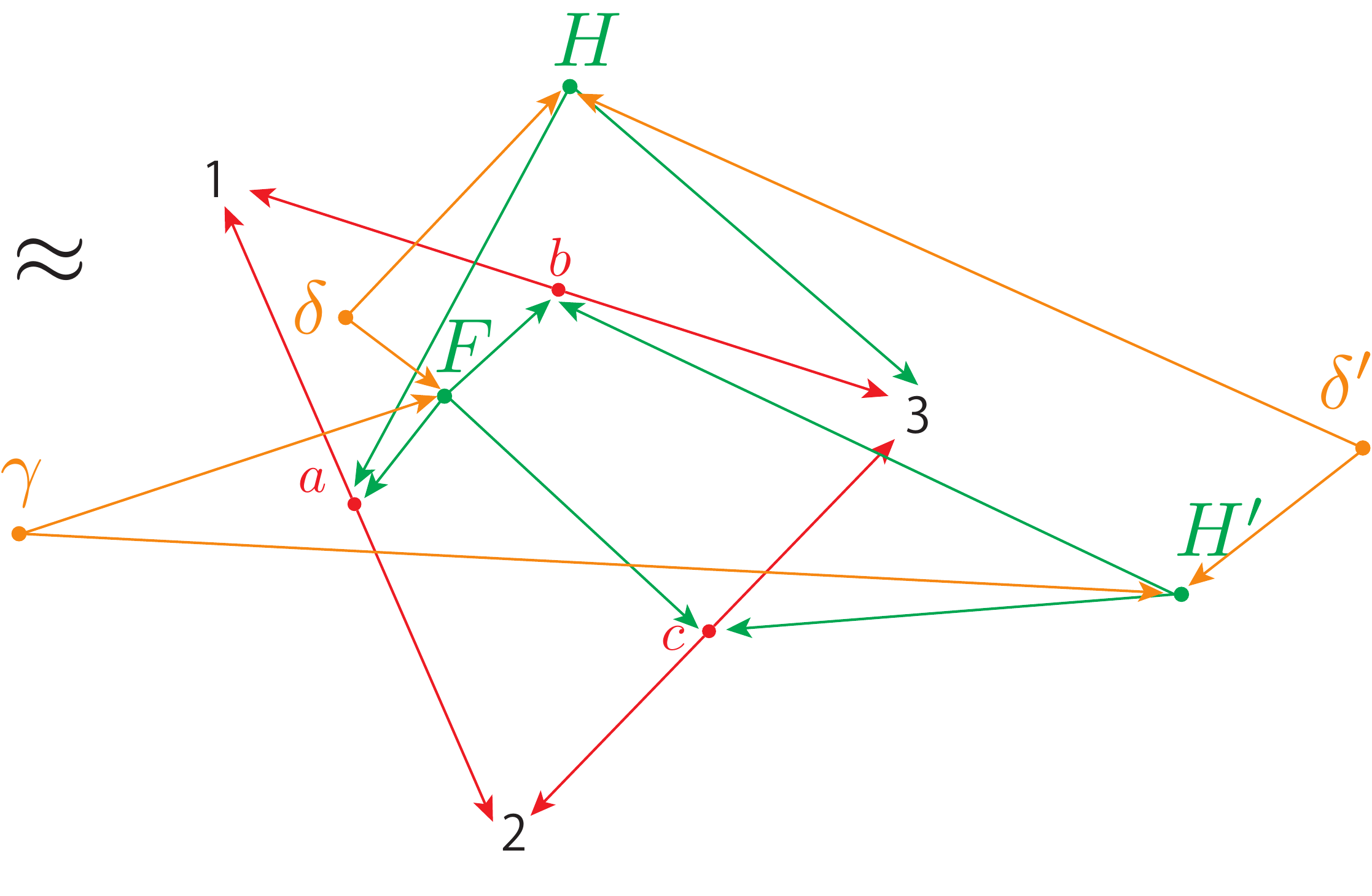}
\end{center}
\caption{The resulting Hesse diagram $X_5$ of the opposite order of the specialization preorder of the resulting diagram of $X_4$ by removing $\beta$.}
\label{fig:Hessian_grad_reduced_04}
\end{figure}

\begin{figure}
\begin{center}
\includegraphics[scale=0.15]{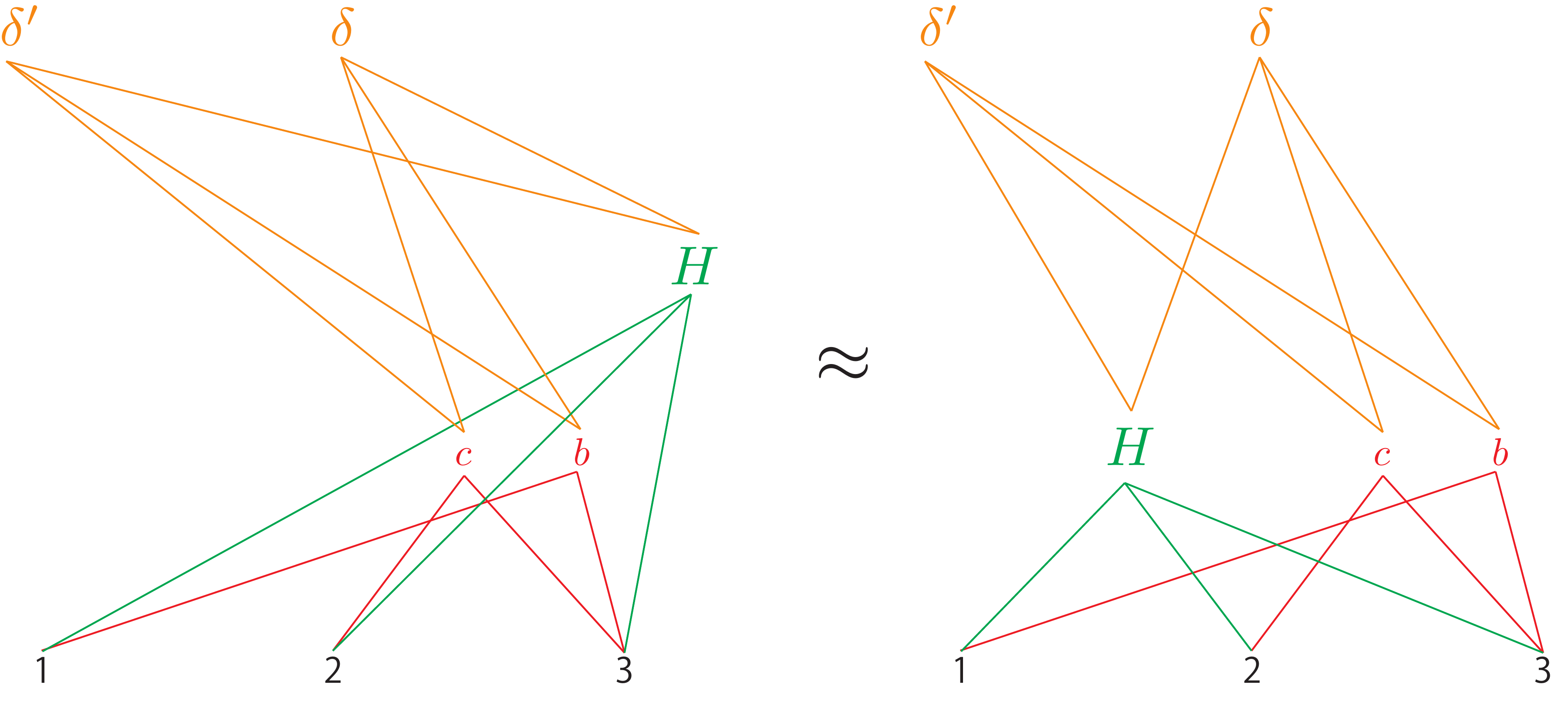}
\includegraphics[scale=0.15]{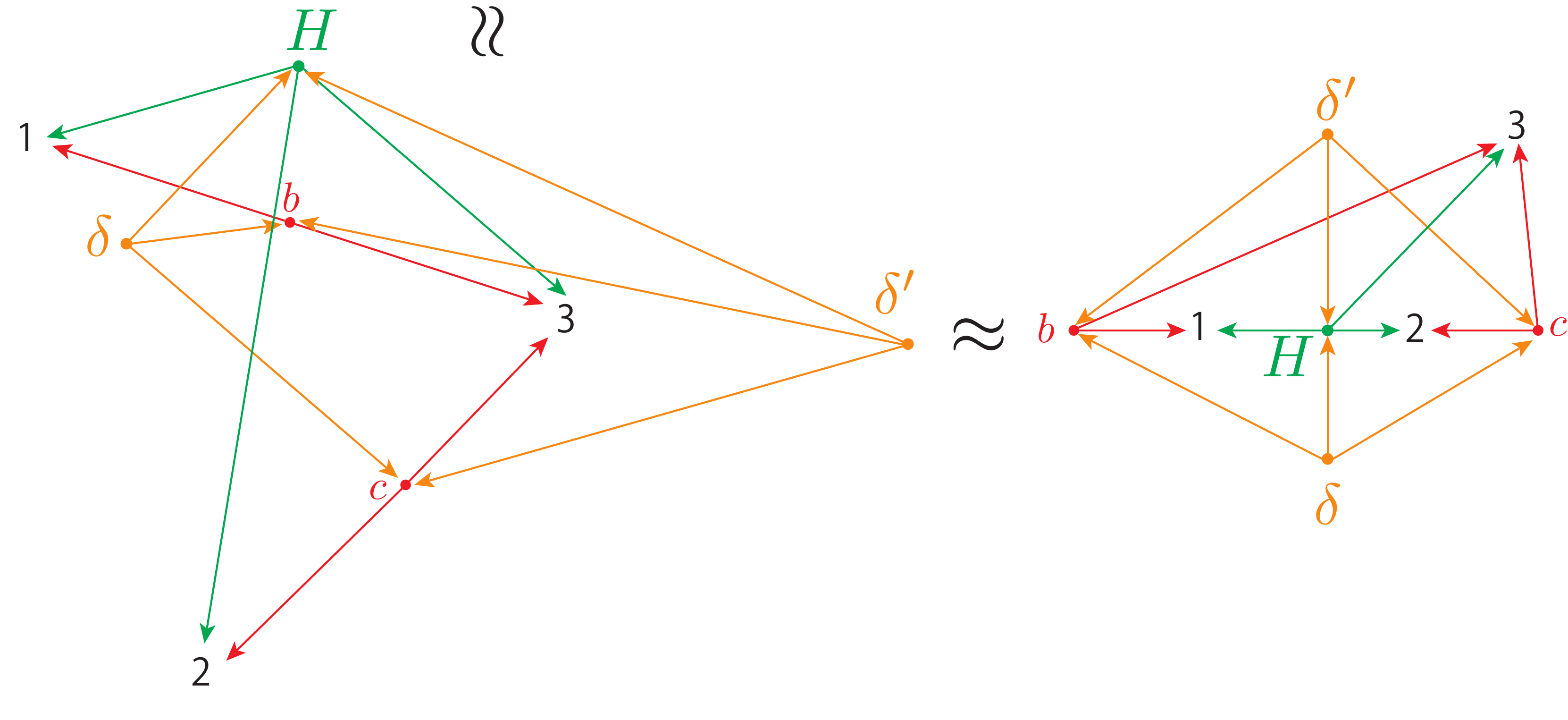}
\end{center}
\caption{The resulting Hesse diagram $X_8$ of the opposite order of the specialization preorder of the resulting diagram of $X_7$ by removing $a$.}
\label{fig:Hessian_grad_reduced_07}
\end{figure}

\begin{lemma}\label{lem:006}
The minimal finite space of the $T_0$-space $[\mathcal{G}^r_{1,0,2,0 >-1}(\Sigma_{0,2})]$ is shown in Figure~\ref{fig:Hessian_grad_reduced_04} and the space $[\mathcal{G}^r_{1,0,2,0 >-1}(\Sigma_{0,2})]$ is weak homotopy equivalent to the finite space whose Hessian is shown in Figure~\ref{fig:Hessian_grad_reduced_07}. 
\end{lemma}

\begin{figure}
\begin{center}
\includegraphics[scale=0.15]{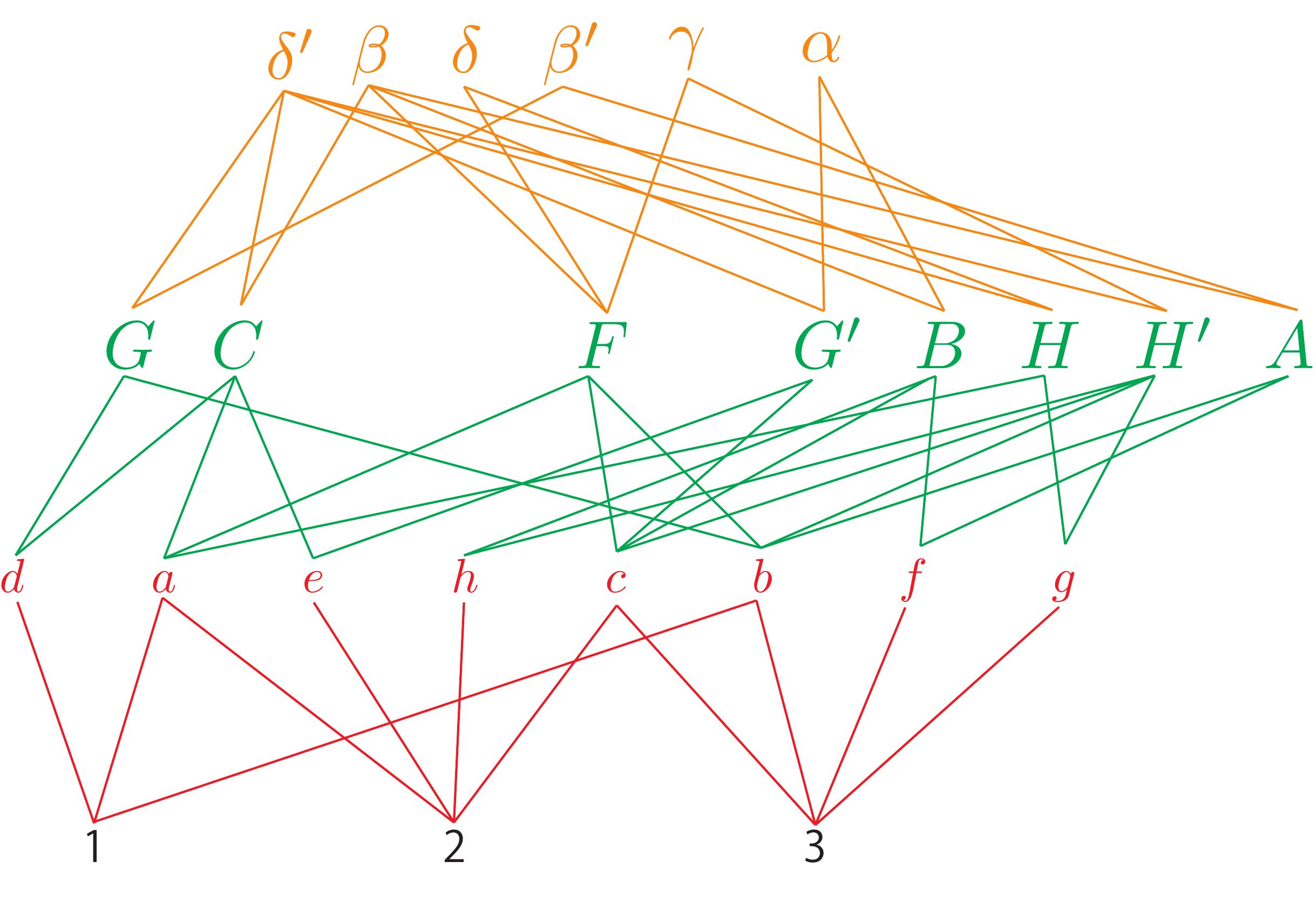}
\includegraphics[scale=0.15]{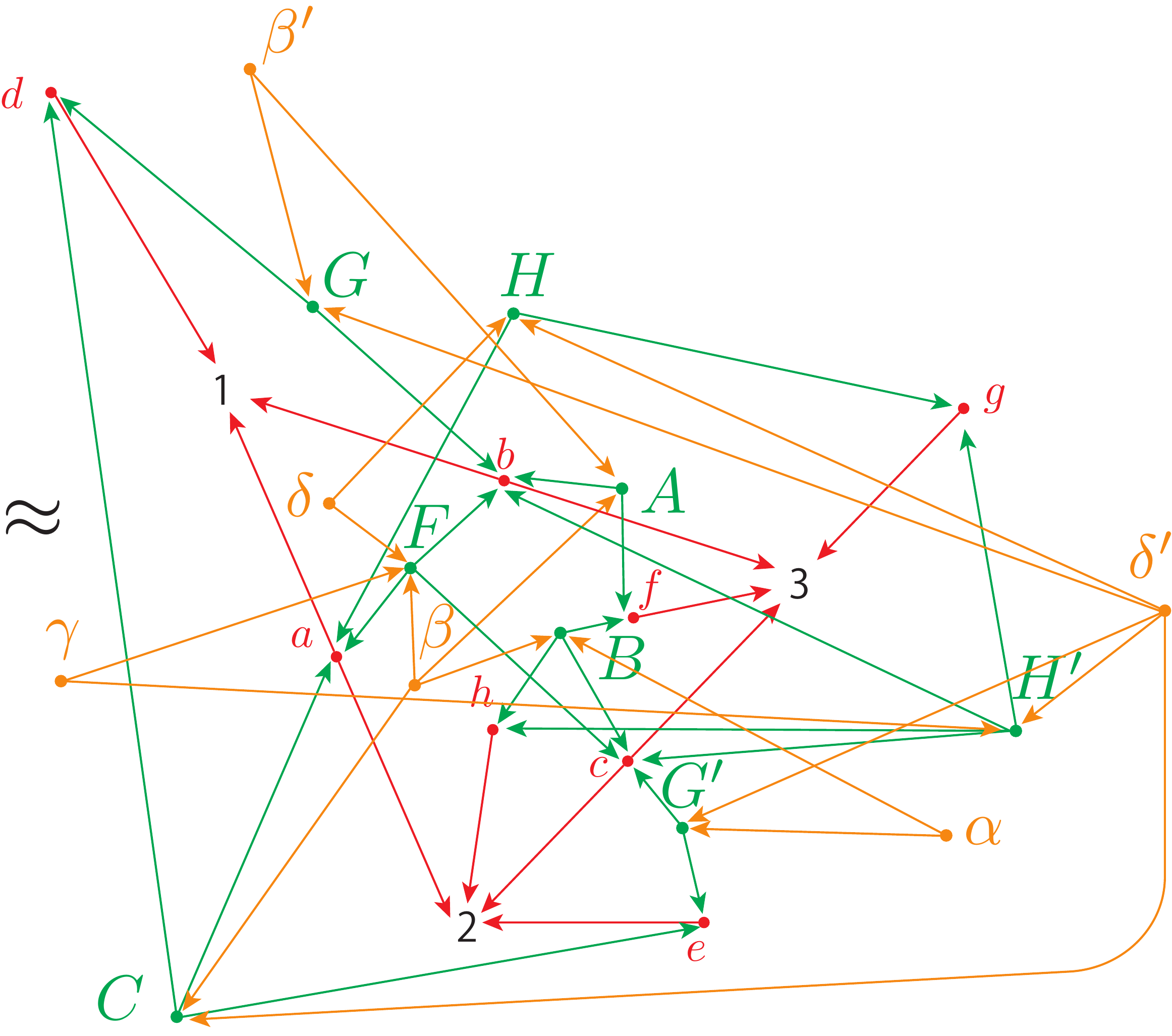}
\end{center}
\caption{The resulting Hesse diagram $X_1$ of the opposite order of the specialization preorder of the resulting diagram of $X_0$ by removing $E$, $C''$, $C'$, and $D$.}
\label{fig:Hessian_grad_reduced_00}
\end{figure}

\begin{figure}
\begin{center}
\includegraphics[scale=0.15]{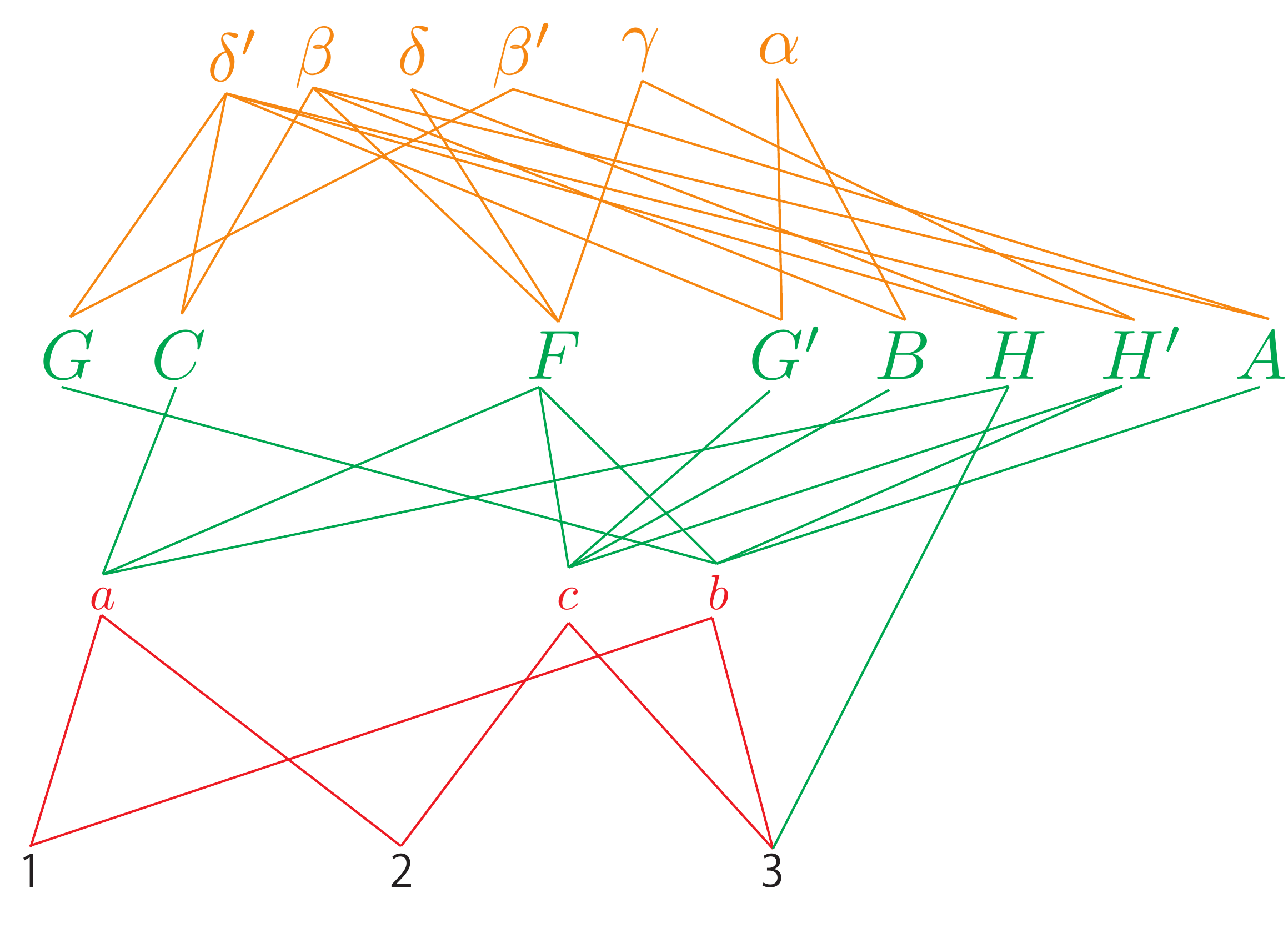}
\includegraphics[scale=0.15]{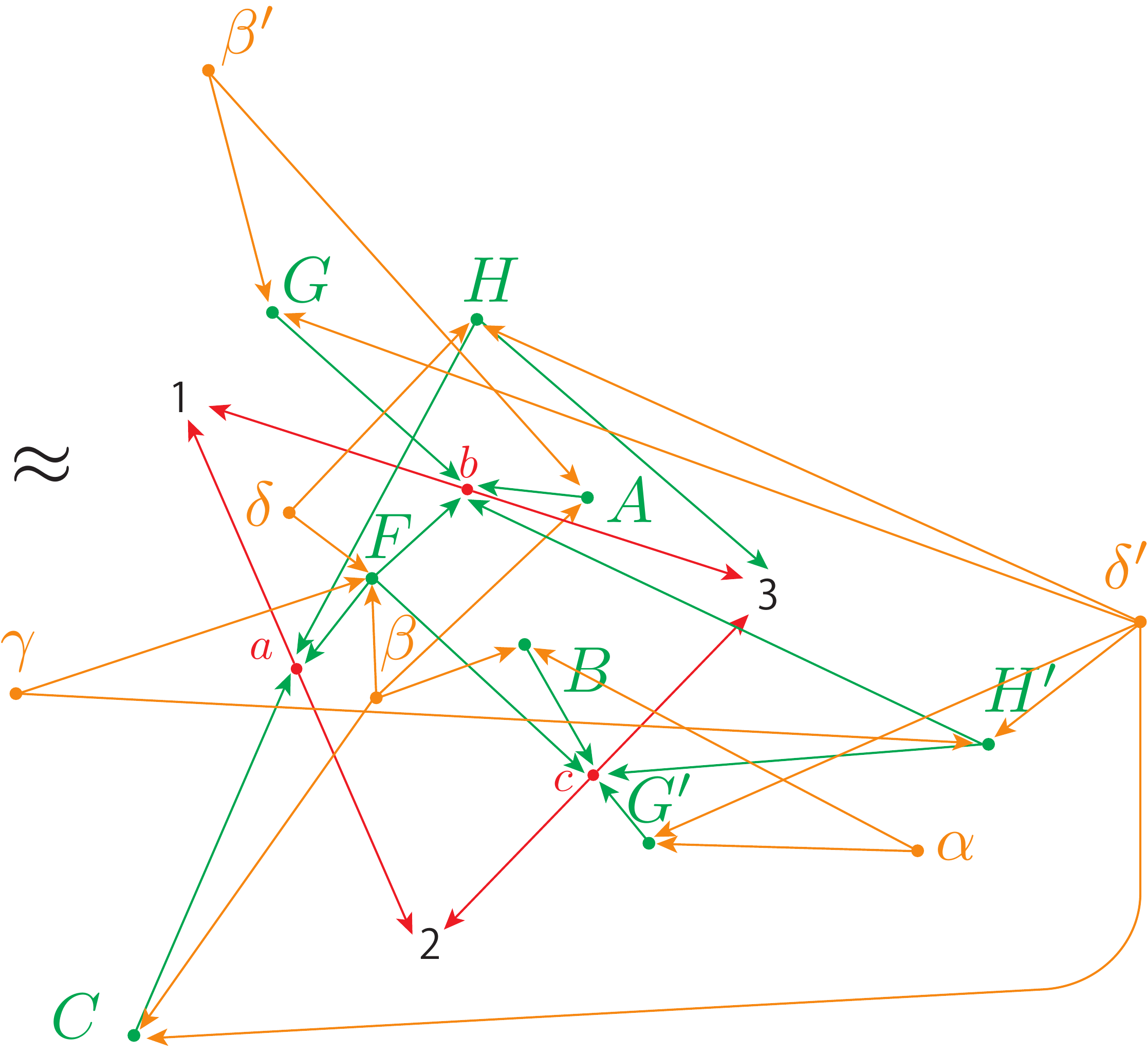}
\end{center}
\caption{The resulting Hesse diagram $X_2$ of the opposite order of the specialization preorder of the resulting diagram of $X_1$ by removing $d, e, h, f$, and $g$.}
\label{fig:Hessian_grad_reduced_01}
\end{figure}

\begin{figure}
\begin{center}
\includegraphics[scale=0.15]{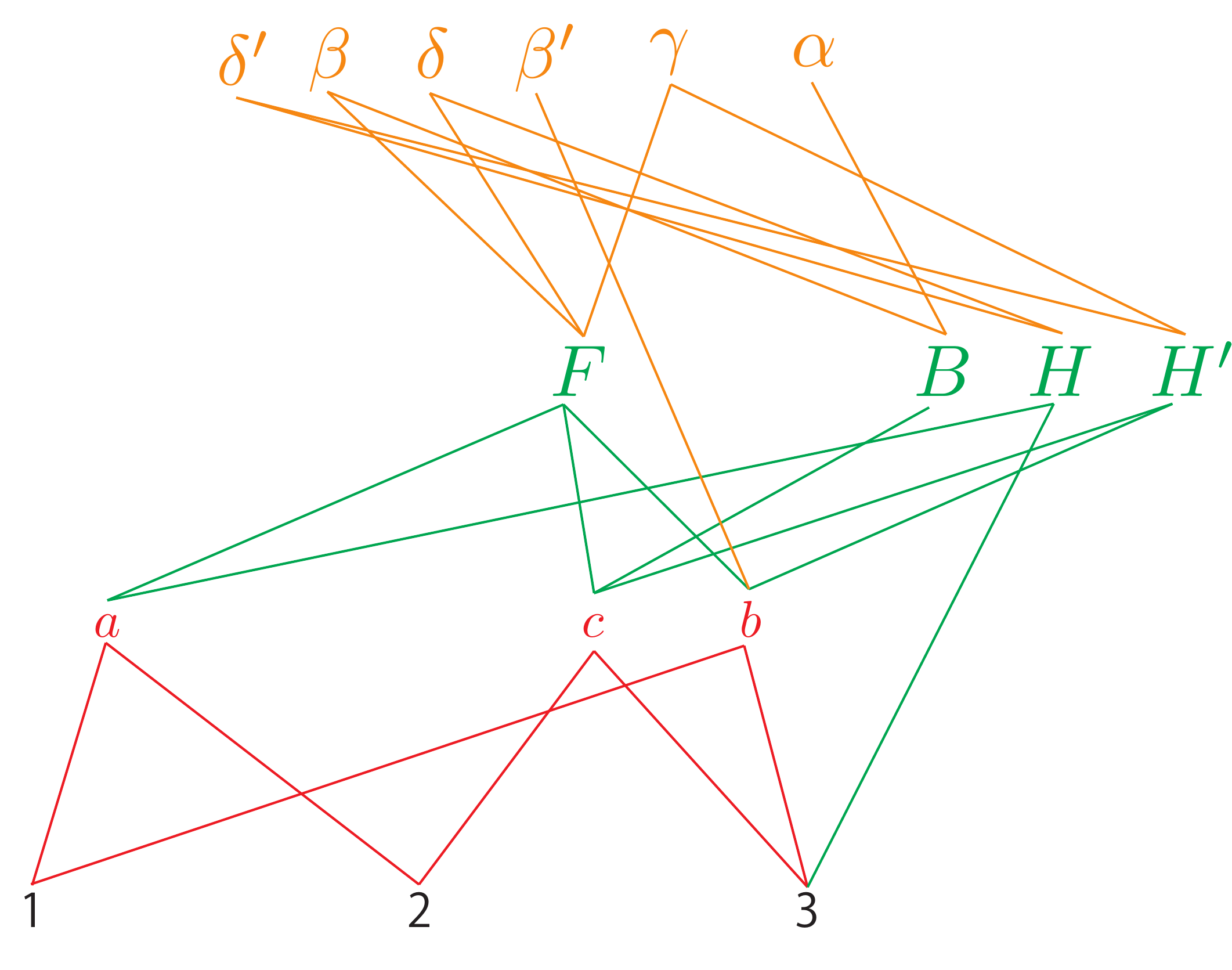}
\includegraphics[scale=0.15]{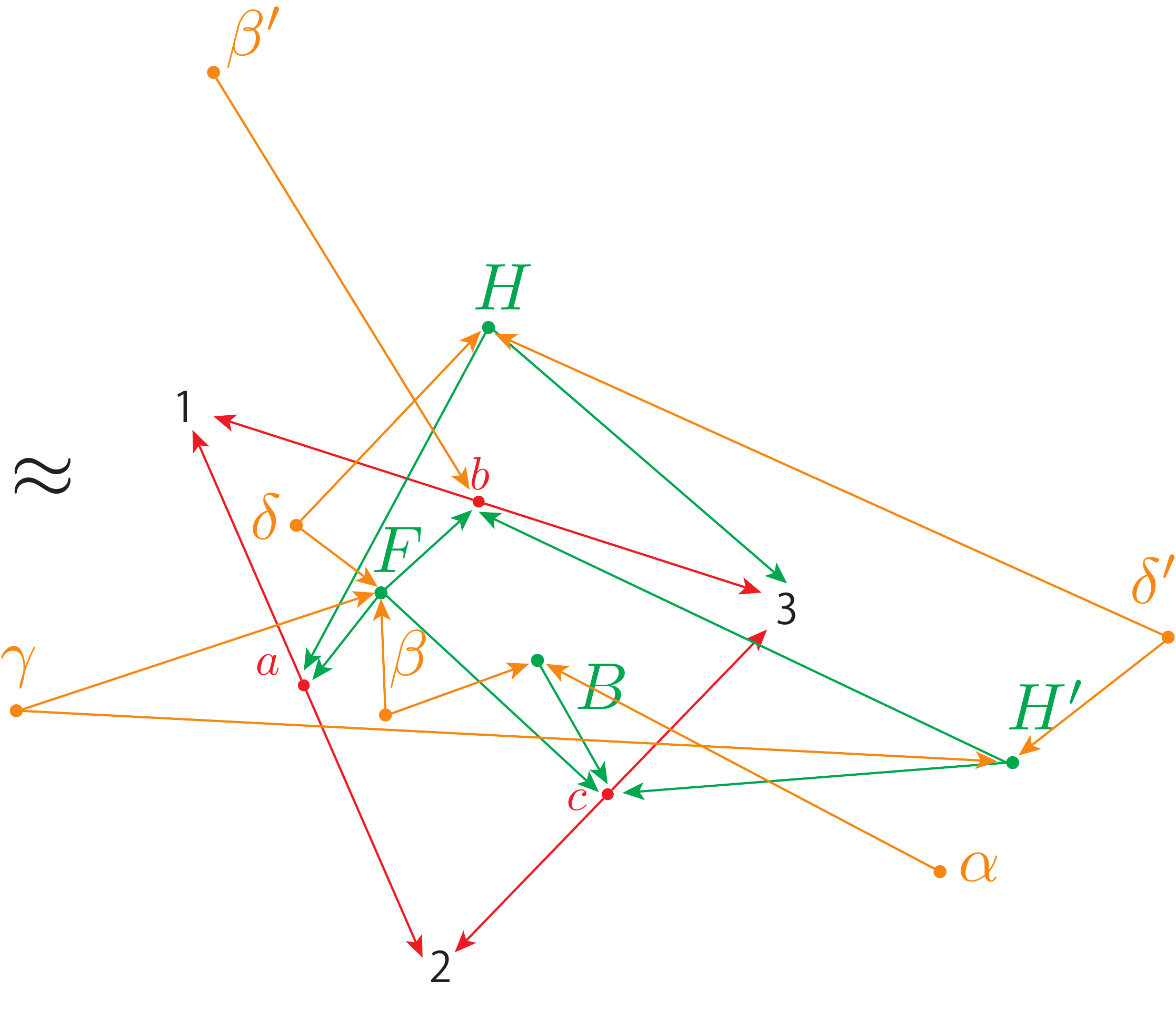}
\end{center}
\caption{The resulting Hesse diagram $X_3$ of the opposite order of the specialization preorder of the resulting diagram of $X_2$ by removing $G$, $C$, $G'$, and $A$.}
\label{fig:Hessian_grad_reduced_02}
\end{figure}

\begin{figure}
\begin{center}
\includegraphics[scale=0.15]{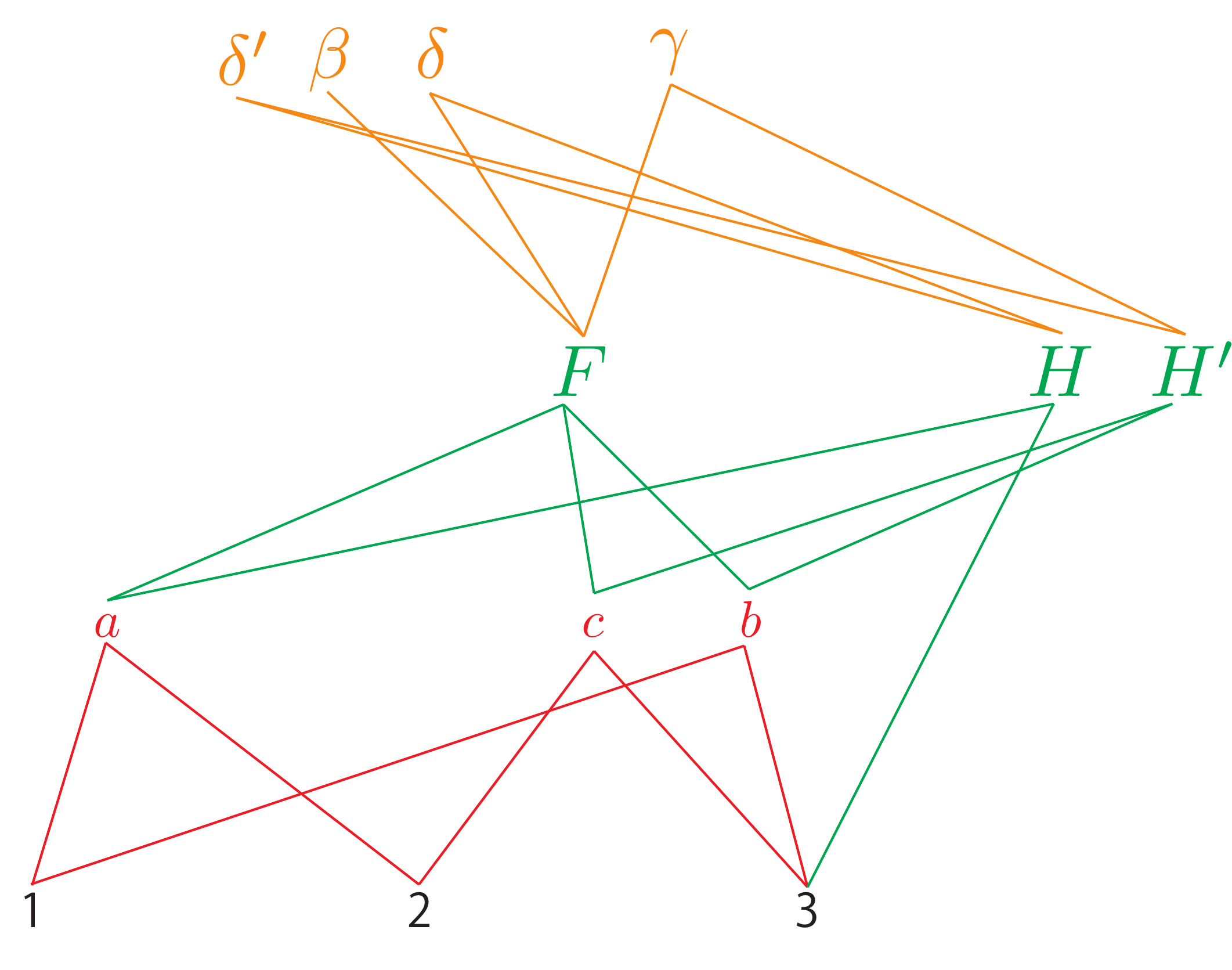}
\includegraphics[scale=0.15]{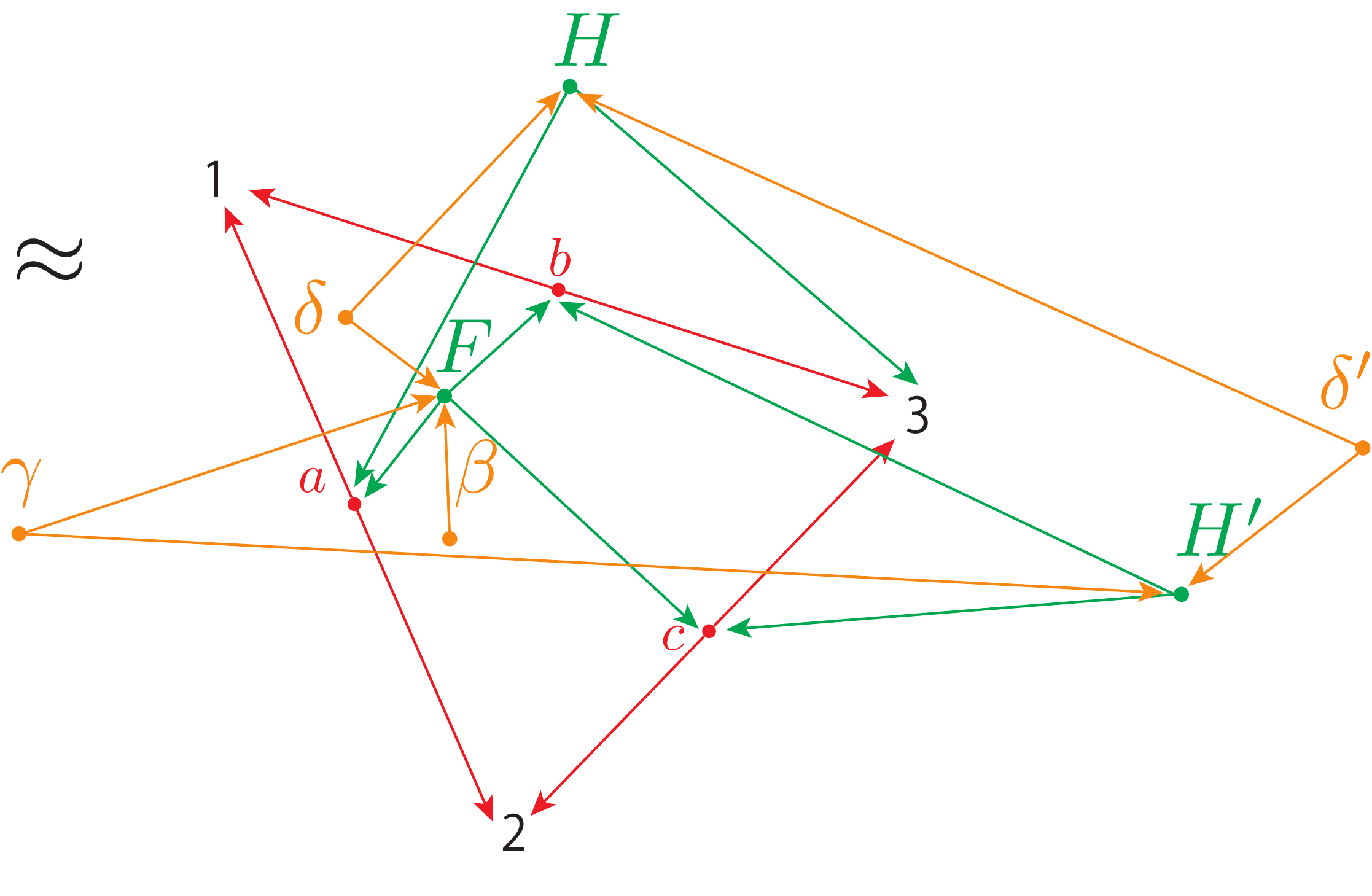}
\end{center}
\caption{The resulting Hesse diagram $X_4$ of the opposite order of the specialization preorder of the resulting diagram of $X_3$ by removing $\alpha$, $\beta'$, and $B$.}
\label{fig:Hessian_grad_reduced_03}
\end{figure}

\begin{figure}
\begin{center}
\includegraphics[scale=0.15]{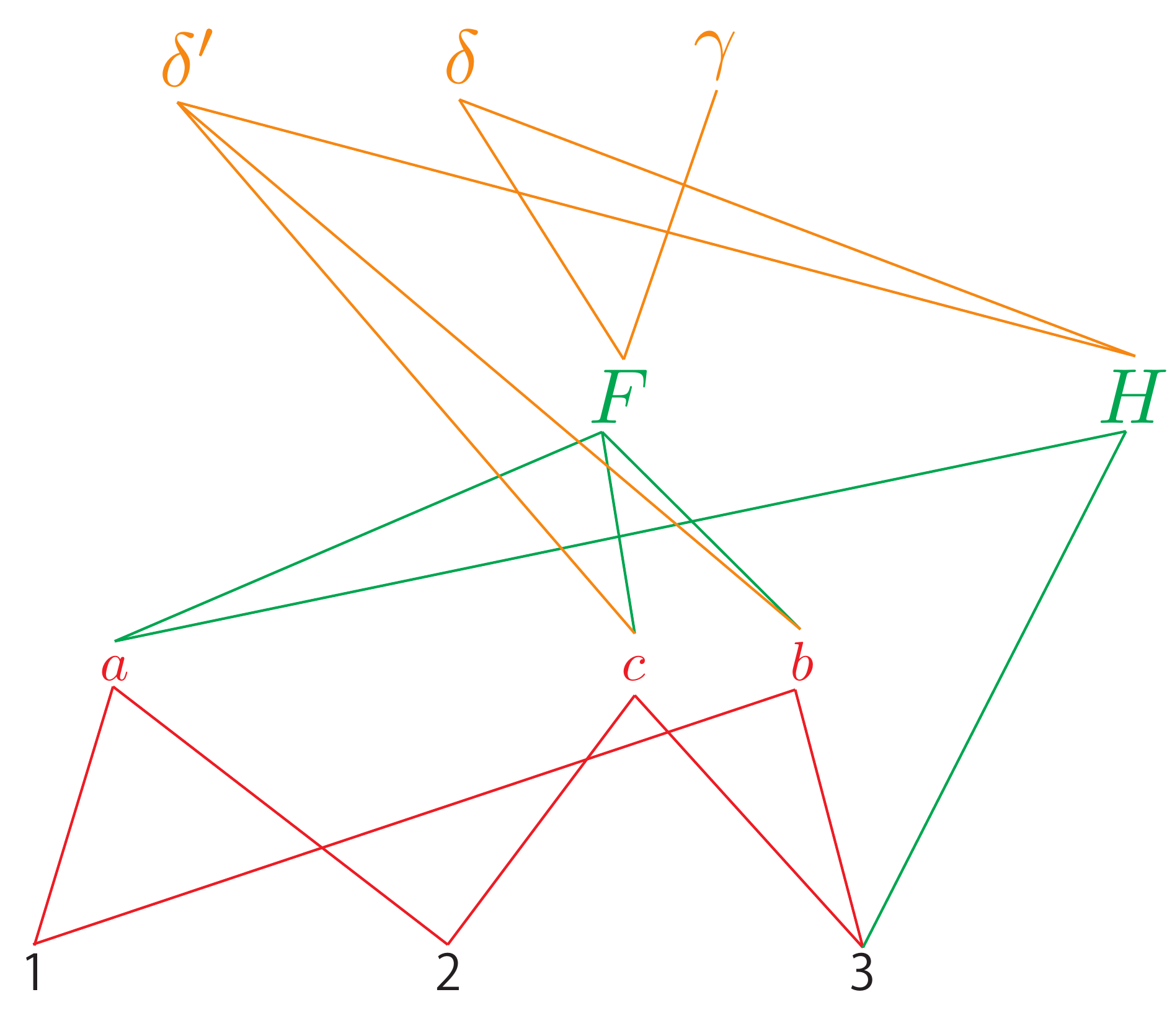}
\includegraphics[scale=0.15]{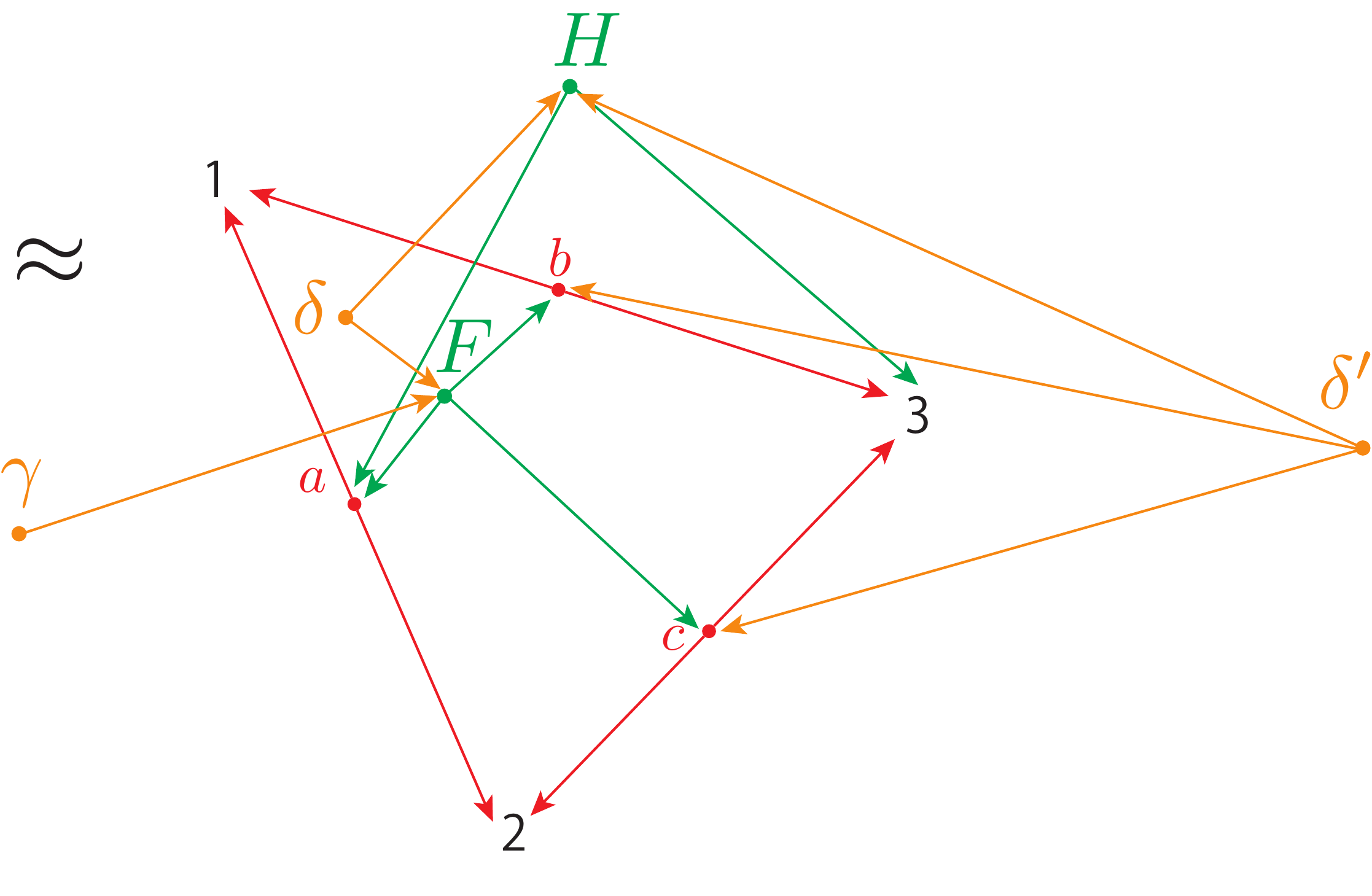}
\end{center}
\caption{The resulting Hesse diagram $X_6$ of the opposite order of the specialization preorder of the resulting diagram of $X_5$ by removing the down weak point $H'$.}
\label{fig:Hessian_grad_reduced_05}
\end{figure}

\begin{figure}
\begin{center}
\includegraphics[scale=0.15]{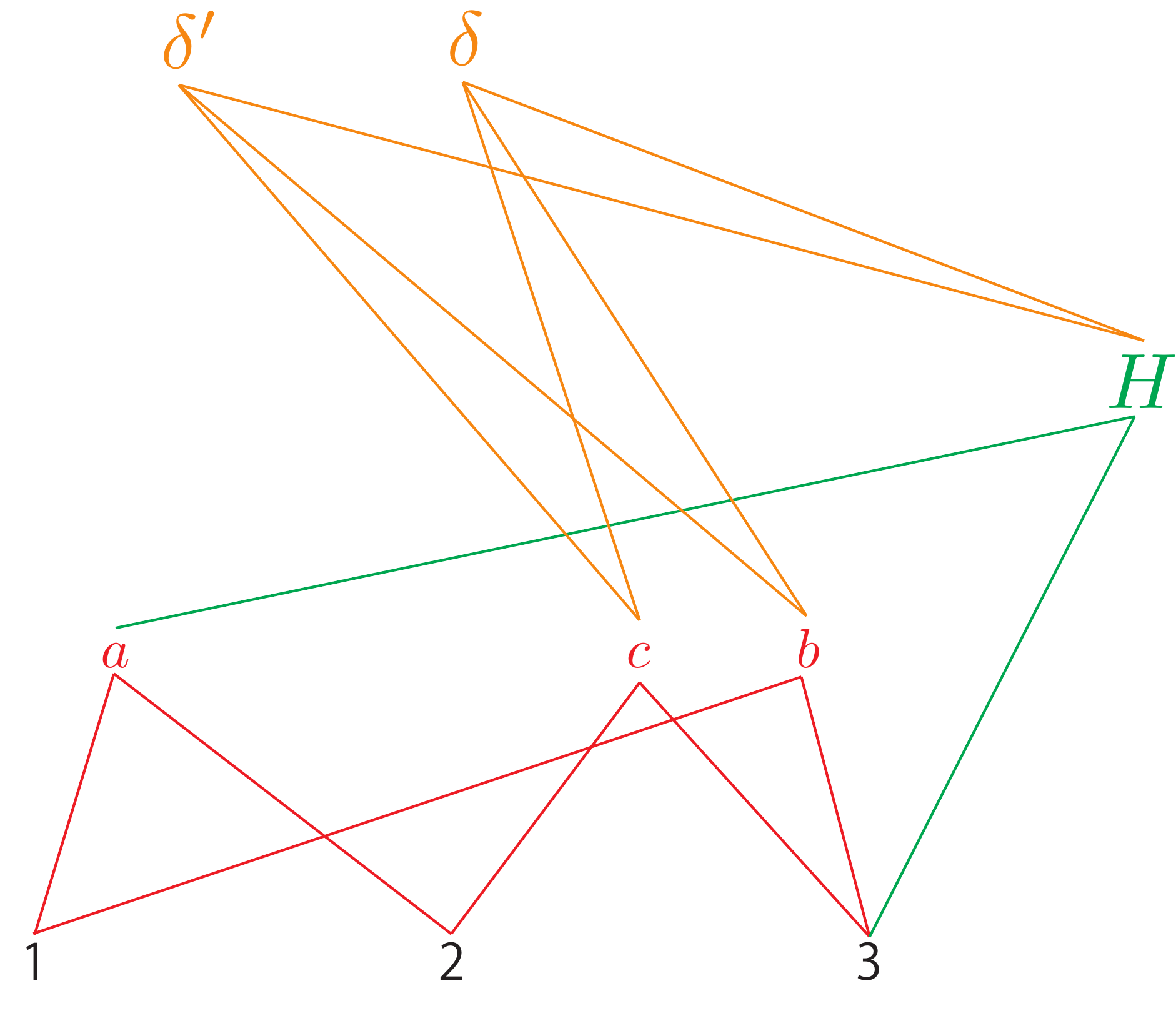}
\includegraphics[scale=0.15]{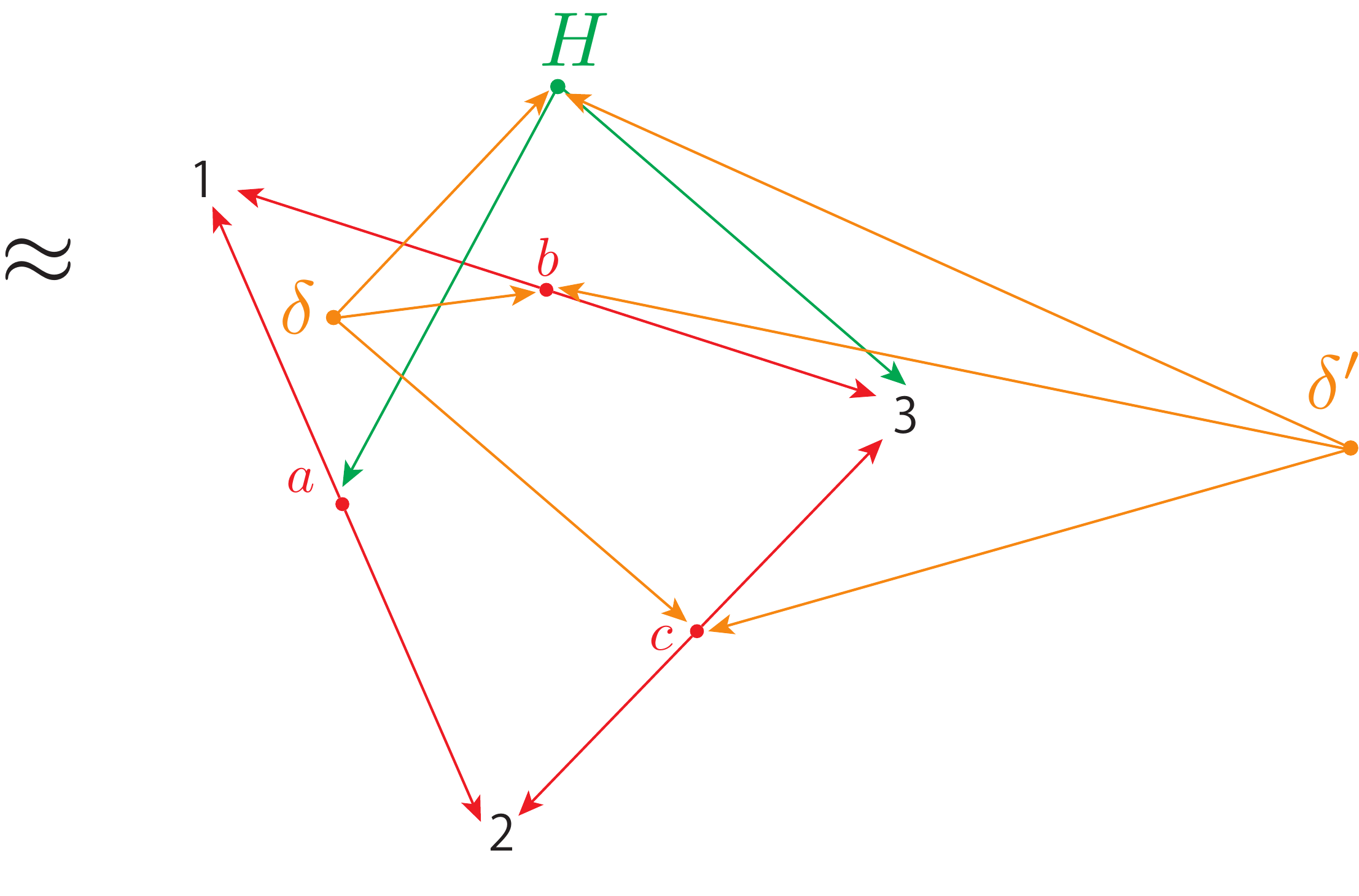}
\end{center}
\caption{The resulting Hesse diagram $X_7$ of the opposite order of the specialization preorder of the resulting diagram of $X_6$ by removing $\gamma$, $F$.}
\label{fig:Hessian_grad_reduced_06}
\end{figure}

\begin{proof}
The Hesse diagram of the opposite order of the specialization preorder of the finite $T_0$-space $[\mathcal{G}^r_{1,0,2,0 >-1}(\Sigma_{0,2})]$ is shown in Figure~\ref{fig:Hessian}.
Remove up beat points $E$, $C''$, $C'$, and $D$. 
Then we have the Hesse diagram as in Figure~\ref{fig:Hessian_grad_reduced_00}. 
Remove down beat points $d, e, h, f,$ and $g$. 
Then we have the Hesse diagram as in Figure~\ref{fig:Hessian_grad_reduced_01}. 
Remove down beat points $G$, $C$, $G'$, and $A$. 
Then we have the Hesse diagram as in Figure~\ref{fig:Hessian_grad_reduced_02}. 
Remove down beat points $\alpha$, $\beta'$, and $B$.
Then we have the Hesse diagram as in Figure~\ref{fig:Hessian_grad_reduced_03}. 
Remove a down beat point $\beta$.
Then we have the Hesse diagram as in Figure~\ref{fig:Hessian_grad_reduced_04} which is the core. 
Remove a down weak point $H'$.
Then we have the Hesse diagram as in Figure~\ref{fig:Hessian_grad_reduced_05}. 
Remove a down beat point $\gamma$.
Remove an up beat points $F$. 
Then we have the Hesse diagram as in Figure~\ref{fig:Hessian_grad_reduced_06}. 
Remove an up beat points $a$. 
Then we have the Hesse diagram as in Figure~\ref{fig:Hessian_grad_reduced_07}. 

This means that $[\mathcal{G}^r_{1,0,2,0 >-1}(\Sigma_{0,2})]$ is weak homotopy equivalent to the finite space whose Hessian is shown in Figure~\ref{fig:Hessian_grad_reduced_07}. 
\end{proof}

\begin{theorem}\label{lem:007}
The finite $T_0$-space $[\mathcal{G}^r_{1,0,2,0 >-1}(\Sigma_{0,2})]$ is weakly homotopy equivalent  to a bouquet $\mathbb{S}^2 \vee \mathbb{S}^2$ of two two-dimensional spheres. 
\end{theorem}

\begin{figure}
\begin{center}
\includegraphics[scale=0.195]{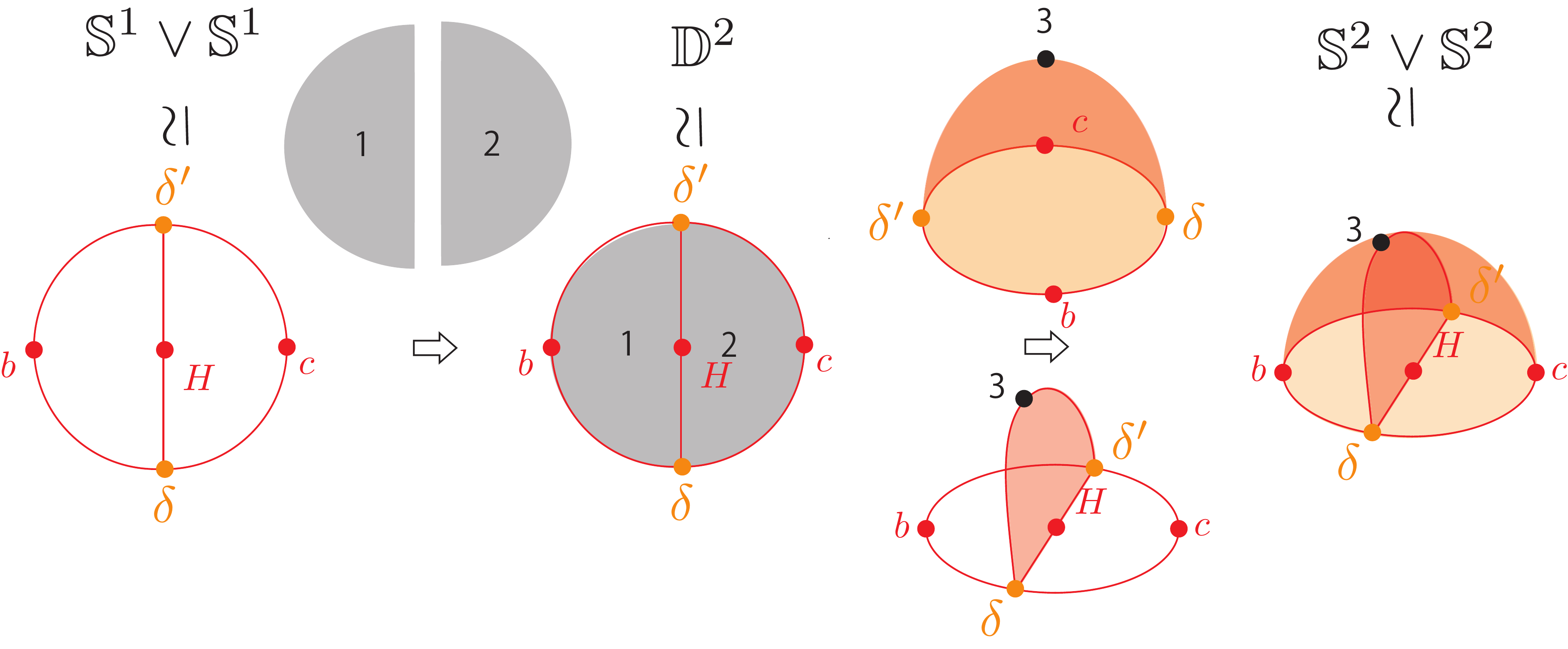}
\end{center}
\caption{Left, a bouquet $\mathbb{S}^1 \vee \mathbb{S}^1$ of two one-dimensional spheres; middle, a disk; right, a bouquet $\mathbb{S}^2 \vee \mathbb{S}^2$ of two two-dimensional spheres.}
\label{fig:order_cpx02}
\end{figure}

\begin{proof}
The order complex $\mathcal{K}(X')$ of the topological space $X' = \{ \delta, \delta', H, b, c \}$ is homotopy equivalent to a bouquet $\mathbb{S}^1 \vee \mathbb{S}^1$ of two one-dimensional spheres as on the left of  Figure~\ref{fig:order_cpx02}, and one of the finite topological space $X'' = \{ \delta, \delta', H, b, c, 1, 2 \}$ is homotopy equivalent to a closed disk as on the middle of  Figure~\ref{fig:order_cpx02}. 
Adding three disks, we obtain the order complex $\mathcal{K}(X_8)$ of generated by the minimal finite space $X_8 = \{ \delta, \delta', H, b, c, 1, 2, 3 \}$ is homotopy equivalent to a bouquet $\mathbb{S}^2 \vee \mathbb{S}^2$ as on the right of Figure~\ref{fig:order_cpx02}. 

From \cite[Proposition~3.3]{barmak2008simple}, the inclusion $X - \{x\} \to X$ into a finite $T_0$-space for any weak beat point $x \in X$ is weak homotopy equivalent. 
Therefore, the subspace $[\mathcal{G}^r_{1,0,2,0 >-1}(\Sigma_{0,2})]$ is weak homotopy equivalent to $X_8$. 
By \cite[Theorem~2]{mccord1966singular}, the order complex $\mathcal{K}(X_8)$ of the finite space $X_8$ is weak homotopy equivalent to the finite space $X_8$. 
This implies that the subspace $[\mathcal{G}^r_{1,0,2,0 >-1}(\Sigma_{0,2})]$ is weak homotopy equivalent to the bouquet $\mathbb{S}^2 \vee \mathbb{S}^2$ of two two-dimensional spheres. 
\end{proof}

The previous lemma implies Proposition~\ref{prop:5.3} and so Theorem~\ref{re_prop:5.3}. 
The author would like to ask the following questions. 
\begin{question}
Is every connected component $[\mathcal{G}^r_{k_{-,1},k_{-,1}, k_{+,1}, k_{+,2}, > -1}(S)]$ simply connected? 
\end{question}

\begin{question}
Does the space of topological equivalence classes of gradient flows on a three-dimensional {\rm(}or higher dimensional{\rm)} manifold have non-contractible connected components, under the non-existence of creations and annihilations of singular points? 
\end{question}

\section{Hierarchical structures of the space of Morse-Smale-like flows}\label{sec:07}

In this section, we show a similar statement of Theorem~\ref{prop:codimension} for ``gradient flow with limit cycles''. 

\subsection{Morse-Smale-like flows on surfaces}

We recall Morse-Smale-like flows to describe gradient flows and  ``gradient flows with limit cycles'' as follows. 
\begin{definition}
A flow four is {\bf Morse-Smale-like} \cite{kibkalo2022topological} if it satisfies the following conditions: 
\\
{\rm(1)} Each recurrent orbit is closed. 
\\
{\rm(2)} There are at most finitely many limit cycles. 
\\
{\rm(3)} Each singular point is finitely sectored.
\\
{\rm(4)} The set of non-recurrent points is open dense in $S$. 
\end{definition}

Notice that quasi-regular Morse-Smale-like flows are quasi-Morse-Smale (i.e. the resulting flows obtained from quasi-regular gradient flows by replacing singular points with limit cycles and pasting limit cycles) (see \cite{kibkalo2022topological} for definition details). 
This means that quasi-regular Morse-Smale-like flows can be considered as ``gradient flows with limit cycles''. 

\subsubsection{Generic non-Morse-Smale flows on surfaces}

Recall that the set of gradient $C^1$-flows is not open in the set of $C^1$-flows because singular points need not non-degenerate.
Moreover, topologically hyperbolic limit cycles for $C^1$-flows can be bifurcated into topologically non-hyperbolic limit cycles.
However, forbidding the existence of fake multi-saddles and fake limit cycles, fixing the sum of indices of repelling singular points (i.e. sources and $\partial$-sources), the sum of indices of attracting singular points (i.e. sinks and $\partial$-sinks), and the number of limit cycles, we topologically characterize codimension $k$ flows of Morse-Smale-like flows. 
%
%

\subsection{Combinatorial structure of the space of Morse-Smale-like flows}

Recall that a flow on a compact surface is a gradient flow with finitely many singular points if and only if the flow is Morse-Smale-like without elliptic sectors or non-trivial circuits \cite[Theorem~B]{kibkalo2022topological}.
Therefore, Morse-Smale-like flows without elliptic sectors can be considered as ``gradient flow with limit circuit'' and with finitely many singular points.
In this subsection, we describe the topology of the space of ``gradient flow with limit cycles'' under the non-existence of creations and annihilations of singular points and physical boundaries. 
%
Note that any Morse-Smale-like flows without elliptic sectors are quasi-regular.

\subsubsection{Codimensions of quasi-regular Morse-Smale-like flows}
Let $w$ be a quasi-regular Morse-Smale-like flow on a surface $T$. 
The {\bf codimension} of a $k$-saddle $y$ ($k \geq 1$) is $2(k-1) = -2(1+ \mathrm{ind}(y))$. 
The {\bf codimension} of a $\partial$-$(l/2)$-saddle $y$ ($l \geq 1$) is $l-1 = -2(1/2 + \mathrm{ind}(y))$. 
Equivalently, the codimension of a $\partial$-$k$-saddle $y$ ($k \geq 1/2$) is $2k-1 = -2(1/2 + \mathrm{ind}(y))$. 
\begin{definition}
The codimension $\bm{\mathrm{codim}_{\mathrm{m}}(w)} \in \Z_{\geq 0}$ of $w$ with respect to multiplicity is defined as the sum of the codimensions of multi-saddles of $w$. 
\end{definition}

\begin{definition}
The codimension $\bm{\mathrm{codim}_{\mathrm{h}}(w)} \in \Z_{\geq 0}$ with respect to heteroclinicity is defined as the number of multi-saddle separatrices outside of the boundary $\partial T$. 
\end{definition}

\begin{definition}
The codimension $\bm{\mathrm{codim}_{\mathrm{f}}(w)} \in \Z_{\geq 0}$ with respect to fakeness is defined as the sum of the number of fake multi-saddles, fake parabolic sectors, and fake limit cycles. 
\end{definition}

The  {\bf codimension} $\bm{\mathrm{codim}(w)}$ of the quasi-regular Morse-Smale-like flow $w$ is defined as follows: 
\[
\mathrm{codim}(w) := \mathrm{codim}_{\mathrm{m}}(w) + \mathrm{codim}_{\mathrm{h}}(w) + \mathrm{codim}_{\mathrm{f}}(w)
\]

\subsubsection{Subspaces}

For any $r \in \mathbb{Z}_{\geq 0} \sqcup \{ \infty \}$, denote by $\mathcal{Q}_{**}^r(S)$ the set of quasi-regular Morse-Smale-like flows without fake saddles or fake limit cycles. 
For any $k_-,k_+ \in \mathbb{Z}_{\geq 0}$, $r \in \mathbb{Z}_{\geq 0} \sqcup \{ \infty \}$, and $l \in \mathbb{Z}_{\geq 0}$, denote by $\mathcal{Q}_{k_-/2,k_+/2,l,**}^r(S)$ the set of quasi-regular Morse-Smale-like flows without fake saddles or fake limit cycles whose sum of indices of attracting (resp. repelling) singular points (i.e. sources and $\partial$-sources (resp. sinks and $\partial$-sinks)) is $k_-/2$ (resp. $k_+/2$) and whose numbers of limit circuits are $l$.
From \cite[Theorem~F]{kibkalo2022topological}, any flow in $\mathcal{Q}_{k_-/2,k_+/2,l,**}^r(S)$ is quasi-Morse-Smale (i.e. resulting flows obtained from quasi-regular gradient flows by replacing singular points with limit cycles and pasting limit cycles).

For any $k_-, k_+,l \in \mathbb{Z}_{\geq 0}$, denote by $\mathcal{Q}_{k_-/2,k_+/2,l,>-1}^r(S) \subseteq \mathcal{Q}_{k_-/2,k_+/2,l,**}^r(S)$ the set of flows in $\mathcal{Q}_{k_-/2,k_+/2,l,**}^r(S)$ without non-periodic circuits. 
Notice that the subspace $\mathcal{Q}_{k_-/2,k_+/2,l,>-1}^r(S)$ is a connected component of $\mathcal{Q}_{**}^r(S)$ under the non-existence of creations and annihilations of singular points, limit cycles, and boundary components.
Then, we have the following observation. 

\begin{lemma}
For any $k_-, k_+,l \in \mathbb{Z}_{\geq 0}$, the subspace $\mathcal{Q}_{k_-/2,k_+/2,l,>-1}^r(S)$ is the set of Morse-Smale-like flows without elliptic sectors, fake saddles, fake limit cycles, or non-periodic circuits such that the sum of indices of attracting {\rm(resp.} repelling{\rm)} singular points {\rm(i.e.} sources and $\partial$-sources {\rm(resp.} sinks and $\partial$-sinks{\rm))} is $k_-/2$ {\rm(resp.} $k_+/2${\rm)} and whose numbers of limit circuits are $l$.
\end{lemma}

By \cite[Theorem~B]{kibkalo2022topological}, any flow in $\mathcal{Q}_{k_-/2,k_+/2,0,>-1}^r(S)$ is gradient as follows. 

\begin{lemma}\label{lem:corr_msl_grad}
For any $k_-, k_+ \in \mathbb{Z}_{\geq 0}$, we have the following equality:  
\[
\mathcal{Q}^r_{k_-/2, k_+/2,0, > -1}(S) = \mathcal{G}^r_{k_-/2, k_+/2, >-1}(S)
\]
\end{lemma}

%
%

\subsubsection{Codimension of Morse-Smale-like flows}

For any non-negative integer $q \in \mathbb{Z}_{\geq 0}$, denote by $\bm{\mathcal{Q}_{k_-/2,k_+/2,l, q}(S)} \subseteq \mathcal{Q}_{k_-/2,k_+/2,l,>-1}^r(S)$ the subspace of such Morse-Smale-like flows whose codimensions are $q$ and define $\bm{\mathcal{Q}_{k_-/2,k_+/2,l, >q}(S)} := \bigcup_{q' > q}\mathcal{Q}_{k_-/2,k_+/2,l, q'}(S)$. 
Notice that $\mathcal{Q}_{k_-/2,k_+/2,l, >q}(S)$ is the subspace of such Morse-Smale-like flows whose codimensions are more than $q$. 
%
Therefore, Theorem~\ref{re_prop:5.3} and Lemma~\ref{lem:corr_msl_grad} imply the following statement. 

\begin{corollary}\label{cor:noncontractible}
For any $r \in \mathbb{Z}_{> 0} \sqcup \{ \infty \}$ and any integers $k_-, k_+,l \in \Z_{\geq 0}$, there is the quotient space of $\mathcal{Q}^r_{k_-/2, k_+/2,l, > -1}(S)$ by topologically equivalence classes which has a non-contractible connected component. 
\end{corollary}

%
We have the following hierarchical structure. 

\begin{theorem}\label{prop:codimension_ms}
The following statements holds for any $r \in \mathbb{Z}_{> 0} \sqcup \{ \infty \}$, any integers $k_-, k_+,l \in \Z_{\geq 0}$ and any $q \in \Z_{\geq -1}$: 
\\
{\rm(1)} The subspace $\mathcal{Q}^r_{k_-/2, k_+/2,l, q+1}(S)$ is open dense in the space $\mathcal{Q}^r_{k_-/2, k_+/2,l, > q}(S)$. 
\\
{\rm(2)} The subspace $\mathcal{Q}^r_{k_-/2, k_+/2,l, q+1}(S)$ consists of structurally stable in the space $\mathcal{Q}^r_{k_-/2, k_+/2,l, > q}(S)$. 
\end{theorem}

Theorem~\ref{prop:codimension} and the previous theorem generalize  \cite[Lemma~7.7]{kibkalo2022topological}\footnote{Though there are typos in \cite[Lemma~7.7]{kibkalo2022topological}, the correctly modified statement is contained in Theorem~\ref{prop:codimension} and Theorem~\ref{prop:codimension_ms}.}. 
To show the previous theorem, we observe the following statement.

\begin{lemma}\label{lem:limit_cycle}
For any $r \in \mathbb{Z}_{> 0} \sqcup \{ \infty \}$, any integers $k_-, k_+,l \in \Z_{\geq 0}$ and any $q \in \Z_{\geq -1}$, let $v$ be a flow in $\mathcal{Q}^r_{k_-/2, k_+/2,l, >-1}(S)$ and $A_1, \ldots , A_l$ basins of attracting or repelling limit cycles each of whose boundaries $\partial A_i$ consists of two closed transversals $T_{-,i}$ and $T_{+,i}$. 
Then there is a $C^r$-\nbd $\mathcal{U} \subset \mathcal{Q}^r_{k_-/2, k_+/2,l, >-1}(S)$ of $v$ such that, for any $w \in \mathcal{U}$, the annuli $A_i$ are basins of the attracting or repelling limit cycles with respect to $w$ and the boundary components $T_{-,i}$ and $T_{+,i}$ are two closed transversals with respect to $w$. 
\end{lemma}

\begin{proof}
Let $L_1, \ldots , L_l$ be the limit cycles of $v$. 
By the non-existence of fake limit cycles, the limit cycles $L_i$ are topologically hyperbolic (i.e. attracting or repelling). 
Then 
$A_i \cap \mathop{\mathrm{Sing}}(v) = \emptyset$. 

Fix any $i \in \{ 1, \ldots , l \}$. 
By time reversion if necessary, we may assume that $L_i$ is attracting. 
Since the transversality of $C^r$ flows is an open condition, there is a $C^r$-\nbd $\mathcal{U}_i \subset \mathcal{Q}^r_{k_-/2, k_+/2,l, >-1}(S)$ of $v$ such that $\bigcup_{w \in \mathcal{U}_i} \mathop{\mathrm{Sing}}(w) \cap A_i = \emptyset$ and $\bigcup_{w \in \mathcal{U}_i} \bigcup_{x \in A_i} \omega(x) \subset A_i$, and that the loops $T_{-,i}$ and $T_{+,i}$ are closed transversals with respect to any flow $w \in \mathcal{U}_i$. 
From the generalization of the Poincar\'e-Bendixson theorem for a flow with finitely many singular points (cf.  \cite[Theorem~2.6.1]{nikolaev1999flows}), the $\omega$-limit set of each point of $A_i$ contained in $A_i$ is a limit cycle. 

Put $\mathcal{U} := \bigcap_{i=1}^l \mathcal{U}_i$. 
For any $w \in \mathcal{U}$, since $w$ contains exactly $l$ limit cycles, any annuli $A_i$ contains exactly one limit cycle, which is topologically hyperbolic, such that $A_i$ is a basin of the attracting or repelling limit cycle with respect to $w$. 
%
%
%
%
\end{proof}

We demonstrate Theorem~\ref{prop:codimension_ms} as follows. 

\begin{proof}[Proof of Thorem~\ref{prop:codimension_ms}]
Fix a flow $v \in \mathcal{Q}^r_{k_-/2, k_+/2,l, >-1}(S)$ and basins of attracting or repelling limit cycles $A_1, \ldots , A_l$ as in Lemma~\ref{lem:limit_cycle}. 
From \cite[Theorem~F]{kibkalo2022topological}, the flow $v$ is quasi-Morse-Smale (i.e. resulting flows obtained from quasi-regular gradient flows by replacing attracting or repelling singular points with limit cycles and pasting limit cycles).
In other words, cutting limit cycles of $v$, taking the end completion, and collapsing the boundary components into singletons, we obtain the resulting flow $v'$ which is a gradient flow in $\mathcal{G}^r_{k_-/2, k_+/2, >-1}(S)$ up to topologically equivalent via a homeomorphism whose restriction to the complement of the union of small basins of the singletons is identical, such that $\mathrm{codim}(v') = \mathrm{codim}(v)$. 
Therefore, any perturbation of $v$ whose restriction to ${\bigsqcup_{i=1}^l A_i}$ is locally topologically equivalent to the restriction $v\vert_{\bigsqcup_{i=1}^l A_i}$ can be reduced into a perturbation of $v'$ in $\mathcal{G}^r_{k_-/2, k_+/2, >-1}(S)$. 
From Lemma~\ref{lem:limit_cycle}, the restriction to $\bigsqcup_{i=1}^l A_i$ of any small perturbation of $v$ is locally topologically equivalent to the restriction $v\vert_{\bigsqcup_{i=1}^l A_i}$. 
Therefore, the proof of Lemma~\ref{lem:codimension_open} implies the assertion. 
\end{proof}

Theorem~\ref{prop:codimension_ms} implies the following finite abstract cell complex structures. 

\begin{theorem}\label{main:03}
The following statements hold for any $r \in \mathbb{Z}_{> 0} \sqcup \{ \infty \}$, any integers $k_-, k_+,l \in \Z_{\geq 0}$ and any $q \in \Z_{\geq -1}$:
\\
{\rm(1)} The quotient space $ [\mathcal{Q}^r_{k_-/2, k_+/2,l, q+1}(S)]$ of $\mathcal{Q}^r_{k_-/2, k_+/2,l, q+1}(S)$ by topological equivalence classes is a finite $T_0$-space and is an abstract cell complex with respect to the opposite order of the specialization preorder and the codimension. 
\\
{\rm(2)} $\overline{[\mathcal{Q}^r_{k_-/2, k_+/2,l, q+1}(S)]_{q+1}} = [\mathcal{Q}^r_{k_-/2, k_+/2,l, >-1}(S)]_{\geq q+1} = [\mathcal{Q}^r_{k_-/2, k_+/2,l, > q}(S)]$. 
\end{theorem}

\section{Final remarks}

We state the following future work and the following open question. 

\subsection{Transition via non-periodic limit circuits between gradient and non-gradient Morse-Smale-like flows} 

Notice that the ``generic'' intermediate flow with a non-periodic limit circuit, which is not a gradient flow, must appear between Morse flows and non-Morse Morse-Smale-like flows in $\bigsqcup_{l \geq 0} \mathcal{Q}_{k_-/2,k_+/2,l,>-1}^r(S)$. 
To state such transitions via saddle connections containing non-periodic limit circuits, we need the ``codimension'' of the multi-saddle connection diagram for flows in the space $\bigsqcup_{l \geq 0} \mathcal{Q}_{k_-/2,k_+/2,l,**}^r(S)$ like the ``codimension'' of the multi-saddle connection diagram for Hamiltonian flows \cite{yokoyama2021combinatorial}. 
We will report such transitions in the near future. 


\subsection{Simple connectivity of the space of ``gradient flows with limit circuits''}

Though we found a non-contractible component for gradient flows, the author does not know whether non-simply-connected connected components exist for gradient flows and Morse-Smale-like flows. 
In other words, one would like to know about the following question.


\begin{question}
Does the space of topological equivalence classes of gradient flows {\rm(}or more generally Morse-Smale-like flows{\rm)} on compact manifolds have non-simply-connected connected components under the non-existence of creations and annihilations of singular points and boundary components? 
\end{question}

\appendix

\section{Invariance of sums of indices of isolated singular points under small perturbations}

We demonstrate the continuity of the partial map $V^v_B \colon \R \times \partial B \rightharpoonup \mathbb{S}^{n-1}$ defined in Definition~\ref{def:conti}, the invariance of indices of singular points, the local stability of basins of repellers (i.e. Lemma~\ref{lem:inv_sink/source}), and the invariance of isolated singular points (Lemma~\ref{lem:inv_index} and Corollary~\ref{cor:inv_index_02}).

\subsection{Continuity of the partial map $V^v_B \colon \R \times \partial B \rightharpoonup \mathbb{S}^{n-1}$}

We have the following statement. 

\begin{lemma}\label{lem:inv_index03}
Let $v$ be a $C^1$-flow on $\R^n$, $B \subset \R^n$ a closed disk with $\mathop{\mathrm{Sing}}(v) \cap \partial B = \emptyset$, and $\mathbb{S}^{n-1} \subset \R^n$ the unit sphere.  
The partial map $V^v_B \colon \R \times \partial B \rightharpoonup \mathbb{S}^{n-1}$ defined in Definition~\ref{def:conti} is a continuous map. 
\end{lemma}

\begin{proof}
By the definition of $C^1$-flow, the derivative $\left. \dfrac{\partial v}{\partial t} \cdot \left\vert \dfrac{\partial v}{\partial t} \right\vert^{-1} \right|_{t=0} \colon \partial B \to \mathbb{S}^{n-1}$ is continuous, and the difference $g_v := v - \mathrm{proj_2} \colon (\R - \{ 0 \}) \times \partial B \to \R^n$ (i.e. $g_v(t,x) = v(t,x) - \mathrm{proj_2}(t,x) = v(t,x) - x$) is continuous. 
By the non-existence of singular points in the complement $\partial B$, we have that $0 \notin \mathrm{Im}(g_v)$ and so that $g_v/ \vert g_v \vert  \colon (\R - \{ 0 \}) \times \partial B \to \mathbb{S}^{n-1}$ is well-defined and continuous.  
Then $V^v_B(t,x) = g_v(t,x)/ \vert g_v(t,x) \vert$ for any $t > 0$, and $V^v_B(t,x) = - g_v(t,x)/ \vert g_v(t,x) \vert$ for any $t < 0$. 
By $\mathop{\mathrm{Sing}}(v) \cap \partial B = \emptyset$, since the flow $v$ is of $C^1$, we obtain the following equality: 
\[
\lim_{t \to 0} \dfrac{g_v(t,x)}{t} = \lim_{t \to 0} \dfrac{v(t,x) - v(0,x)}{t} = \dfrac{\partial v(0,x)}{\partial t} \neq 0
\]
Therefore, we have the following equality: 
\[
\begin{split}
\lim_{t \to 0\pm} V^v_B(t,x) = \lim_{t \to 0\pm} \pm \dfrac{g_v(t,x)}{ \vert g_v(t,x) \vert } 
&= \lim_{t \to 0\pm} \dfrac{g_v(t,x)}{t}  \left| \dfrac{g_v(t,x)}{t} \right|^{-1}
\\
&= \lim_{t \to 0\pm} \dfrac{v(t,x) - v(0,x)}{t} \left| \dfrac{v(t,x) - v(0,x)}{t} \right|^{-1} \\
&= \dfrac{\partial v(0,x)}{\partial t} \cdot \left| \dfrac{\partial v(0,x)}{\partial t} \right|^{-1}
\end{split}
\]
This means that $V^v_B \colon \R \times \partial B \to \mathbb{S}^{n-1}$ is a continuous map. 
\end{proof}

In the rest of this subsection, we state some properites for indices to show the local stability of basins of repellers and the invariance of isolated singular points. 

\subsubsection{Invariance of the sum of indices of continuations of isolated singular points of $C^1$-flows on manifolds}

%
Moreover, we have the following statement. 

\begin{corollary}\label{cor:inv_index_sum}
Let $w$ be a $C^1$-flow on a manifold and for any closed disk $B$ whose boundary contains no singular points and which contains at most finitely many singular points. 
Then the partial map $V^w_B \colon \R \times \partial B \rightharpoonup \mathbb{S}^{n-1}$ defined in Definition~\ref{def:conti} is continuous and the domain of $V^w_B$ is a \nbd of $\{ 0 \} \times \partial B$ in $\R \times \partial B$. 
\end{corollary}

\begin{proof}
By the proof of the previous lemma, it suffices to show that there is a positive number $T$ such that $w(t, x) \neq x$ for any $(t,x) \in ((-T, T) - \{ 0 \}) \times \partial B)$.  
Assume that there are no such $T$. 
Then there are sequences $\{ x_n \}_{n=1}^\infty$ of points $x_n \in \partial B$ and $\{ T_n \}_{n=1}^\infty$ of positive numbers $T_n$ with $\lim_{n \to \infty} T_n = 0$ and $w(T_n,x_n) = x_n$. 
Then $x_n \in \mathop{\mathrm{Per}}(w)$ for any $n \in \Z_{> 0}$. 
By the compactness of $\partial B$, there is a point $x_\infty \in \partial B$ with $x_\infty = \lim_{n \to \infty} x_n$. 
Since $\lim_{n \to \infty} T_n = 0$ and $T_n >0$ for any $n \in \Z_{> 0}$, for any $m \in \Z_{>0}$, there is a large number $N_m \in \Z_{>0}$ and a sequence $\{ k_{m,n} \}_{n = N_m}^\infty$ of natural numbers $k_{m,n} \in \Z_{>0}$ such that $k_{m,n} T_n \in [1/(m+1),1/m]$. 

Fix any natural number $m \in \Z_{>0}$. 
Take an increasing subsequence $\{ i_{m,n} \}_{n=1}^\infty$ of the sequence $\{ n \}_{n=1}^\infty$ of natural numbers such that the sequence $\{ k_{m,i_{m,n}} T_{i_{m,n}} \}_{n = 1}^\infty$ of numbers converges into a number $S_m \in [1/(m+1),1/m]$. 
Since $\lim_{n \to \infty} x_{i_{m,n}} = x_\infty$, $w(k_{m,i_{m,n}} T_{i_{m,n}},x_{i_{m,n}}) = x_{i_{m,n}}$, and $\lim_{n \to \infty} k_{m,i_{m,n}} T_{i_{m,n}} = S_m$, the continuity of $w$ implies that $w(S_m,x_\infty) = \lim_{n \to \infty} w(k_{m,i_{m,n}} T_{i_{m,n}},x_{i_{m,n}}) = \lim_{n \to \infty} x_{i_{m,n}} = x_\infty$. 

Then $x_\infty$ is a closed point whose periods are $S_n \in [1/(n+1),1/n]$ for any $n \in \Z_{>0}$. 
Since $x_\infty$ has arbitrarily small periods, the point $x_\infty \in \partial B$ is a singular point, which contradicts $\partial B$ contains no singular points.
\end{proof}

We have the following statement. 

\begin{lemma}\label{lem:eq_index}
Let $w$ be a flow on a surface $S$, $\nu$ a loop transverse at all but finitely many tangencies, and $B$ the closed disk whose boundary is $\nu$ with $ \vert B \cap \mathop{\mathrm{Sing}}(w) \vert  < \infty$. 
Then 
\[
\dfrac{n_{i} - n_{o} + 2}{2} = \deg (V^w_B(0, \cdot)) = \sum_{y \in B \cap \mathop{\mathrm{Sing}}(w)} \mathrm{ind}_w(y)
\]
where $n_i$ {\rm(resp.} $n_o${\rm)} is the number of inner {\rm(resp.} outer{\rm)} tangencies of $\nu$ with respect to $B$. 
\end{lemma}

\begin{proof}
From $ \vert B \cap \mathop{\mathrm{Sing}}(w) \vert   < \infty$ and $\mathop{\mathrm{Sing}}(w) \cap \partial B = \emptyset$, there is a finitely many closed disks $B_1, \ldots, B_k$ each of which contains at most one singular point such that the intersection $B_i \cap B_j$ for any $i \neq j$ are either the empty set or closed intervals, and that $B = \bigcup_{i=1}^k B_i$ and $\mathop{\mathrm{Sing}}(w) \cap \bigcup_{i=1}^k \partial B_i = \emptyset$. 
Then the image of $V^w_B(0, \cdot)$ is the connected sum of  the images of $V^w_{B_1}(0, \cdot), \ldots, V^w_{B_k}(0, \cdot)$ such that each $B_i$ is an isolated \nbd of a singular point or contains no singular point. 
Lemma~\ref{lem:equiv_index} and Corollary~\ref{cor:inv_index_sum} imply that 
\[
\begin{split}
\deg (V^w_B(0, \cdot)) &= \sum_{i=1}^k \deg (V^w_{B_i}(0, \cdot)) 
\\
&= \sum_{i=1}^k \sum_{y \in B_i \cap \mathop{\mathrm{Sing}}(w)} \mathrm{ind}_w(y) = \sum_{y \in B \cap \mathop{\mathrm{Sing}}(w)} \mathrm{ind}_w(y)= \dfrac{n_{i} - n_{o} + 2}{2}
\end{split}
\]
because $\deg (V^w_{B_i}(0, \cdot)) = \mathrm{ind}_w(y_i)$ if $B_i$ contains exactly one singular point $y_i$, and $\deg (V^w_{B_i}(0, \cdot)) = 0$ if $B_i$ contains no singular points. 
\end{proof}

\subsubsection{Degrees of attracting/repelling subsets}

Recall that $A \Subset C$ means $\overline{A} \subset \mathop{\mathrm{int}} C$. 
We observe the following equality. 

\begin{lemma}\label{lem:inv_attracting}
Let $v$ be a $C^1$-flow on $\R^n$, $B$ closed disk, and $V^v_B \colon \R \times \partial B \rightharpoonup \mathbb{S}^{n-1}$ the partial map defined in Definition~\ref{def:conti}. 
If there is a nonzero number $T \in \R$ such that $v(T, B) \Subset B$, then $V^v_B$ is a continuous map satisfying the following equality: 
\[
\deg (V^v_B(0, \cdot)) = \deg (V^v_B(T, \cdot)) = \begin{cases}
 1 & \text{if } T < 0 \\
 (-1)^n & \text{if } T > 0
\end{cases}
\]
\end{lemma}

\begin{proof}
Suppose that $T<0$. 
Since $v(T, v(nT, B)) = v((n+1)T, B)$, we have that 
\[
B \Supset v(T, B) \Supset v(2T, B) \Supset v(3T, B) \Supset \cdots
\]
and so that $\partial B$ contains no closed points. 
Lemma~\ref{lem:inv_index03} implies that $V^v_B$ is a continuous map with $\deg (V^v_B(0, \cdot)) = \deg (V^v_B(T, \cdot))$. 
%
Then there are disjoint open \nbds $U_i$ and $U_o$ of $v(T, B)$ and $\R^n - \mathop{\mathrm{int}} B$ respectively. 
By coordinate transformation, we may assume that $\mathop{\mathrm{int}} B$ contains the origin. 
Since $v(T, B) \subset \mathop{\mathrm{int}} B$, collapsing $v(T, B)$ into a singleton in $\mathop{\mathrm{int}} B$, there is a homotopy $H \colon [0,1] \times \R^n \to \R^n$ from $H(0, \cdot ) = v(T, \cdot)$ such that $H(1,v(T, B)) = \{ 0\}$ and $H(t, \cdot) \vert _{U_o} = v(T, \cdot) \vert _{U_o}$ for any $t \in [0,1]$. 
Since $H(t, \cdot) \vert _{U_o} = v(T, \cdot) \vert _{U_o}$ for any $t \in [0,1]$, we have $\bigcup_{t \in [0,1]} H(t, \partial v(T, B)) \cap \partial B = \emptyset$. 
Then $V_H \colon [0,1] \times \partial B \to \mathbb{S}^{n-1}$ defined as follows is well-defined and continuous: 
\[
V_H(t,x) :=  - \dfrac{H(t,x) - x}{ \vert H(t,x) - x \vert } 
\]
By definition, we obtain that 
\[
\begin{split}
V_H(0,x) &=  - \dfrac{H(0,x) - x}{ \vert H(0,x) - x \vert } =  - \dfrac{v(T,x) - x}{ \vert v(T,x) - x \vert } = V^v_B(T, x)
\\
V_H(1,x) &=  - \dfrac{H(1,x) - x}{ \vert H(1,x) - x \vert } = \dfrac{x}{\vert x \vert} 
\end{split}
\]
for any $x \in \partial B$. 
The continuity of $V_H$ implies that $\deg (V^v_B(T, \cdot)) = \deg (V_H(0,\cdot)) = \deg (V_H(1,\cdot)) = 1$. 

Suppose that $T>0$. 
The assertion for $T<0$ implies that  $\deg (V^{-v}_B(t, \cdot)) = 1$. 
By time reversion, we have $\deg (V^v_B(t, \cdot)) = (-1)^n \deg (V^{-v}_B(t, \cdot)) = (-1)^n$. 
\end{proof}

%

\subsection{Local stability of basins of repellers under small perturbations}\label{sec:loc_stb}

We have the following persistence. 

\begin{lemma}\label{lem:inv_transversal}
Let $v$ be a flow on a surface $S$ and $x \in S$ a source or sink. 
For any closed disk $B$ which is an isolated \nbd of $x$, there is a $C^0$-\nbd $\mathcal{U}$ of $v$ satisfying the following statement: 
\\
{\rm(1)} For any $w \in \mathcal{U}$, there is a closed transversal $\mu$ to $w$ in the interior $\mathop{\mathrm{int}}B$ of $B$ such that $\mathop{\mathrm{Sing}}(w) \cap B = \mathop{\mathrm{Sing}}(w) \cap D_\mu$, where $D_\mu \subset \mathop{\mathrm{int}}B$ is the closed disk whose boundary is $\mu$.  
\end{lemma}

\begin{proof}
Since $x \in S$ is a source {\rm(resp.} sink{\rm)}, there is a positive (resp. negative) number $T$ such that $v(T,B) \Subset B$.  
Take an open \nbd $U \Subset B$ of $v(T,B)$ and define $\mathcal{U}:= \{ w : \text{flow on } S \mid w(T,B) \subset U \}$. 
Fix any $w \in \mathcal{U}$. 
Applying \cite[Lemma~3.1]{kibkalo2022topological} to the loop $\partial v(T,B)$, there is a loop $\mu \Subset B$ that is smooth and is transverse except finitely many tangencies with respect to $w$.
By Gutierrez's smoothing theorem~\cite{gutierrez1978structural}, we may assume that $w$ is of $C^1$. 
From $w(T,B) \Subset B$, Lemma~\ref{lem:inv_attracting} implies that the index of $x$ with respect to $w$ is one. 
Since the Euler characteristic of a closed disk is one, by Lemma~\ref{lem:eq_index}, the Poincar{\'e}-Hopf theorem implies that the number of inner tangencies and one of outer tangency of $w$ are the same. 
By a small perturbation of $\mu$ as in Figure~\ref{elim_inner_tangency} if necessary, we may assume that $\mu$ has no inner tangencies and so no tangencies.
\begin{figure}
\begin{center}
\includegraphics[scale=0.4]{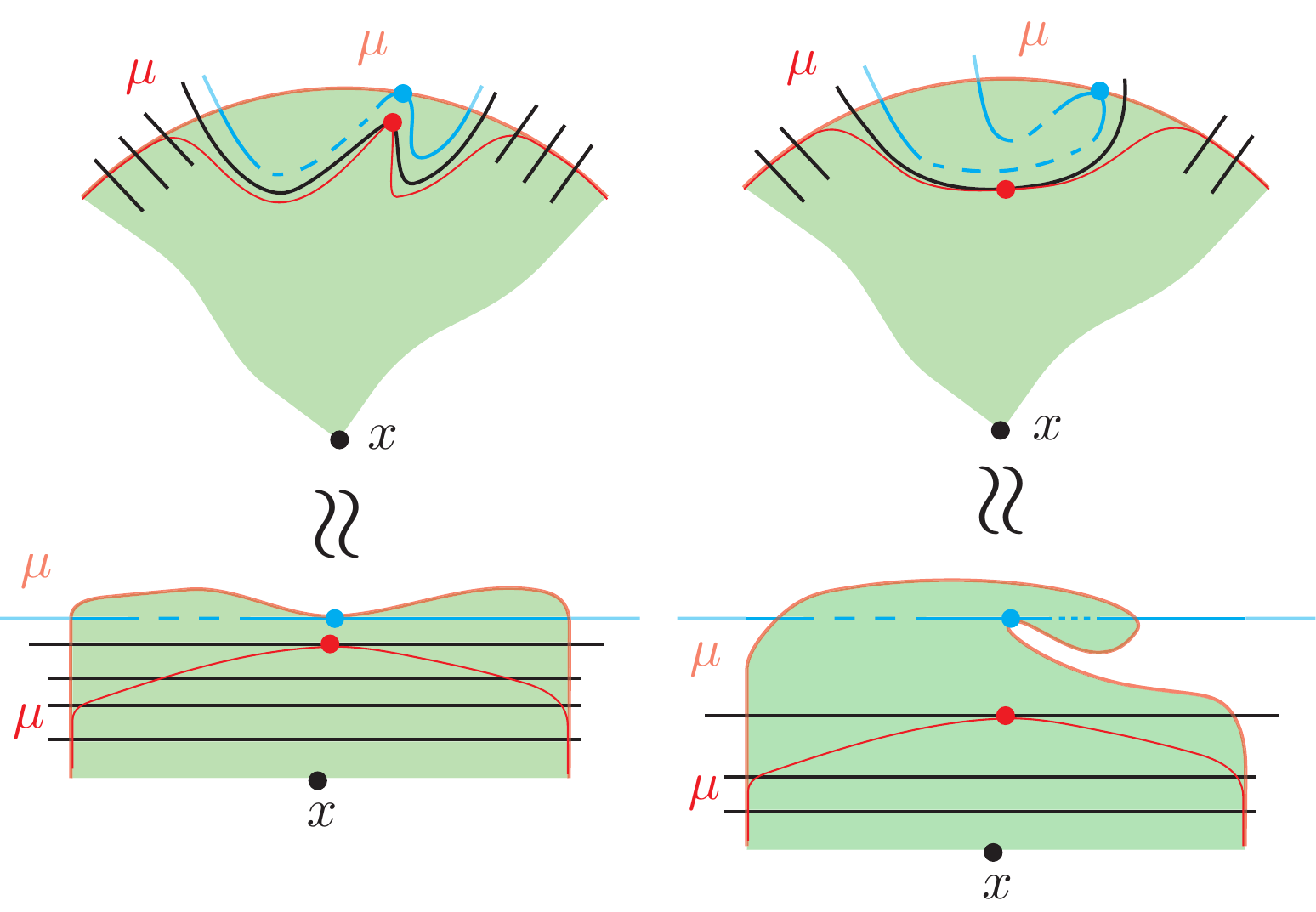}
\end{center}
\caption{Deformations to eliminate inner tangencies.}
\label{elim_inner_tangency}
\end{figure}
This means that the resulting loop $\mu$ is a closed transversal of $w$ in $\mathop{\mathrm{int}}B$. 
\end{proof}

We demonstrate the local stability of basins of repellers as follows. 

\begin{proof}[Proof of Lemma~\ref{lem:inv_sink/source}
]
Let $v$ be a flow on a compact surface $S$ and $x \in S$ a sink or a source. 
By time reversion if necessary, we may assume that $x \in S$ is a source. 
Since the index of the isolated singular point is invariant under topological equivalence, by Gutierrez's smoothing theorem~\cite{gutierrez1978structural}, we may assume that $v$ is of $C^1$. 

Put $K := \{ T \} \times B$ and take an open \nbd $U \Subset B$ of $v(T, B)$. 
Consider an open $C^0$-\nbd $\mathcal{U}(K, U) := \{ w : \text{flow on } S \mid v(K) \subset U \}$ of $v$.
For any flow $w \in \mathcal{U}(K, U)$, we have that $w(T,B) = w(K) \subset U \Subset B$. 
Moreover, for any flow $w \in \mathcal{U}(K, U)$, since $w(T,\cdot)$ is a homeomorphism on $S$, we obtain that  $w(2T,B) = w(T,w(T,B)) \Subset w(T,B)$ and so $w(nT,B) \Subset w(T,B) \Subset B$ for any $n \in \Z_{>0}$. 
Similarly, put $K' := \{ T \} \times B'$ and take an open \nbd $U' \Subset B'$ of $v(T, B')$. 
Then $w(K') = w(T,B') \subset U' \Subset B'$ for any flow $w \in \mathcal{U}(K', U')$. 
Put $K'_0 := [0,T] \times B'$ and take an open \nbd $U_0 \Subset B$ of $v([0,T], B')$. 
Consider an open $C^0$-\nbd $\mathcal{U}(K'_0, U_0) := \{ w : \text{flow on } S \mid v(K'_0) \subset U_0 \}$ of $v$.
For any flow $w \in\mathcal{U}(K, U) \cap  \mathcal{U}(K'_0, U_0)$, we have that $w([0,T],B') = w(K'_0) \subset U_0 \Subset B$ and so that $w(nT,w([0,T],B')) \Subset w(nT,B) \Subset w(T,B)$ for any $n \in \Z_{>0}$. 
For any flow $w \in \mathcal{U}(K, U) \cap \mathcal{U}(K'_0, U_0)$, since $w(\R_{\geq 0},B') = \bigcup_{n \in \Z_{\geq 0}} w(nT,w([0,T],B')) $, we have the following inequality: 
\[
\begin{split}
w(\R_{\geq 0},B') &= w([0,T],B') \cup \bigcup_{n \in \Z_{> 0}} w({nT},w([0,T],B'))
\\
& \subseteq w([0,T],B') \cup w(T,B) \Subset B
\end{split}
\]
%
Lemma~\ref{lem:inv_index02} implies that there is an open $C^0$-\nbd $\mathcal{U} \subseteq \mathcal{U}(K, U) \cap \mathcal{U}(K', U') \cap \mathcal{U}(K'_0, U_0)$ such that $\mathop{\mathrm{Sing}}(w) \cap \partial B = \emptyset$ for any $w \in \mathcal{U}$. 
By Lemma~\ref{lem:inv_transversal}, for any flow $w \in \mathcal{U}$, there is a closed disk $D_\mu \subset \mathop{\mathrm{int}} B$ whose boundary is a closed transversal $\mu$ of $w$ such that $\mathop{\mathrm{Sing}}(w) \cap B = \mathop{\mathrm{Sing}}(w) \cap D_\mu$. 

Fix any flow $w \in \mathcal{U}$ with $ \vert B \cap \mathop{\mathrm{Sing}}(w) \vert   < \infty$ if exists. 
Since the index of the isolated singular point is invariant under topological equivalence, by Gutierrez's smoothing theorem~\cite{gutierrez1978structural}, we may assume that $w$ is of $C^1$. 
By Lemma~\ref{lem:inv_attracting}, we obtain $\deg (V^w_B(T, \cdot)) = 1$, where the continuous map $V^w_B \colon \R \times \partial B \to \mathbb{S}^{1}$ is defined in Definition~\ref{def:conti}. 
Since $B$ contains only finitely many singular points of $w$, we have that $\mathrm{ind}_v(x) = 1 =  \deg (V^w_B(0, \cdot))$. 
By $w(T,B) \Subset B$, we have $\mathop{\mathrm{Sing}}(w) \cap \partial B = \emptyset$. 
From $ \vert B \cap \mathop{\mathrm{Sing}}(w) \vert   < \infty$, Lemma~\ref{lem:eq_index} implies that $1 = \mathrm{ind}_v(x) = \deg (V^w_B(0, \cdot)) =  \sum_{y \in B \cap \mathop{\mathrm{Sing}}(w)} \mathrm{ind}_w(y)$. 
\end{proof}

\subsection{Invariance of isolated singular points}

\subsubsection{Existence of closed transversals}

We observe the following statement. 

\begin{lemma}\label{lem:inv_multi-saddle_01}
Let $v$ be a flow on a compact surface $S$ and $\nu$ a closed transversal. 
For any \nbd $U_{\nu}$ of $\nu$, then there is a small closed transverse annulus $A \Subset U_{\nu}$ of $v$ which is a \nbd of $\nu$ and there is a $C^0$-\nbd $\mathcal{U}$ of $v$ satisfying the following statements for any flow $w \in \mathcal{U}$:  
\\
{\rm(1)} 
There is a closed transversal $\nu' \Subset A$ parallel to $\nu$ in $A$ with respect to $w$ whose flow direction is the same as the flow direction of $\nu$ with respect to $v$ in the orientable closed annulus $A$, and there is a closed annulus $A' \Subset A$ which is a \nbd of $\nu'$ such that the restriction $w \vert _{A'}$ is a transverse annulus. 
\\
{\rm(2)} If there is a closed disk $B'$ whose boundary is a boundary component of $A$ such that the union $A \cup B'$ is a closed disk $B$, then there are a closed disk $D_w \Subset B$ with $\partial D_w \Subset A$ and $\mathop{\mathrm{Sing}}(w) \cap B = \mathop{\mathrm{Sing}}(w) \cap D_w$ such that $\partial D_w$ is a closed tranvsercal of $w$. 
\\
{\rm(3)} If the disk $B$ as above exists and satisfies $\vert B \cap \mathop{\mathrm{Sing}}(w) \vert  < \infty$, then $1 = \sum_{y \in B \cap \mathop{\mathrm{Sing}}(w)} \mathrm{ind}_w(y)$. 
\end{lemma}

\begin{proof}
Fix any \nbd $U_{\nu}$ of the closed transversal $\nu$. 
Choose any small $T \neq 0 \in \R$ with $v(t,\nu) \cap v(s,\nu) = \emptyset$ for any $t\neq s \in [-T,3T]$ 
 such that $A_0 := v([-T,3T],\nu)$ is a closed annulus and $A_1 := v([-T,T],\nu)$ is a closed annulus in $U_{\nu}$. 
Then the closed annulus $A := v([-T/2,T/2],\nu) \Subset A_1 \subset A_0$ contains $\nu$. 

Put $K := \{ T \} \times A_1$ and take an open \nbd $U \Subset A_0$ of $v(T, A_1)$.
Consider an open $C^0$-\nbd $\mathcal{U}(K, U) := \{ w : \text{flow on } S \mid w(K) \subset U \}$ of $v$.
Then $w(K) = w(T,A_1) \subset U \Subset A_0$ for any flow $w \in \mathcal{U}(K, U)$. 
Similarly, put $K' := \{ -T \} \times v(2T, A_1)$. 
Consider an open $C^0$-\nbd $\mathcal{U}(K', U) := \{ w : \text{flow on } S \mid w(K') \subset U \}$ of $v$.
Then $w(K') = w(T,A_1) \subset U \Subset A_0$ for any flow $w \in \mathcal{U}(K', U)$. 
In addition, write $K'' := \{ 2T \} \times A$. 
Since the loop $v(2T,A) = v([3T/2,5T/2],\nu) \Subset A_0 - A_1$ and the closed annulus $A_1 = v([-T,T],\nu)$ are disjoint closed, there is a closed annulus $V \Subset A_0 - A_1$ which is a \nbd of $v(2T,A)$. 
Consider an open $C^0$-\nbd $\mathcal{U}(K'', V) := \{ w : \text{flow on } S \mid w(K'') \subset V \}$.  
Lemma~\ref{lem:inv_index02} implies that there is an open $C^0$-\nbd $\mathcal{U} \subseteq \mathcal{U}(K, U) \cap \mathcal{U}(K', U) \cap \mathcal{U}(K'', V)$ such that $\mathop{\mathrm{Sing}}(w) \cap v([-T,3T],\nu)  = \emptyset$ for any $w \in \mathcal{U}$. 
In particular, since $A_0 = v([-T,3T],\nu)$, for any $w \in \mathcal{U}$, we have $\mathop{\mathrm{Sing}}(w) \cap A_0 = \emptyset$. 
Furthermore, there is a \nbd $U_{A_1}$ of $A_1$ with $U_{A_1} \cap (V \cup v([2T,3T],\nu)) = \emptyset$.

Fix any $w \in \mathcal{U}$. 
Then $w(T,A_1) \subset U \Subset A_0$, $\mathop{\mathrm{Sing}}(w) \cap A_0  = \emptyset$,  and $w(2T,A) \subset V \Subset A_0 - A_1$.
By Gutierrez's smoothing theorem~\cite{gutierrez1978structural}, we may assume that $w$ is of $C^1$. 
Applying \cite[Lemma~3.1]{kibkalo2022topological} to the loop $\nu \Subset A$, there is a loop $\nu' \Subset A \Subset A_1$ that is smooth and is transverse except finitely many tangencies with respect to $w$.
Then $w(2T,\nu') \Subset w(2T,A) \subset V \Subset A_0 - A_1$.
Let $A_\partial \subset A_0$ be the closure of the connected component of $A_0 - w(2T,\nu')$ contining $v(3T,\nu)$. 
Then $A_\partial \cap U_{A_1} = \emptyset$.


\begin{claim}\label{claim:0a}
By a small perturbation of $\nu'$ if necessary, we may assume that $\nu'$ is a closed transversal of $w$. 
\end{claim}


\begin{proof}[Proof of Claim~\ref{claim:0a}]
Multiplying a bump $\varphi$ function on $A_0$ with $\varphi^{-1}(0) = A_\partial$ and $U_{A_1} \subset \varphi^{-1}(1)$ to $dw(0,\cdot)/dt$, we have the resulting vector field $X' : = \varphi \cdot dw(0,\cdot)/dt$ on $A_0$
which satisfies that $A_0 \cap \mathop{\mathrm{Sing}}(X') = A_\partial$. 
Collapsing $A_\partial$ into a singleton $x_\partial$, the resulting surface $B$ from $A_0$ is a closed disk and  the resulting vector field $X''$ on $B$ from $X'$ satisfies that $X'' = X' = dw(0,\cdot)/dt$ on the annulus $A_1$ and $B \cap \mathop{\mathrm{Sing}}(X'') = \{ x_\partial \}$. 
This means that the subset $A_1$ can identified with the closed disk $B$ and  the restriction $w\vert_{A_1}$ can be identified with the restriction to $A_1 \subset B$ of a flow $w''$ on a closed surface containing the disk $B$. 

From $w(T,A_1) \Subset A_0$, since the annulus $A_0$ is an isolated \nbd of a singular point of $w''$, Lemma~\ref{lem:inv_attracting} implies that the index of $x$ with respect to $w''$ is one. 
Since the Euler characteristic of the closed disk $B$ is one, the Poincar{\'e}-Hopf theorem implies that the number of inner tangencies and one of outer tangency of $\nu'$ with respect to $w''$ and so $w$ are the same. 
By a small perturbation of $\nu'$ as in Figure~\ref{elim_inner_tangency} if necessary, we may assume that $\nu'$ has no inner tangencies and so no tangencies.
This means that the resulting loop $\nu' \Subset A \Subset A_1$ is a closed transversal of $w$. 
\end{proof}

The previous claim implies assertion {\rm(2)}. 
By construction of $\nu'$, the closed transversal $\nu' \Subset A$ of $w$ is parallel to $\nu$ in $A$ with respect to $w$ whose flow direction is the same as the flow direction of $\nu$ with respect to $v$ in the orientable closed annulus $A$. 
From the flow box theorem (cf. \cite[Theorem~4.2.6, p.95]{markley2023flows}), there is a closed annulus $A' \Subset A \Subset A_1$ which is a \nbd of $\nu'$ such that the restriction $w \vert _{A'}$ is a transverse annulus. 
Thus assertion (1) holds. 
Lemma~\ref{lem:inv_sink/source} implies assertion {\rm(3)}. 
\end{proof}

\subsubsection{Existence of closed transverse arcs}

We state the existence of transverse arcs in the following way.

\begin{lemma}\label{lem:existence_transversal}
Let $v$ be a flow on a surface $S$, $C$ a closed arc whose boundary consists of two points $x,y$, and $B$ a closed disk containing no singular points. 
Suppose that there is a positive number $T$ satisfying the following conditions: 
\\
{\rm(1)}  $v([0,3T],x) \cap C = \{x \}$ and $v([0,3T],y) \cap C = \{ y \}$. 
\\
{\rm(2)} $C \sqcup v(T,C) \sqcup v(2T,C) \sqcup v(3T,C) \subset B$.
\\
{\rm(3)} There is a positive number $\varepsilon >$ such that $v(\sigma_i [0,T], C) \cap B_{\varepsilon}(v(3T,C)) = \emptyset$. 
\\
{\rm(4)} The boundary $\partial B$ is a loop $v((0,3T), x) \sqcup v((0,3T), y) \sqcup C \sqcup v(3T, C)$.

Then there is a closed transverse arc $\gamma \subset B - (v(T,C) \sqcup v(2T,C))$ between $v((T,3T), x)$ and $v((T,3T), y)$ of $v$ such that the intersection $\gamma \cap \partial B$ consists of two points in $v((T,3T), x) \sqcup v((T,3T), y)$.  
\end{lemma}

\begin{figure}
\begin{center}
\includegraphics[scale=0.6]{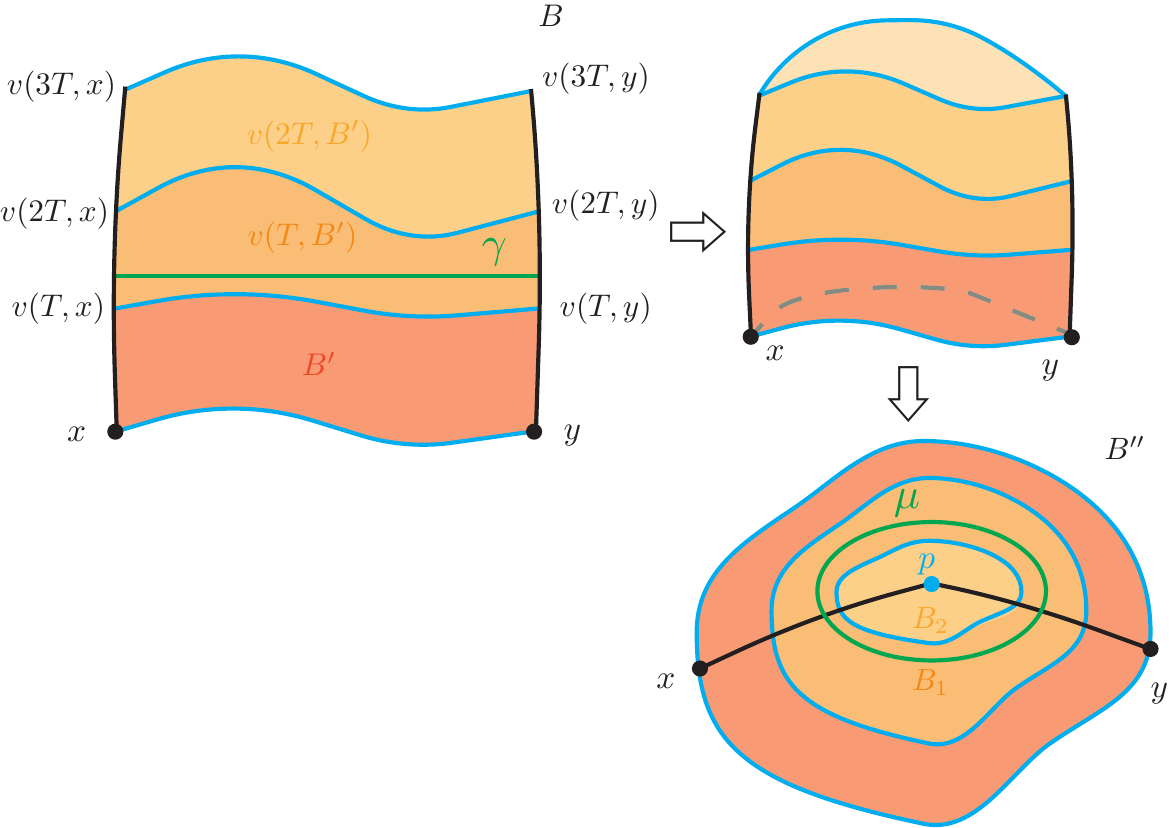}
\end{center}
\caption{Two closed disks and closed annulus.}
\label{floxbox_B}
\end{figure}

\begin{proof}
%
%
%
By Gutierrez's smoothing theorem~\cite{gutierrez1978structural}, we may assume that $v$ is of $C^1$. 
Considering a smooth bump function $\varphi \colon B \to [0,1]$ with $\varphi^{-1}(0) = v(3T,C)$ and $v([0,3T],C) \setminus B_{\varepsilon}(v(3T,C)) \subseteq \varphi^{-1}(1)$, let $w$ be the flow generated by the vector field $\varphi \dfrac{d v(0, \cdot)}{dt}$. 
Then the restriction $w \vert _{B - v(3T,C)}$ is locally topologically equivalent to $v \vert _{B - v(3T,C)}$. 
Since $v([0,3T],x) \cap C = \{x \}$ and $v([0,3T],y) \cap C = \{ y \}$, taking two copies of $B$ with $w$ and pasting along $v([0,3T], x)$ and $v([0,3T], y)$ respectively, we obtain the resulting space which is a closed annulus with the resulting flow $w'$. 
Collapsing the boundary component of the closed annulus corresponding to two copies of $v(3T,C)$ into a singleton $p$, we obtain a flow $w''$ on a closed disk $B''$ as in the bottom of Figure~\ref{floxbox_B}. 
Put $\widetilde{B}_1 :=v([T,3T],C)$. 
Let $B_1$ (resp. $B_2$) be the closed disk in $B''$ corresponding to two copies of $\widetilde{B}_1 = v([T,3T],C)$ (resp. $v([2T,3T],C)$) such that the boundary $\partial B_1$ (resp. $\partial B_2$) corresponds to the union of two copies of $v(T,C)$ (resp. $v(2T,C)$). 
There is a continuous mapping $h \colon \widetilde{B}_1 \to B_1$ via which the restriction $v \vert _{\mathop{\mathrm{int}} \widetilde{B}_1}$ is locally topologically equivalent to $w'' \vert _{h(\mathop{\mathrm{int}} \widetilde{B}_1)}$. 
By $v(\sigma_i [0,T], C) \subseteq v([0,3T],C) \setminus B_{\varepsilon}(v(3T,C)) \subseteq \varphi^{-1}(1)$, we have $w''(T,B_1) \Subset B_1$. 
Then the restriction $w'' \vert _{B_1}$ has the unique isolated singular point $p$. 
Lemma~\ref{lem:inv_attracting} implies that the index of $p$ is one. 
Applying \cite[Lemma~3.1]{kibkalo2022topological} to the loop in $B_1 - B_2$, there is a loop $\mu \Subset B_1 - B_2$ that is smooth and is transverse except finitely many tangencies with respect to $w''$.
Since the index of $p$ is one, the number of inner tangencies and one of outer tangency are same. 
By a small perturbation of $\mu$ as in Figure~\ref{elim_inner_tangency} if necessary, we may assume that $\mu$ has no inner tangencies and so no tangencies.
This means that the resulting loop $\mu$ is a closed transversal of $w''$ in $\mathop{\mathrm{int}}B_1$.
Therefore, the inverse image $\gamma' := h^{-1}(\mu) \subset \widetilde{B}_1 - (v(T,C) \sqcup v(3T,C)) \subset B - (v(T,C) \sqcup v(3T,C))$ is a closed transversal in $\mathop{\mathrm{int}}\widetilde{B}_1$. 
This means that the restriction $\gamma := \gamma' \cap B$ to $B$ of $\gamma'$ is desired. 
\end{proof}

\subsubsection{Local despription of loops and closed disks}

We have the following statement. 

\begin{lemma}\label{lem:inv_multi-saddle02}
Let $v$ be a flow on a compact surface $S$ and $\nu$ a loop transverse at all but finitely many tangencies.
For any \nbd $U_{\nu}$ of $\nu$, there is a small closed annulus $A_v \Subset U_{\nu}$ which is a \nbd of $\nu$ and there is a $C^0$-\nbd $\mathcal{U}$ of $v$ satisfying the following statements for any flow $w \in \mathcal{U}$:  
\\
{\rm(1)} There is a loop $\nu' \Subset A_v$ parallel to $\nu$ in $A_v$ with exactly the same number of inner tangencies and outer tangencies with respect to $w$ in the same order as $\nu$. 
Moreover, there is a closed annulus $A_w \Subset A_v$ with $\nu' \Subset A_w$ such that the restriction $w \vert _{A_w}$ is locally topologically equivalent to the restriction $v \vert _{A_v}$. 
\\
{\rm(2)} Moreover, if the loop $\nu$ bounds a closed disk $B$ such that $ \vert B \cap \mathop{\mathrm{Sing}}(w) \vert  < \infty$, then $\mathrm{ind}_v(x) = \sum_{y \in B \cap \mathop{\mathrm{Sing}}(w)} \mathrm{ind}_w(y)$. 
\end{lemma}

\begin{proof}
Fix any small \nbd $U_{\nu}$ of $\nu$. 
By Lemma~\ref{lem:inv_multi-saddle_01}, we may assume that $\nu$ is not any closed transversal. 
Denote by $z_1, z_2, \ldots , z_k$ the tangencies in $\nu$ for some even natural number $k > 0$. 
Put $z_0 =: z_k$.
Then there is a closed annulus $A' \Subset U_{\nu}$ which is a \nbd of $\nu$, contains no singular points, and is a finite union $k$ closed disks $V_1, \ldots , V_{k} =: V_0$ with $z_i \in V_i$ to which the restrictions are locally topologically equivalent to the restriction to a rectangle as in Figure~\ref{tang_sectors} such that $V_i$ intersects only $V_{i\pm1}$ and that the intersection $V_i \cap V_{i-1}$ is a trivial flow box and the intersection $\nu \cap V_i \cap V_{i-1}$ is a closed transverse arc. 
\begin{figure}
\begin{center}
\includegraphics[scale=0.9]{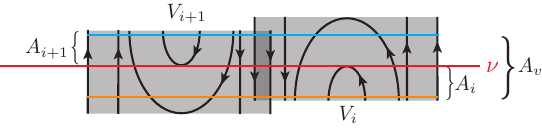}
\end{center}
\caption{Two serial rectangles each of which contains exactly one tangency of the loop $\nu$ with respect to $w$.}
\label{tang_sectors}
\end{figure} 
Choose a closed annulus $A_v \Subset A'$ which is a \nbd of $\nu$ and each of whose boundary components is a loop parallel to $\nu$ in $A'$ with exactly the same number of inner tangencies and outer tangencies in the same order as $\nu$. 
Denote by $A_i$ the closure of the connected component of $V_i \cap (A_v - \nu)$ which contains the orbit arc in $V_i$ containing $z_i$. 
Let $W_i$ be the connected component of $V_i \cap (A_v - \nu)$ which does not intersect $A_i$. 
Choose points $x_i \in W_i$ for any $i \in \{ 1, \ldots , k \}$ as in Figure~\ref{fig:perturbations_tang}. 
For any $i \in \{ 1, \ldots , k \}$, there are numbers $s_{i-}, s'_{i-} < 0$ and  $s_{i+}, s_{i+} > 0$ such that $v(s_{i\pm},x_i) \in \mathop{\mathrm{int}} A_i$ and $W_i \cap v([s_{i-},s_{i+}],x_i) = v((s'_{i-},s'_{i+}),x_i)$. 
There are small numbers $\varepsilon, \delta > 0$ such that the following compact subsets $K_{i,j}$ and open subsets $U_{i,j}$: 
\[
K_{i,1} := [s_{i-},s_{i+}] \times \{ x_i \}
\]
\[
U_{i,1} := B_{\varepsilon}(v(K_{i,1}))
\]
\[
K'_{i,2} := \{s_{i-}\} \times \{ x_i \}
\]
\[
K_{i,2} := [s_{i-}, s_{i-} + \delta] \times \{ x_i \}\]
\[
U_{i,2} := B_{\varepsilon}(v(K'_{i,2})) \Subset A_i
\]
\[
K'_{i,3} := \{s_{i+}\} \times \{ x_i \}
\]
\[
K_{i,3} := [s_{i+} - \delta, s_{i+}] \times \{ x_i \}
\]
\[
U_{i,3} := B_{\varepsilon}(v(K'_{i,3})) \Subset A_i
\]
\[
v(K_{i,2}) \Subset U_{i,2} \hspace{10pt} v(K_{i,3}) \Subset U_{i,3}
\]
Let $C'_i$ be the connected component of $\nu \setminus \overline{U_{i,1}}$ containing $z_i$ and  contained in $A_i$, and $\gamma'_i$ the closure of the connected component of $\nu - (\bigsqcup_j C'_j)$ contained in $V_i \cup V_{i-1}$ such that $z_i \notin \gamma'_i$. 
Then $\gamma'_i$ are closed arcs. 
There are small numbers $T \in (0, \delta/3)$, $T_i = \sigma_i T \in \{-T, T\}$, $\varepsilon'>0$ such that the following compact subsets and open subsets: 
\[
K_{i,4} := [-3T, 3T] \times \gamma'_i
\]
\[
U_{i,4} := B_{\varepsilon'}(v(K_{i,4})) \Subset A_v - \overline{U_{i,2} \sqcup U_{i,3}}
\]
\[
K_{i,5} := \{ 3T_i \} \times \gamma'_i 
\]
\[
U_{i,5} := B_{\varepsilon'}(v(K_{i,5})) \Subset A_v - B_{\varepsilon'}(v(\{0,T_i, 2T_i \}, \gamma'_i))
\]
\[
K_{i,6} := \{ 2T_i \} \times \gamma'_i 
\]
\[
U_{i,6} := B_{\varepsilon'}(v(K_{i,6})) \Subset A_v - B_{\varepsilon'}(v(\{0,T_i, 3T_i \}, \gamma'_i))
\]
\[
K_{i,7} := \{ T_i \} \times \gamma'_i 
\]
\[
U_{i,7} := B_{\varepsilon'}(v(K_{i,7})) \Subset A_v - B_{\varepsilon'}(v(\{0,2T_i, 3T_i \}, \gamma'_i))
\]
\[
K_{i,8} := \sigma_i [0,T] \times \gamma'_i 
\]
\[
U_{i,8} := B_{\varepsilon'}(v(K_{i,8})) = B_{\varepsilon'}(v(\sigma_i [0,T],\gamma'_i) ) \Subset A_v - B_{\varepsilon'}(U_{i,5})
\]

Taking $\varepsilon'$ small if necessary, we may assume that the subsets $U_{i,4}$ are open disks and $\overline{U_{i,4}} \cap \mathop{\mathrm{Sing}}v = \emptyset$. 
Define an open $C^0$-\nbd $\mathcal{V}$ of $v$ by 
\[
\mathcal{V} := \bigcap_{i=1}^{k} \bigcap_{l=1}^7 \mathcal{U}(K_{i,j}, U_{i,j})
\]
where $\mathcal{U}(K, U) := \{ w : \text{flow on } S \mid w(K) \subset U \}$ is an open $C^0$-subset.
Lemma~\ref{lem:inv_index02} implies that there is an open $C^0$-\nbd $\mathcal{U} \subseteq \mathcal{V}$ such that $\mathop{\mathrm{Sing}}(w) \cap \overline{U_{i,4}}  = \emptyset$ for any $w \in \mathcal{U}$. 

Fix any flow $w \in \mathcal{U}$ and any $i \in \{ 1, \ldots , k \}$. 
Put $C^w_i := w( [s_{i-},s_{i+}], x_i )$. 
Then $C^w_i \subset U_{i,1} = B_{\varepsilon}(v(K_{i,1}))$ and $w(\sigma_i [0,T], \gamma'_i) \cap B_{\varepsilon'}(U_{i,5}) = \emptyset$. 
Moreover, the closed arces $\gamma'_i$, $w(T_i,\gamma'_i)$, $w(2T_i,\gamma'_i)$, and $w(3T_i,\gamma'_i)$ are pairwise disjoint. 
Let $\gamma_i$ be the connected component of $\gamma'_i \setminus \bigcup_{i'} C^w_{i'} \subset \nu$ connecting $C^w_i$ and $C^w_{i-1}$. 
Since $\gamma_i$ is also a connected component of $\overline{\gamma'_i} \setminus \bigcup_{i'} C^w_{i'}$, the boundary $\partial \gamma_i \subset C^w_i \sqcup C^w_{i-1}$ consists of a point $p_{i-} \in C^w_{i-1} \cap \partial \gamma_i$ and a point $p_{i+} \in C^w_i \cap \partial \gamma_i$. 
From $\gamma_i \subset \gamma'_i$, we have $w([-3T, 3T], \gamma_i) \cap \overline{U_{i,2} \sqcup U_{i,3}} = \emptyset$.

\begin{claim}\label{claim:08}
$w([-3T, 3T], \gamma_i) \cap (C^w_i \sqcup C^w_{i-1}) = \emptyset$. 
\end{claim}

\begin{proof}[Proof of Claim~\ref{claim:08}]
Assume that there are a point $p \in \gamma_i$ and a number $t_p \in [-3T, 3T]$ with $
w(t_p,p) \in C^w_i$. 
From $w([-3T, 3T], \gamma_i) \cap \overline{U_{i,2} \sqcup U_{i,3}} = \emptyset$ and 
\[
w( [s_{i-}, s_{i-}+\delta] \sqcup [s_{i+}- \delta, s_{i+}], \{ x_i \}) \subset U_{i,2} \sqcup U_{i,3},
\]
we obtain that $w([-3T, 3T], p) \cap w( [s_{i-}, s_{i-}+\delta] \sqcup [s_{i+}- \delta, s_{i+}], x_i ) = \emptyset$. 
By $t_p \in [-3T, 3T] \subset (-\delta, \delta)$ and $C^w_i - w( [s_{i-}, s_{i-}+\delta] \sqcup [s_{i+}- \delta, s_{i+}], x_i)  = w( (s_{i-}+ \delta,s_{i+} - \delta), x_i)$, we have that $w(t_p,p) \in w( (s_{i-}+ \delta,s_{i+} - \delta), x_i)$ and so that $p \in w( (s_{i-}+ \delta - t_p,s_{i+} - \delta - t_p), x_i) \subseteq w( (s_{i-},s_{i+}), x_i) \subset C^w_i$, 
which contradicts that $\gamma_i \cap C^w_i = \emptyset$. 
Thus, $w([-3T, 3T], \gamma_i) \cap C^w_i = \emptyset$. 

By symmetry, the same argument implies that $w([-3T, 3T], \gamma_i) \cap C^w_{i-1} = \emptyset$. 
\end{proof}

Denote by $I_i$ the interval between $0$ and $3T_i$. 
The union $\Gamma_i := w(I_i,p_{i-}) \sqcup w(I_i,p_{i+}) \sqcup \overline{\gamma_i} \sqcup \overline{w(3T_i,\gamma_i)}$ is a loop in the open disk $U_{i,4}$. 
By the Jordan curve theorem, there is a closed disk $B_i \Subset U_{i,4}$ whose boundary is $\Gamma_i$. 
Since $\mathop{\mathrm{Sing}}(w) \cap \overline{U_{i,4}}  = \emptyset$, we have $\mathop{\mathrm{Sing}}(w) \cap B_i = \emptyset$.  
From the flow box theorem, there are trivial flow boxes of $w(I_i,p_{i-})$ and $w(I_i,p_{i+})$. 
This implies that the open arcs $w(T_i,\gamma_i)$ and $w(2T_i,\gamma_i)$ intersects $B_i$. 
Since closed arces $\overline{\gamma_i}$, $w(T_i,\overline{\gamma_i})$, $w(2T_i,\overline{\gamma_i})$, and $w(3T_i,\overline{\gamma_i})$ are pairwise disjoint, Claim~\ref{claim:08} implies that $\overline{\gamma_i} \sqcup w(T_i,\overline{\gamma_i}) \sqcup w(2T_i,\overline{\gamma_i}) \sqcup w(3T_i,\overline{\gamma_i}) \subset B_i$. 
Applying Lemma~\ref{lem:existence_transversal} to the closed disk $B_i$ with respect to $w$, there is a closed transverse arc $\mu_i \subset w((-3T,3T), \overline{\gamma_i})$ between $w((-3T,3T), p_{i+})$ and $w((-3T,3T), p_{i-})$ of $w$ such that the intersection $\mu_i \cap \partial B$ consists of two points $p_{i-}$ and $p_{i+}$.
Then there are numbers $s'^w_{i-}<0$ and  $s'^w_{i+}>0$ such that $w(\{ s'^w_{i-}, s'^w_{i+} \}, x_i) = \{ p_{i-}, p_{i+} \} = C^w_i \cap (\partial \gamma_i \sqcup \partial \gamma_{i-1})$. 
Put $C'^w_i := w( [s'^w_{i-}, s'^w_{i+}], \{ x_i \}) \subset C^w_i$. 
Then $\mu := \bigcup_{i=1}^{k} C'^w_i \cup \mu_i$ is a loop which consists of orbit arcs $C'^w_i$ and closed transverse arcs $\mu_i$. 
By arbitrarily small perturbation to $\mu$ near tangencies as in Figure~\ref{fig:perturbations_tang}, 
\begin{figure}
\begin{center}
\includegraphics[scale=0.6]{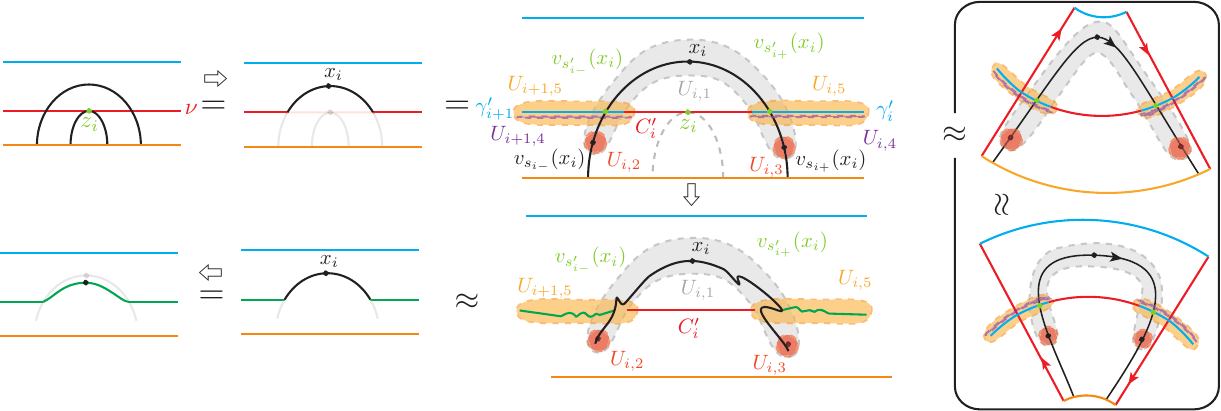}
\end{center}
\caption{Small perturbations of a transverse arc except for one tangency of $v$ into one of $w$.}
\label{fig:perturbations_tang}
\end{figure} 
there are a loop $\mu'$ which is transverse except finitely many tangencies with exactly the same number of inner tangencies and outer tangencies with respect to $w$ in the same order as $\nu$, and there is a closed annulus $A \Subset A_v$ which is a \nbd of the loop $\mu'$ such that the loop $\nu'$ is parallel to $\nu$ in $A_v$ and that the restriction $w \vert _{A}$ is locally topologically equivalent to the restriction $v \vert _{A_v}$. 
This implies assertion (1). 

Suppose that the loop $\nu$ bounds a closed disk $B$ such that $ \vert B \cap \mathop{\mathrm{Sing}}(w) \vert  < \infty$. 
By Lemma~\ref{lem:equiv_index}, the index $\mathrm{ind}_v(x)$ is the number of inner tangencies of the loop $\nu$ minus the number of outer tangencies of the loop $\nu$.  
Lemma~\ref{lem:eq_index} implies $\mathrm{ind}_v(x) = (2 + n_{i} - n_{o})/2 =  \sum_{y \in B \cap \mathop{\mathrm{Sing}}(w)} \mathrm{ind}_w(y)$. 
%
\end{proof}

\subsubsection{Proofs of Lemma~\ref{lem:inv_index} and Corollary~\ref{cor:inv_index_02}}

We show Lemma~\ref{lem:inv_index} and Corollary~\ref{cor:inv_index_02} as fllows. 

\begin{proof}[Proof of Lemma~\ref{lem:inv_index}]\label{prf:lem_inv_ind}
Let $v$ be a flow on a surface $S$ and $x \in S$ an isolated singular point. 
Fix a loop $C$ which is transverse at all but finitely many points and bounds an open disk $D$ whose closure is an isolated \nbd of $x$, and a small positive number $\varepsilon >0$. 
Lemma~\ref{lem:inv_multi-saddle02} implies that there is a desired $C^0$-\nbd $\mathcal{U}$ of $v$. 
\end{proof}

\begin{proof}[Proof of Corollary~\ref{cor:inv_index_02}]\label{prf:cor_inv_ind}
Since the resulting flows $v'_{\mathrm{dbl}}$ of any perturbations $v'$ of $v$ by taking the double $S_{\mathrm{dbl}}$ of $S$ are perturbations of the resulting $v_{\mathrm{dbl}}$ flow on $S_{\mathrm{dbl}}$ of $v$, we may assume that $S$ is closed.  
Then $\partial D$ is a closed transversal with inward (resp. outward) flow direction.
By the compactness of $\partial D$, the existence of trivial flow boxes implies that there is an open \nbd $U$ of $\overline{D}$ which is an open disk and is an isolated \nbd of $x$.  
By Lemma~\ref{lem:inv_index02}, there is a $C^0$-\nbd $\mathcal{U}$ of $v$ as in Lemma~\ref{lem:inv_index02}. 
Fix flow $w \in \mathcal{U}$. 

Suppose that $D \cap \mathop{\mathrm{Sing}}(w)$ contains exactly one singular point $x_w$. 
By \cite[Theorem~A]{kibkalo2022topological}, the isolated singular point $x_w$ consists of finitely many hyperbolic and parabolic sectors.
Then one minus half of the number of hyperbolic sectors of $x_w$ is the index of $x_w$. 
On the other hand, Lemma~\ref{lem:inv_index} implies that $1 = \mathrm{ind}_v(x) = \mathrm{ind}_w(x_w)$. 
Therefore, the isolated singular point $x_w$ consists of exactly one parabolic sector and so is either attracting or repelling. 
By Lemma~\ref{lem:inv_multi-saddle_01}, the attracting {\rm(resp.} repelling{\rm)} property is preserved and so the isolated singular point $x_w$ is attracting {\rm(resp.} repelling{\rm)}. 
\end{proof}

\bibliographystyle{abbrv}
\bibliography{yt20211011}

\begin{thebibliography}{10}

\bibitem{andronov1937rough}
A.~Andronov and L.~Pontryagin.
\newblock Rough systems.
\newblock In {\em Dokl. Akad. Nauk SSSR}, volume~14, pages 247--250, 1937.

\bibitem{aref1998stagnation}
H.~Aref and M.~Br{\o}ns.
\newblock On stagnation points and streamline topology in vortex flows.
\newblock {\em Journal of Fluid Mechanics}, 370:1--27, 1998.

\bibitem{barmak2008one}
J.~Barmak and E.~Minian.
\newblock {One-point reductions of finite spaces, h--regular CW--complexes and
  collapsibility}.
\newblock {\em Algebraic \& Geometric Topology}, 8(3):1763--1780, 2008.

\bibitem{barmak2008simple}
J.~A. Barmak and E.~G. Minian.
\newblock Simple homotopy types and finite spaces.
\newblock {\em Advances in Mathematics}, 218(1):87--104, 2008.

\bibitem{bolsinov2004integrable}
A.~V. Bolsinov and A.~T. Fomenko.
\newblock {\em Integrable Hamiltonian systems: geometry, topology,
  classification}.
\newblock CRC press, 2004.

\bibitem{gutierrez1978structural}
C.~Guti{\'e}rrez.
\newblock Structural stability for flows on the torus with a cross-cap.
\newblock {\em Transactions of the American Mathematical Society},
  241:311--320, 1978.

\bibitem{izydorek1996note}
M.~Izydorek, S.~Rybicki, and Z.~Szafraniec.
\newblock A note on the poincar{\'e}-bendixson index theorem.
\newblock {\em Kodai Mathematical Journal}, 19(2):145--156, 1996.

\bibitem{kibkalo2022topological}
V.~Kibkalo and T.~Yokoyama.
\newblock Topological characterizations of morse-smale flows on surfaces and
  generic non-morse-smale flows.
\newblock {\em Discrete and Continuous Dynamical Systems}, 42(10):4787--4822,
  2022.

\bibitem{kidambi2000streamline}
R.~Kidambi and P.~K. Newton.
\newblock Streamline topologies for integrable vortex motion on a sphere.
\newblock {\em Physica D: Nonlinear Phenomena}, 140(1-2):95--125, 2000.

\bibitem{labarca1990stability}
R.~Labarca and M.~Pacifico.
\newblock Stability of morse-smale vector fields on manifolds with boundary.
\newblock {\em Topology}, 29(1):57--81, 1990.

\bibitem{ma2005geometric}
T.~Ma and S.~Wang.
\newblock {\em Geometric theory of incompressible flows with applications to
  fluid dynamics}.
\newblock Number 119. American Mathematical Soc., 2005.

\bibitem{markley2023flows}
N.~G. Markley and M.~Vanderschoot.
\newblock {\em Flows on Compact Surfaces: The Weil--Hedlund--Anosov Program}.
\newblock Springer Nature, 2023.

\bibitem{may2003finite}
J.~May.
\newblock Finite topological spaces.
\newblock {\em Notes for REU}, 2003.

\bibitem{mccord1966singular}
M.~C. McCord.
\newblock Singular homology groups and homotopy groups of finite topological
  spaces.
\newblock {\em Duke Mathematical Journal}, 33(3):465--474, 1966.

\bibitem{moffatt2001topology}
H.~Moffatt.
\newblock The topology of scalar fields in 2d and 3d turbulence.
\newblock In {\em IUTAM symposium on geometry and statistics of turbulence},
  pages 13--22. Springer, 2001.

\bibitem{nikolaev1999flows}
I.~Nikolaev, E.~Zhuzhoma, and E.~V. {\v{Z}}u{\v{z}}oma.
\newblock {\em Flows on 2-dimensional manifolds: an overview}.
\newblock Number 1705. Springer Science \& Business Media, 1999.

\bibitem{peixoto1962structural}
M.~M. Peixoto.
\newblock Structural stability on two-dimensional manifolds.
\newblock {\em Topology}, 1(2):101--120, 1962.

\bibitem{sakajo2014unique}
T.~Sakajo, Y.~Sawamura, and T.~Yokoyama.
\newblock Unique encoding for streamline topologies of incompressible and
  inviscid flows in multiply connected domains.
\newblock {\em Fluid Dynamics Research}, 46(3):031411, 2014.

\bibitem{sakajo2015transitions}
T.~Sakajo and T.~Yokoyama.
\newblock Transitions between streamline topologies of structurally stable
  hamiltonian flows in multiply connected domains.
\newblock {\em Physica D: Nonlinear Phenomena}, 307:22--41, 2015.

\bibitem{sakajo2018tree}
T.~Sakajo and T.~Yokoyama.
\newblock {Tree representation of topological streamline patterns of
  structurally stable 2D Hamiltonian vector fields in multiply conected
  domains}.
\newblock {\em The IMA Journal of Applied Mathematics}, 83:380--411, 2018.

\bibitem{sakajo2020discrete}
T.~Sakajo and T.~Yokoyama.
\newblock Discrete representations of orbit structures of flows for topological
  data analysis.
\newblock {\em Discrete Math. Algorithms Appl.}, 15(6):Paper No. 2250143, 38,
  2023.

\bibitem{smale1961gradient}
S.~Smale.
\newblock On gradient dynamical systems.
\newblock {\em Annals of Mathematics}, pages 199--206, 1961.

\bibitem{stong1966finite}
R.~E. Stong.
\newblock Finite topological spaces.
\newblock {\em Transactions of the American Mathematical Society},
  123(2):325--340, 1966.

\bibitem{yokoyama2017decompositions}
T.~Yokoyama.
\newblock Decompositions of surface flows.
\newblock {\em arXiv preprint arXiv:1703.05501}, 2017.

\bibitem{yokoyama2023dependency}
T.~Yokoyama.
\newblock Dependency of the positive and negative long-time behaviors of flows
  on surfaces.
\newblock {\em Journal of Differential Equations}, 368:376--393, 2023.

\bibitem{yokoyama2021combinatorial}
T.~Yokoyama.
\newblock Combinatorial structures of the space of hamiltonian vector fields on
  compact surfaces.
\newblock {\em Topology and its Applications}, page 109439, 2025.

\bibitem{yokoyama2013word}
T.~Yokoyama and T.~Sakajo.
\newblock Word representation of streamline topologies for structurally stable
  vortex flows in multiply connected domains.
\newblock {\em Proceedings of the Royal Society A: Mathematical, Physical and
  Engineering Sciences}, 469(2150):20120558, 2013.

\bibitem{yokoyama2021complete}
T.~Yokoyama and T.~Yokoyama.
\newblock {Complete transition diagrams of generic Hamiltonian flows with a few
  heteroclinic orbits}.
\newblock {\em Discrete Mathematics, Algorithms and Applications},
  13(02):2150023, 2021.

\bibitem{yokoyama2021cot}
T.~Yokoyama and T.~Yokoyama.
\newblock {COT representations of 2D Hamiltonian flows and their computable
  applications}, {\it preprint}, 2021.

\end{thebibliography}

\end{document}